\newif\ifTR

\TRtrue

\newif\ifTS
\TSfalse

\ifTR%
\documentclass[a4paper,11pt]{article}
\usepackage{authblk}
\usepackage{natbib}
\usepackage[cm]{fullpage}
\usepackage{caption}
\usepackage{subcaption}
\usepackage{amsthm}

\newtheorem{proposition}{Proposition}
\newtheorem{theorem}{Theorem}

\newtheorem{corollary}{Corollary}
\newtheorem{definition}{Definition}
\else

\OneAndAHalfSpacedXI
\usepackage{natbib}
\bibpunct[, ]{(}{)}{,}{a}{}{,}%
\TheoremsNumberedThrough
\ECRepeatTheorems
\EquationsNumberedThrough

\fi
\usepackage{amsmath,amsfonts,amssymb}
\usepackage{graphicx}
\usepackage{paralist}
\usepackage[hidelinks]{hyperref}
\usepackage{comment}
\usepackage{bm}
\usepackage{graphics}
\usepackage[capitalize,noabbrev]{cleveref}
\usepackage{multirow,longtable}
\usepackage{rotating}
\usepackage{adjustbox}
\usepackage{algorithm}
\usepackage{algpseudocode}
\usepackage{soul}
\usepackage{xspace}
\usepackage{booktabs}
\usepackage{pdflscape}
\usepackage{tikz}
\usepackage{float}
\usepackage[normalem]{ulem}
\usepackage{nicefrac}
\usepackage{color} 
\usepackage{pgfplots}
\usepackage{booktabs}
\usepackage{pifont}
\usepackage{eurosym}
\usepackage{makecell}

\pgfplotsset{compat = newest}
\usepackage{booktabs} 
\usepackage{subcaption}
\usetikzlibrary{arrows, positioning, fit, calc, backgrounds, quotes, angles}
\usetikzlibrary{arrows.meta, positioning}

\allowdisplaybreaks

\newcommand{\blue}[1]{{\color{black}#1}}

\newcommand{\red}[1]{{\color{red}#1}}

\newcommand{\sw}[1]{{\color[rgb]{0,0.5,0}#1}}
\newcommand\crowdshipper{occasional driver\xspace}
\newcommand\crowdshippers{occasional drivers\xspace}

\newcommand\Crowdshippers{Occasional drivers\xspace}
\newcommand{\apdx}{electronic companion\xspace}

\newcommand{\var}{z_{wb}}

\newcommand{\repTH}[3]{{
		
		\ifTS
		
		\begin{repeattheorem}[#1]
			
			#2
			
		\end{repeattheorem}
		
		\begin{proof}{Proof.}
			
			#3 \Halmos
			
		\end{proof}
		
		\else
		
		\begin{proof}[\textbf{Proof of #1}]
			
			#3
			
		\end{proof}
		
		\fi
		
}}

\ifTR

\DeclareMathOperator*{\argmax}{\mathrm{arg\,max}}
\fi

\newcommand{\myabstract} 
{
Challenges in last-mile delivery such as high costs, rising customer expectations, and congested urban traffic have encouraged innovative solutions like crowdsourced delivery, in which occasional drivers undertake delivery tasks along their pre-planned trips in exchange for compensation. 
A key challenge is that occasional drivers’ acceptance behavior towards offered tasks is uncertain and influenced by task properties and compensation. The current literature lacks formulations that address this challenge \blue{with a joint and exact approach}. Hence, we formulate an integrated optimization problem that maximizes total expected cost savings by offering bundles of tasks to occasional drivers. To this end, we simultaneously determine the optimal set of bundles their assignment to occasional drivers, and personalized compensations for each bundle-driver pair while considering bundle- and compensation-dependent acceptance probabilities, which are captured via generic logistic functions. The vast number of potential bundles, combined with incorporating acceptance probabilities leads to a mixed-integer nonlinear programming (MINLP) formulation with exponentially many variables. Using mild assumptions, we address these complexities by exploiting properties of the problem, leading to an exact linearization of the MINLP which we solve via a tailored exact column generation algorithm. Our algorithm uses a variant of the elementary shortest path problem with resource constraints (ESPPRC) that features a non-linear and non-additive objective function as its subproblem, for which we develop tailored dominance and pruning strategies. We introduce several heuristic and exact algorithm variants, and perform an extensive set of computational experiments evaluating the performance of algorithms and the structure of the solutions. The results demonstrate the efficiency of exact and heuristic algorithm variants for instances with up to 120 tasks and 60 drivers, and highlight the advantages of integrated decision-making over sequential approaches in line with the recent literature. The sensitivity analysis indicates that sensitivity to compensation is the most influential factor in shaping the bundle structure.
}

\title{Pricing, bundling, and driver behavior in crowdsourced delivery}	

\ifTR

\author{Alim Bu\u{g}ra \c{C}{\i}nar\textsuperscript{1}, Claudia Archetti\textsuperscript{2}, Wout Dullaert\textsuperscript{1}, Markus Leitner\textsuperscript{1}, and Stefan Waldherr\textsuperscript{1}} 
\affil{\textsuperscript{1}Department of Operations Analytics, Vrije Universiteit Amsterdam, The Netherlands \\ 
\textsuperscript{2}Department of Economics and Management, University of Brescia, Italy}

\else

\fi

\begin{document}	

\ifTR
\maketitle
\begin{abstract}
	\myabstract 
	\textbf{Crowdsourced delivery, occasional drivers, compensation schemes, acceptance uncertainty, \crowdshipper behavior, logistic regression} 
\end{abstract}

\else




\RUNAUTHOR{\c{C}{\i}nar et al.}
%
\RUNTITLE{Pricing, bundling, and driver behavior in crowdsourced delivery}
%
\TITLE{Pricing, bundling, and driver behavior in crowdsourced delivery}
%
\ARTICLEAUTHORS{%
	\AUTHOR{Alim Bu\u{g}ra \c{C}{\i}nar\textsuperscript{1}, Claudia Archetti\textsuperscript{2}, Wout Dullaert\textsuperscript{1}, Markus Leitner\textsuperscript{1}, Stefan Waldherr\textsuperscript{1}} 
	\AFF{1. Department of Operations Analytics, Vrije Universiteit Amsterdam, The Netherlands \\ 
	2. Department of Economics and Management, University of Brescia, Italy}	
	\EMAIL{a.b.cinar@vu.nl}
	\EMAIL{claudia.archetti@unibs.it}	
	\EMAIL{w.e.h.dullaert@vu.nl}
	\EMAIL{m.leitner@vu.nl}
	\EMAIL{s.m.g.waldherr@vu.nl}	
}
\ABSTRACT{%
	\myabstract
}%
%


\pagenumbering{arabic}
\maketitle
%
%
%
%

\fi

\newcommand{\compensationvalues}{%
	There exists an optimal solution to formulation~\eqref{eq:minlp} in which each (optimal) compensation value $C^*_{wb}$ for every bundle $b$ offered to \crowdshipper $w$ is calculated as
	\begin{equation*} 
		C^*_{wb} = \argmax_{C \ge 0} ~ P_w(b,C)\left(\sum_{i \in b}c_i-C\right).
	\end{equation*}%
}

\newcommand{\optimalcompensation}{
	For the logistic acceptance probability function in \cref{def:logistic_prob}, the optimal compensation value offered to \crowdshipper $w\in N$ for the bundle $b\in B$ is equal to 

	\begin{equation*} 
	C^*_{wb} = \begin{cases}
		- \frac{W(e^{\alpha^w + B^w X_{wb} + D^w Y_{w} + \gamma^w\bar{C} - 1})-\gamma^w \bar{C} + 1}{\gamma^w}, & \text{if } \gamma_w\bar C > 1+e^{\alpha^w + B^w X_{wb} + D^w Y_{w}},\\[2mm]
		0 & \text{otherwise},
	\end{cases}
	\end{equation*}
	in which $W(\cdot)$ is the Lambert $W$ function, and $\bar{C}=\sum_{i \in b}c_i$. \blue{The case $C^*_{wb}=0$ occurs when the saving relative to using the 3PL is too small to justify offering any compensation.}
	
}

\newcommand{\reducedcost}{
Assume that the acceptance behavior of \crowdshipper $w
\in N$ follows the logistic acceptance probability function~\eqref{eq:acceptance_probability}. Then, the reduced costs~\eqref{eq:sp:obj} of a bundle $b\in B$ offered to \crowdshipper $w\in N$ are equal to 
	\begin{equation*}
		\tilde{c}_{wb}=\frac{W(e^{\alpha^w + B^w X_{wb} + D^w Y_{w} + \gamma^w\sum_{i\in b} c_i - 1})}{\gamma^w} - \sum_{i \in b}\pi_i - \mu_w.  
	\end{equation*} 
}

\newcommand{\dominance}{
	Let $L_w^p$ and $L_w^f$, $L_w^{p} \neq L_w^{f}$ be two labels associated with driver $w\in N$ whose bundles end at the same task, i.e., $b^p=(i^p_1,\dots,i^p_\ell)$ and $b^f=(i^f_1, \dots, i^f_{\ell'})$ such that $i^p_\ell=i^f_{\ell'}$.  
Label $L^p$ dominates $L^f$ if the following conditions hold:	 
\begin{enumerate}[(i)]
	\item $B^wX^p +\gamma^w\bar{C}^p \geq B^w X^f+\gamma^w\bar{C}^f$
	\item $\sum_{i \in {b^p}}\pi_i \leq \sum_{i \in {b^f}}\pi_i$
	\item $R^p \supseteq R^f$
	\item $q^p\leq q^f$
\end{enumerate}
}

\newcommand{\rcbasedpruning}{
	Let $L_w^p$ be an arbitrary label associated to driver $w\in N$ and bundle $b^p=(i_1, \dots, i_\ell)$.
Let 
$c^\mathrm{max}=\max_{t \in R^p}c_t$,
$c^\mathrm{min}=\min_{t \in R^p}c_t$,
$\pi^\mathrm{max}=\max_{t \in R^p}\pi_t$ and 
$\pi^\mathrm{min}=\min_{t \in R^p}\pi_t$.
The	following conditions hold for any label $L_w^f$ extending $L_w^p$ by exactly $k$ items:	
\begin{enumerate}[(i)]
	\item $X^p_i+ k U^\mathrm{max}_i \ge X^f_i \ge X^p_i+ k U^\mathrm{min}_i$ for all $i\in \{1, \dots n \}$
	\item $\bar{C}^p + k c^\mathrm{max} \ge  \bar{C}^f \ge \bar{C}^p + k c^\mathrm{min}$		
	\item $\pi^p+k\pi^\mathrm{max} \ge \pi^f \ge \pi^p+k\pi^\mathrm{min}$
\end{enumerate}
}

\newcommand{\rcbasedpruningcorr}{
	Let $L_w^p$ be an arbitrary label associated to driver $w\in N$ and bundle $b^p=(i_1, \dots, i_\ell)$. Consider $\hat{X}^k=(\hat{X}^k_1, \dots, \hat{X}^k_n)\in \mathbb{R}^n_+$, $\hat{C}^k\ge 0$, and $\hat{\pi}^k\ge 0$ such that 
\begin{enumerate}[(i)]
	\item $\hat{X}_i^k=X^p_i+ k U^\mathrm{max}_i$ for all $i\in \{1, \dots n : B_i^w > 0\}$
	\item $\hat{X}_i^k=X^p_i+ k  U^\mathrm{min}_i$ for all $i\in \{1, \dots n : B_i^w <  0\}$
	\item $\hat{C}^k=\bar{C}^p+k c^\mathrm{max}$			
	\item $\hat{\pi}^k=\pi^p+k\pi^\mathrm{min}$
\end{enumerate}

The reduced costs $\tilde{c}^f$ of any label $L_w^f$, obtained by extending $L_w^p$ by exactly $k$ items, are smaller than or equal to $\hat{c}^k = \frac{W(e^{\alpha^w + B^w \hat{X}^k + D^w Y_w + \gamma^w\hat{C}^k - 1})}{\gamma^w} - \hat{\pi}^k - \mu_w$.	
}

\newcommand{\conditionforpruning}{
	Let $L_w^p$ be an arbitrary label associated to driver $w\in N$ and bundle $b^p=(i_1, \dots, i_\ell)$ to which at most $k_\mathrm{max}$ tasks can be added. If the values $\hat{X}^k=(\hat{X}^k_1, \dots, \hat{X}^k_n)\in \mathbb{R}^n_+$, $\hat{C}^k\ge 0$ and $\hat{\pi}^k\ge 0$ satisfying the conditions of \cref{corr:rcbasedpruning} imply that $\hat{c}^k \le 0$ for every $k\in \{1, \dots k_\mathrm{max}\}$, then no extension of $L_w^p$ can correspond to a variable with positive reduced cost.
}

\newcommand{\dominancesimple}{
	Let $L^p_w$ and $L_w^f$, $L_w^p\ne L_w^f$, be two labels associated with driver $w\in N$ that visit the same set of tasks. Label $L^p_w$ dominates $L_w^f$ if $\tilde{c}^p > \tilde{c}^f$ and $R^p\supseteq R^f$. 
	}

\newcommand{\dominancecorr}{
	The logistic acceptance probability function introduced in \cref{eq:acceptance_model} satisfies the conditions for dominance introduced in \cref{prop:dominancerules}. 
}

\newcommand{\detourupperbound}{
	Let $L_w^p$ be a label associated with driver $w\in N$ and bundle $b^p=(i_1, \dots, i_\ell)$ to which at most $k_\mathrm{max}$ tasks can be added. Let $L_w^k$ be any extension of $L_w^p$ with $k \in \{1, \dots k_\mathrm{max}\}$ additional tasks with total detour $\Delta^k$. 
	Consider $\hat{\pi}^k$, $\hat{C}^k$, and $\hat{b}^k$ as defined in \cref{corr:rcbasedpruning}. Let 
	$\bar{\Delta}^k=\frac{\ln(\gamma^w(\hat{\pi}^k + \mu_w )) + \gamma^w(\hat{\pi}^k + \mu_w - \hat{C}^k) - \alpha^w - \beta_2^w (\hat{b}^k) + 1 }{\beta_1^w}$. If $\Delta^k \ge \bar{\Delta}^k$, then $\tilde{c}^k\le 0$. 
	
%
}

\newcommand{\detourbasedpruning}{
	Consider a label $L_w^p$ associated with driver $w\in N$ and bundle $b^p=(i_1, \dots, i_\ell)$ whose detour is equal to $\Delta^p$ and to which at most $k_\mathrm{max}$ tasks can be added. Let $t\in R^p$ and $\Delta^t=\Delta^p + \mathcal{U}_{\Delta}(\mathcal{T}^p,t)$ be the total detour after extending label $p$ to task $t$. Task $t$ can be removed from $R^p$ if $\Delta^t \ge \max_{k \in \{1,\dots,k_\mathrm{max}\}}\bar{\Delta}^k$.
	
}

\section{Introduction} \label{sec:intro}

Rising demand and growing customer expectations in last-mile delivery, coupled with congested urban traffic, increasing emissions, and environmental concerns, make innovative solutions more critical than ever. A novel strategy is crowdsourced delivery in which tasks are outsourced to occasional drivers, i.e., independent individuals that use their residual capacity by integrating delivery tasks into their pre-planned trips in exchange for compensation \citep{archetti_vehicle_2016, mancini_bundle_2022, wang_joint_2023, cinar2024role}.
Thereby, the term occasional emphasizes sporadic, opportunistic participation rather than regular, shift-based delivery work.
Involving occasional drivers in delivery services can contribute to reductions in traffic congestion and emissions, and offers several advantages over traditional methods including greater flexibility in delivery capacity and potential cost savings which are crucial to deal with growing customer expectations in increasingly congested urban networks \citep{kaspi_directions_2022}. 
In some cases, multiple deliveries are bundled into single offers to enhance driver utilization, provide higher earning opportunities, and reduce per-task delivery costs \citep{mancini2025dynamic}.
\blue{}

Crowdsourced delivery introduces unique planning and operational challenges that are not present in traditional delivery settings. These include 
the uncertainty in individuals' responses to delivery offers and the need to determine appropriate compensation 
\citep{savelsbergh2024challenges}.
Furthermore, operator decisions, such as task bundling and compensation, directly influence 
individuals’ responses \citep{mohri_modeling_2024}, highlighting the need for models that explicitly incorporate interdependencies among these decisions for improved \blue{decision making}. 

Still, existing studies tend to \blue{address} these challenges \blue{individually, partially or (at best) sequentially}, leaving a research gap for models that \blue{handle} these challenges \blue{jointly}~\citep{savelsbergh2024challenges}.
%
To address this gap, \blue{we study a static, single-decision-epoch problem in which} an operator receives delivery requests from customers, creates bundles of tasks, and offers \blue{at most one} bundle to \blue{each} occasional driver who may accept or reject these offers \blue{as a whole. Treating the bundle as the unit of choice enables the operator to exploit the cost subadditivity \citep{mancini_bundle_2022}}. We assign delivery tasks to drivers' pre-planned trips, making use of spare transportation capacity with limited additional effort. \blue{Such a static problem may represent one decision epoch within a periodic decision-making process where the operator collects the delivery tasks and the occasional drivers who have declared their availability, and designs all offers for this batch in a single offer round ahead of the delivery. Customer demand and driver participation are therefore exogenous inputs.}

\blue{Our setting} contrasts with settings where drivers perform dedicated trips only when delivery tasks are available, such as in-house delivery services or gig economy platforms like UberEats, Instacart, and DoorDash. Although these platforms also fall under the umbrella of crowdsourced delivery, many participating drivers, while independent in determining availability and accepting tasks, typically commit themselves to multi-hour shifts and undertake dedicated trips to complete deliveries \blue{and compete with each other to collect higher compensations}. \blue{Consistent with this distinction, we exclude from our scope offers involving multiple pickup locations, preparation and depot-waiting times, strategic interactions among market participants commonly observed on the aforementioned gig economy platforms, and the associated routing, behavioral, and game-theoretic extensions.}


\blue{To capture drivers’ acceptance behavior in our setting, we use a logistic probability function that incorporates the additional detour required, the number of extra stops, and the compensation. We further illustrate how this function can be extended to include other driver-specific attributes and bundle-related factors reflecting how well an offer aligns with their originally planned trip.}
These probability functions are integrated 
 into a unified decision making problem that 
 simultaneously optimizes 
 task bundling, bundle assignment, and compensation decisions while maximizing  
 the total expected cost savings of an
  operator.


We model this decision making problem as a novel 
mixed-integer non-linear program (MINLP) with an exponential number of variables. 
We show how to linearize the MINLP and develop an exact solution algorithm based on column generation. Our computational experiments demonstrate that our algorithm can solve large instances to optimality. We further implement heuristic variants that 
achieve high-quality solutions in significantly reduced times. Moreover, our analysis reveals that drivers’ sensitivity to compensation is the most critical factor influencing bundle structure, with driver heterogeneity amplifying differences in offers made to different driver groups.


Our main contributions can be summarized as follows:
\begin{itemize}
	\item 
	\blue{We are the first to propose an exact approach that jointly integrate task bundling, bundle assignment, and compensation decisions for offers made to occasional drivers, incorporating their compensation- and bundle-dependent acceptance probabilities into a unified optimization framework.}
	Our model can handle generic 
	probability functions modeling the acceptance behavior of occasional drivers and is applicable to many variants 
	of crowdsourced delivery, such as 
		using in-store customers or 
		assigning pick-up and delivery tasks to \crowdshippers.
	\item We propose a mixed-integer nonlinear programming (MINLP) model with an exponential number of variables, show how to linearize it, and present an effective exact solution algorithm based on column generation. \blue{To this end, we develop a tailored pricing algorithm that generalizes the single-task compensation optimization of \citet{cinar2024role} to bundled offers.}
	\item We show that the pricing subproblem can be solved as a variant of the elementary shortest path problem with resource constraints (ESPPRC)
	featuring a non-linear, non-additive objective function. These properties prohibit the use of classical labeling approaches 
		as they rely 
		on dominance rules that assume additivity \citep{sadykov2021bucket, costa2019exact}. 
		Therefore, we introduce novel dominance and pruning rules that explicitly address the non-linearity and non-additivity of the objective function, leading to a computationally efficient algorithm that advances beyond previous approaches to similar problems \citep{rostami2021branch}.
	\item We develop several exact and heuristic variants of our solution algorithm to evaluate the effectiveness of its individual components and to address scalability in large instances.
	\item We evaluate all considered variants of the proposed solution algorithm through extensive computational experiments and analyze sensitivity across varying numbers of \crowdshippers, tasks, and behavioral class distributions. The results show that
	 (i) the newly introduced dominance and pruning rules are crucial to achieve substantial performance improvements in our algorithms, (ii) our exact algorithms are able to solve large instances within a reasonable amount of time, (iii) our heuristic variants provide high quality solutions in a short amount of time, and (iv) the integrated approach leads to significantly better solutions than a sequential method in line with the recent literature, where bundling, compensation, and assignment decisions are made in succession.
	\item \blue{We evaluate personalized compensation offer design from both the operator's and the drivers' perspectives. Our integrated personalized offer design yields the highest operator savings and, when a high monetary value is assigned to driver effort, a more balanced cross-class welfare distribution than approaches that partially or fully neglect driver heterogeneity.}
	
\end{itemize}


%

The paper is organized as follows. In \cref{sec:literature}, we review the related literature. \cref{sec:problem} formally presents the problem, followed by the solution methodology in \cref{sec:methodology}. \cref{sec:experiments} describes the experimental setting, while \cref{sec:results} presents the results of the computational study. Finally, \cref{sec:conclusion} concludes the paper and provides future research directions.

\section{Literature review} \label{sec:literature}

Major contributions of this paper include the simultaneous consideration of task bundling and compensation decisions while explicitly modeling the compensation- and bundle-dependent acceptance probabilities of \crowdshippers and the development of an exact solution algorithm that integrates these aspects. To the best of our knowledge, prior crowdsourced delivery literature has addressed these elements \blue{individually, partially or (when covering all of them) sequentially rather than in a joint and exact manner.}
Throughout this review, we consistently use the term “occasional driver” to refer to all types of crowdsourced drivers, regardless of the original terminology employed by their authors. A broader overview of the crowdsourced delivery literature is given in \citet{savelsbergh2024challenges}. 

Several studies consider settings in which the operator allocates 
individual delivery tasks to occasional drivers. 
\citet{gdowska_stochastic_2018, santini_probabilistic_2022} and \citet{ausseil_online_2024} 
consider acceptance probabilities but neglect the compensation decisions.
These studies also assume static acceptance probabilities within the allocation decisions, though \citet{ausseil_online_2024} introduce a dynamic mechanism to update probabilities based on past driver responses. A group of more closely related works (\citet{barbosa_data-driven_2023, hou_optimization_2022, cinar2024role}) explicitly model compensation-dependent acceptance probabilities together with compensation decisions. Specifically, \citet{barbosa_data-driven_2023} consider identical acceptance probabilities and compensations for all drivers, employing approximation methods for compensation optimization. \citet{hou_optimization_2022} optimize task assignment and compensation decisions sequentially in a heuristic manner, while \citet{cinar2024role} provide an exact approach integrating these decisions. In contrast to our model, these studies do, however, not consider bundling decisions.

A second set of studies 
focus primarily on bundling decisions while assuming 
that acceptance behavior of occasional drivers and compensation amounts are either fully known or deterministic. 
One common approach, as adopted by \citet{kafle_design_2017, mancini_bundle_2022, mancini_bundle_2024, mancini2025dynamic}, is to use
auction-based mechanisms where occasional drivers submit bids specifying both the bundles of tasks they are willing to deliver and the compensations they request. In these models, 
the operator selects 
bundles from a limited, predefined set
representing driver bids and associated compensations.
In contrast, 
\citet{wang_joint_2023} and \citet{torres_vehicle_2022} employ column generation based approaches and are, therefore, not restricted to a predefined bundle set.
As in our approach, this has
the advantage of systematically evaluating the complete bundle space to identify optimal task bundles.
 However, these studies differ fundamentally from our work by assuming 
 deterministic
 driver acceptance and fixed compensations
 calculated as linear functions of bundle characteristics, typically combining fixed service fees and variable costs dependent on delivery distance, time, or number of stops.
 Unlike all the aforementioned approaches, our research explicitly incorporates stochastic, compensation-dependent acceptance probabilities into the simultaneous optimization of bundling, compensation, and bundle allocation decisions.

Another stream of literature addresses pickup-and-delivery problems, where occasional drivers may sequentially serve multiple tasks. Among these studies, \citet{le_designing_2021} explicitly optimize compensation decisions but rely on a threshold-based acceptance model, assuming drivers automatically accept any compensation above a given threshold. 
In contrast, \citet{horner_optimizing_2021} consider drivers who may select multiple tasks from a menu offered by the operator, modeling acceptance probabilities individually per task rather than collectively. They treat both acceptance probabilities and compensations as fixed parameters and  
assume a linear dependency of compensation on 
service time and distance, and of acceptance probabilities on compensation. 
Lastly, \citet{behrendt_task_2024} explicitly model individual drivers' task acceptance probabilities using a multinomial logit model dependent on compensation and travel distance, and optimize compensation accordingly. However, they assume 
independence between tasks and, therefore, calculate individual acceptance probabilities rather than modeling drivers' acceptance of task bundles as a holistic choice. In contrast to these approaches, we
explicitly consider 
bundles as 
entities and 
captures the interdependencies between tasks of a bundle. 
By modeling acceptance probabilities at the individual driver-bundle level, our methodology significantly enhances realism in representing actual driver decision-making.

The studies most closely related to ours, as they incorporate compensation decisions and/or acceptance probability modeling within bundle allocation contexts, are those by \citet{torres_crowdshipping_2022, cerulli2024bilevel, macrina2024bundles}.
\citet{torres_crowdshipping_2022} employ a bundle (route)-centric approach
where the probability that a bundle is accepted by any occasional driver is positively correlated with 
the compensation offered. In their model, compensation is, however, not a decision variable. Instead, it 
is calculated as a linear function of fixed costs, distance, and the number of deliveries. 
In contrast, our approach treats 
compensation as a decision variable and explicitly models acceptance probabilities of individual occasional drivers towards specific bundles, enabling a more detailed representation of individual driver preferences and behavior.
\citet{cerulli2024bilevel} consider compensation decisions as a margin optimization problem
where
the compensation paid to occasional drivers is a margin of the price paid by the customers (which is a fixed parameter).
They consider the acceptance behavior of occasional drivers without explicitly incorporating their acceptance probability functions, thus obtaining a fully deterministic model. Instead, they
model the behavior of occasional drivers as the follower's problem that depends on the platform decisions for the bundle and the margin in a bi-level formulation. In contrast, our methodology considers compensation optimization decisions for each bundle and occasional driver pair individually, which increases the flexibility of the decision maker. Furthermore, we explicitly incorporate individual acceptance probability functions that are influenced by the platform's decisions. 
Finally, \citet{macrina2024bundles} solve the bundle generation problem using a greedy algorithm and the compensation calculation problem considering the acceptance probabilities in a sequential way. They consider an infinite supply of occasional drivers whose acceptance behaviors are identical.
In contrast, our methodology optimizes bundling and compensation decisions simultaneously and explicitly models a finite, heterogeneous group of occasional drivers. Our approach captures each driver's individual acceptance probabilities as functions of compensation, bundle attributes, and driver-specific characteristics, thus significantly enhancing both the realism and practical applicability.

In summary, the existing literature lacks studies that \blue{propose exact solution approaches to jointly} optimize task bundling, bundle assignment, and compensation decisions, while explicitly modeling acceptance probabilities at the individual driver-bundle level as functions of these interdependent decisions.
This is also highlighted by \cref{tab:literature} showing which of the studies reviewed in this section consider acceptance probability, compensation and bundling, see \cref{sec:literature_apdx}.
Our comprehensive approach aims to address this gap, significantly enhancing operational relevance and realism, and thereby advancing the state-of-the-art in crowdsourced delivery optimization.

\section{Problem definition and mathematical formulation} \label{sec:problem} 

This study examines a delivery optimization problem faced by an operator that fulfills delivery tasks using external transport services. 
The operator can outsource tasks to a third-party logistics (3PL) company or offer sets of tasks (i.e., bundles) to \crowdshippers that may be willing to perform the associated delivery tasks for a compensation proposed by the operator.
\Crowdshippers are autonomous agents who have the freedom to accept or decline offers made by the operator. When occasional drivers accept a bundle, they commit to delivering all tasks within the bundle and receive the agreed compensation.
Any task not accepted by occasional drivers is handled by the 3PL at a predetermined cost \blue{and does not return to the offer pool}. The operator's goal is to maximize the total expected cost savings in comparison to using 3PL by determining an optimal set of bundles and corresponding compensations offered to \crowdshippers.


We now formally describe the optimization problem considered in this article. Each problem instance consists of a 
set of tasks $M$, a set of (pick-up) depots $D$ and a set of \crowdshippers $N$. Task $i \in M$ is identified by its delivery location $t_i$, 
load consumption $q_i\ge 0$, and predetermined 3PL cost $c_i\ge 0$. Each 
\crowdshipper $w\in N$ is identified by its current location $s_w$, destination $e_w$, and 
capacity $Q^w$. Each depot $j\in D$ is associated with a location $d_j$ and a task set $T(d_j) \subseteq M$ that denotes the tasks that can be delivered from $d_j$.
The operator 
creates bundles of tasks and offers them to \crowdshippers. These bundles are represented as ordered tuples of tasks and the set of all possible bundles is denoted as $B=\{\left(i_1,i_2,\ldots, i_\ell\right) | \{i_1, i_2, \ldots, i_\ell\}\in 2^M \}$. 
We assume that the operator can estimate the probability $P_w(b, C)$ that \crowdshipper $w\in N$ accepts bundle $b \in B$ offered in exchange for the proposed \emph{personalized} compensation $C\ge 0$.

A feasible solution to the problem consists of a set of bundles $B_\mathrm{S} \subset B$ corresponding compensations $C(b)$ for each bundle $b \in B_\mathrm{S}$, and the \crowdshipper $w(b)$ to whom bundle $b$ is offered. At most one bundle can be offered to every \crowdshipper and each task can belong to at most one offered bundle. The total load of all tasks in a bundle $b\in B_\mathrm{S}$ must not exceed the capacity of the \crowdshipper to whom the bundle is offered, i.e., $\sum_{i\in b} q_i \le Q^{w(b)}$. Additionally, there must exist at least one depot from which all tasks in the bundle can be delivered, i.e. $\exists d\in D: b \subseteq T(d)$. The total expected cost savings measures the savings obtained by using occasional drivers instead of shipping every task individually through the 3PL, and is calculated as
\begin{equation} \label{eq:expected_savings}
	\sum_{b \in B_\mathrm{S}} P_{w(b)}(b, C(b)) \left( \sum_{i \in b} c_i - C(b) \right).
\end{equation}
The objective is to identify a feasible solution maximizing the total expected cost savings.

\paragraph{Assumptions:}

In the remainder of this article, we assume that an \crowdshipper never accepts a bundle without compensation.
%
We also assume 
that each \crowdshipper can visit every depot, and that there exists a mapping $B\times N\mapsto D$ that indicates an optimal, detour-minimizing, depot $d\in D$ when given a bundle $b\in B$ and an \crowdshipper $w\in N$.
We note that while these three assumptions are made to simplify notation, all results discussed in the following sections are either valid or could easily be generalized to situations in which they do not hold. 
We define the detour of a driver with respect to a bundle as the additional distance incurred from visiting the locations of all tasks in the bundle when traveling from the driver’s origin to the destination.

\begin{figure}
	\captionsetup[subfigure]{aboveskip=-0.5pt,belowskip=-0.5pt}
	\begin{subfigure}[]{0.7\textwidth}
		\centering
		\resizebox{0.6\textwidth}{!}{%
\begin{tikzpicture}[
	scale=0.5,
	->,
	>=stealth,
	node distance=1.5cm and 1.5cm,
	every node/.style={draw, circle, minimum size=0.75cm, inner sep=0pt},
	longdash/.style={dash pattern=on 10pt off 4pt}, 
	dotdash/.style={dash pattern=on 8pt off 3pt on 1pt off 3pt},
	]
	\node[draw, circle, fill=cyan!60] (loc1) {$s_1$};
	\node[below=of loc1, fill=cyan!60] (loc2) {$s_2$};
	\node[right=of loc1, circle, fill=lightgray] (dep1) {$d_1$};
	\node[below=of dep1, circle, fill=lightgray] (dep2) {$d_2$};
	\node[below right=of dep1, fill=orange!50, yshift=.5cm] (task1) {$t_3$};    
	\node[above right=of task1, fill=orange!50, yshift=-0.25cm] (task2) {$t_2$};
	\node[below right=of task1, fill=orange!50, yshift=0.5cm] (task3) {$t_4$};    
	\node[above=of task1, fill=orange!50, yshift=0.25cm] (task4) {$t_1$};
	\node[below=of task1, fill=orange!50, yshift=-0.25cm] (task5) {$t_5$};
	\node[right=of task2, fill=cyan!60, xshift=0.75cm] (dest1) {$e_1$};
	\node[right=of task3, fill=cyan!60, xshift=0.75cm] (dest2) {$e_2$};
	
	\draw (dep1) -- (task1);
	\draw (dep1) -- (task2);
	\draw (dep1) -- (task4);
	
	\draw[] (dep2) -- (task3);
	\draw[] (dep2) -- (task5);
	
	\draw[<->] (task1) -- (task2);
	\draw[<->] (task1) -- (task4);
	\draw[<->] (task2) -- (task4);
	\draw[<->] (task3) -- (task5);
	
	\draw[] (task1) -- (dest1);
	\draw[] (task1) -- (dest2);
	\draw[] (task2) -- (dest1);
	\draw[] (task2) -- (dest2);
	\draw[] (task3) -- (dest1);
	\draw[] (task3) -- (dest2);
	\draw[] (task4) to[out=0, in=165] (dest1);
	\draw[] (task4) -- (dest2);
	\draw[] (task5) -- (dest1);
	\draw[] (task5) to[out=0, in=-165] (dest2);
	
	\draw[] (loc1) -- (dep1);
	\draw[] (loc1) -- (dep2);
	\draw[] (loc2) -- (dep1);
	\draw[] (loc2) -- (dep2); 
	
	\draw[line width=3, solid, red] (loc1) -- (dep1);
	\draw[line width=3, solid, red] (dep1) -- (task4);
	\draw[line width=3, solid, red] (task4) -- (task2);
	\draw[line width=3, solid, red] (task2) -- (dest1);
	
	\draw[line width=3, longdash,blue] (loc2) -- (dep2);
	\draw[line width=3, longdash,blue] (dep2) -- (task5);
	\draw[line width=3, longdash,blue] (task5) -- (task3);
	\draw[line width=3, longdash,blue] (task3) -- (dest2);
\end{tikzpicture}
	}
		\caption{Example problem instance.}
		\label{fig:example_instance}
	\end{subfigure}%
	\begin{subfigure}[]{0.3\textwidth}
		\centering
		\includegraphics[width=0.7\textwidth]{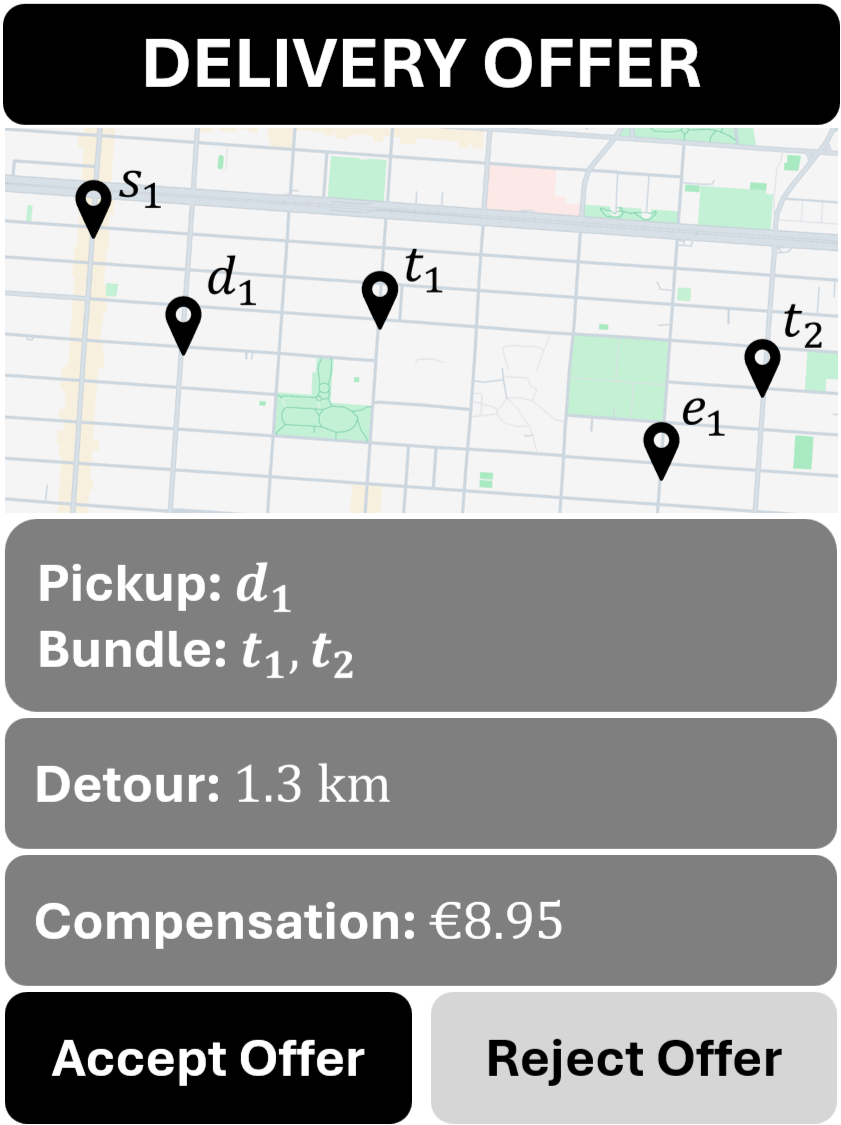}
		\caption{Driver interface}
		\label{fig:offer_interface}
	\end{subfigure}
	\caption{An instance with $D=2, M=5, N=2$, and two offers (indicated in red and blue) as well as the driver interface of the (red) offer to driver one.}
	\label{fig:instance_and_solution}
\end{figure}

\paragraph{Example problem instance} \cref{fig:instance_and_solution} illustrates a problem instance consisting of $|D|=2$ depots, $|M|=5$ tasks, and $|N|=2$ \crowdshippers. The network in \cref{fig:example_instance} contains nodes for driver origins ($s_1, s_2$) and destinations ($e_1, e_2$), depots ($d_1, d_2$), and task delivery locations ($t_1, \dots, t_5$). Directed arcs represent the feasible connections, i.e. origin-depot links, depot-task links for admissible depot-task pairs, task-task links between tasks that share a depot from which they can be delivered, and task-destination links. In particular, tasks $t_1, t_2$ and $t_3$ are restricted to depot $d_1$, and $t_4$ and $t_5$ to depot $d_2$. Consequently, arcs between tasks exist only if they can be served from the same depot.
The solid red and dashed blue arcs illustrate examples of feasible bundle-delivery paths.

A partial solution that correspond to the solid red path 
is depicted on the right, where an \crowdshipper receives a bundle offer (\cref{fig:offer_interface}), which provides relevant information 
such as the pickup depot, the bundle of tasks, the required detour, and the compensation.
The \crowdshipper then reveals his/her decision to accept or reject the bundle. 

A notable special case of this setting, widely examined in the crowdsourced delivery literature, occurs when in-store customers act as occasional drivers \citep{archetti_vehicle_2016, mancini_bundle_2022, cinar2024role}. In this case, customers visit a retail store for their own shopping and may accept an additional offer from the store to deliver packages to other customers located along their way home. 
	This case represents a specific instance of our setting
	in which the drivers’ origins coincide with the store, and the offer consists solely of a bundle of tasks collected at that store and the associated compensation.

\subsection{Mixed-integer nonlinear programming formulation}


In this section, we formulate the operator's problem as a MINLP. The MINLP formulation \eqref{eq:minlp} uses binary variables $\var \in \{0,1\}$ that indicate whether bundle $b\in B$ is offered to \crowdshipper $w\in N$ and continuous variables $C_{wb}\ge 0$ that model the compensation offered to \crowdshipper $w\in N$ for bundle $b\in B$.

\vspace{-1em}

\begin{subequations} \label{eq:minlp}
	\begin{align}
		\max \quad   & \sum_{w \in N}\sum_{b \in B} P_w(b,C_{wb})\left(\sum_{i \in b}c_i-C_{wb} \right)\var \label{eq:minlp:objective}  \\
		\mbox{s.t.}\quad & \sum_{w \in N} \sum_{b\in B : i\in b} \var  \leq 1 && i \in M \label{eq:minlp:constraint1} \\
		& \sum_{b \in B} \var \leq 1 && w \in N \label{eq:minlp:constraint2} \\
		& \var \in \{0,1\} &&  w \in N, b \in B \label{eq:minlp:variable1}  \\
		& C_{wb} \geq 0 &&  w \in N, b \in B \label{eq:minlp:variable2}
	\end{align}
\end{subequations}

The objective function (\ref{eq:minlp:objective}) maximizes the total expected savings obtained by offering bundles to the \crowdshippers. 
Constraints \eqref{eq:minlp:constraint1} and \eqref{eq:minlp:constraint2} ensure that each task is offered to at most one \crowdshipper and at most one bundle is offered to each \crowdshipper, respectively. 
Constraints \eqref{eq:minlp:variable1} and \eqref{eq:minlp:variable2} define the domains of the variables. 
Notice that the depot from which a bundle is delivered does not need to be explicitly represented in the formulation as we assumed that the optimal depot is implied when given a bundle and \crowdshipper. 

Formulation~\eqref{eq:minlp} cannot be solved directly by off-the-shelf software as it involves an exponential number of variables 
and because it features a non-linear objective function. 
In \cref{sec:methodology}, we will propose a solution algorithm that overcomes these challenges within reasonable computation time. 


\subsection{Modeling and fitting acceptance probability functions} \label{sec:acceptance-probability}

Logistic regression (or binary logit) models are typically used to predict binary outcomes, such as ``accept'' or ``reject'' decisions. For brevity, we refer readers to \citet{hosmer2013applied} for a comprehensive introduction. In the crowdsourced delivery literature, several studies employ these models to analyze the acceptance behavior of occasional drivers. For instance, \citet{devari_crowdsourcing_2017, hou_optimization_2022}, and \citet{barbosa_data-driven_2023} utilize logistic regression or binary logit models based on real-world survey data to predict acceptance decisions. \citet{hou_optimization_2022, barbosa_data-driven_2023}, and \citet{cinar2024role} incorporate logistic regression models within optimization frameworks designed to address crowdsourced delivery problems.


Although there is a growing body of literature addressing the acceptance behavior of participants in crowdsourced delivery systems, the aspect of modeling the acceptance behavior of occasional drivers has not received much attention yet.

Since a dedicated behavioral study is beyond the scope of this paper, we base the design of the acceptance probability function on the recent study by \citet{mohri_modeling_2024}, who estimate acceptance probabilities in a crowdsourced delivery system involving public transport passengers using detour, bundle weight, and compensation as predictors.
To illustrate our solution framework, we use the logistic acceptance probability function \eqref{eq:acceptance_model} to model acceptance probabilities of occasional drivers with the relevant predictors from the behavioral study, i.e. the detour $\Delta_{wb}$ required by occasional driver $w\in N$, the number of tasks in the bundle $b\in B$ (as a proxy for the bundle weight), and the offered compensation $C$. 
	\begin{equation} \label{eq:acceptance_model}
		P_w(b,C) = \frac{1}{1+e^{-(\alpha^{w}+ \beta_1^{w} \Delta_{wb} + \beta_2^{w} |b| + \gamma^{w} {C})}}
	\end{equation}

	Thereby, $\alpha^w, \beta_1^w<0, \beta_2^w<0, \gamma^w>0$ represent the coefficients for the intercept, detour, bundle size, and compensation, respectively. We assume that longer detours and larger bundle sizes require more effort from the driver and therefore reduce the acceptance probability while higher compensation makes an offer more attractive. In \cref{eq:acceptance_model}, this is captured by the fact that coefficients $\beta_1^w$ and $\beta_2^w$ are negative while  $\gamma^w$ is positive.


\blue{Our solution methodology is, however, not restricted to this specific function. All results and algorithms presented in the following are valid for} the generic logistic acceptance probability functions defined in \cref{def:logistic_prob}\blue{, which generalize \eqref{eq:acceptance_model} to arbitrary bundle- and driver-dependent predictors}. 
Additional benefits for choosing a logistic regression model is given in \cref{sec:model_fitting} where we also explain how such a model can be fitted.

\begin{definition}[Logistic acceptance probability] \label{def:logistic_prob}	
	The logistic function estimating the acceptance probability of \crowdshipper $w\in N$ for bundle $b\in B$ in exchange for compensation $C$ is defined as
	\begin{equation} \label{eq:acceptance_probability}
		P_w(b,C) = \frac{1}{1+e^{-(\alpha^w + B^w X_{wb} + D^w Y_{w} + \gamma^w{C})}}.
	\end{equation}
	Here, $\alpha^w$ is the coefficient of the intercept,  $B^w=\{\beta_1^w,\dots,\beta_n^w\}$ and $D^w=\{\delta_1^w, \dots, \delta_m^w\}$ are the coefficients of a set of $n+m$ predictors, and $\gamma^w>0$ is the coefficient for the compensation offered,  where the positive sign reflects that higher compensation increases the acceptance probability.
	Furthermore, $X_{wb}$ are $n$ non-negative and additive values of the predictors that depend on the bundle (and possibly the driver), $Y_{w}$ are $m$ values of predictors that depend on properties of the driver alone independent of the offered bundle, and $C$ is the compensation offered. 

	

\end{definition}

\section{Methodology}  \label{sec:methodology}

In this section, we discuss all components of the novel exact solution algorithm for the optimization problem formally introduced in the previous section. The algorithm is based on a linearization of formulation \eqref{eq:minlp} that is derived in \cref{sec:preliminaries} and determines optimal bundles and compensations for \crowdshippers under compensation-dependent acceptance probabilities. An overview of the complete exact algorithm and the interaction of its components is provided in \cref{fig:exact_algorithm_overview}.

\begin{figure}[htbp]
	\centering
	\footnotesize
	\begin{tikzpicture}[
		node distance=6mm,
		box/.style={draw, rounded corners, align=left, inner sep=6pt, text width=0.9\textwidth},
		arrow/.style={-Latex, thick}
		]
		
		
		\node[box] (corr) {
			\textbf{Step 0: Corridor Initialization (\cref{sec:heuristic_algo})} \\ 
			\textbf{Input:} LP relaxation of MILP~(\ref{eq:milp}) (Section 4.1) with empty set of columns (bundle variables)\\
			\textbf{Focus:} Generate initial columns \\ 
			\textbf{Output:} Set of initial columns
		};
		
		\node[box, below=of corr] (cg) {
			\textbf{Step 1: Column Generation (\cref{sec:column_generation,sec:pricing-algorithm})}\\ 
			\textbf{Input:} LP relaxation of MILP~(\ref{eq:milp}) restricted to initial columns\\
			\textbf{Focus:} Generate additional columns to solve LP relaxation of MILP~(\ref{eq:milp})\\
			\textbf{Output:} LP relaxation solution of MILP~(\ref{eq:milp}), set of generated columns, \textbf{upper bound}
		};
		
		\node[box, below=of cg] (milp) {
			\textbf{Step 2: MILP Heuristic (\cref{sec:exact_algo})}\\ 
			\textbf{Input:} MILP~(\ref{eq:milp}) restricted to generated columns\\
			\textbf{Focus:} Obtain heuristic solution by solving MILP~(\ref{eq:milp}) for restricted set of columns\\
			\textbf{Output:} Feasible solution, \textbf{lower bound}
		};
		
		\node[box, below=of milp] (enum) {
			\textbf{Step 3: Variable Enumeration \& Re-optimization (\cref{sec:exact_algo})}\\
			\textbf{Input:} LP relaxation of MILP~(\ref{eq:milp}) restricted to generated columns, upper and lower bounds \\
			\textbf{Focus:} Generate further variables possibly included in optimal solution, solve resulting MILP \\
			\textbf{Output:} \textbf{Optimal solution}
		};
		
		\draw[arrow] (corr) -- (cg);
		\draw[arrow] (cg) -- (milp);
		\draw[arrow] (milp) -- (enum);		
	\end{tikzpicture}
	\caption{Overview of the exact solution algorithm.}	 
	\label{fig:exact_algorithm_overview}

\end{figure}


In the following, we first describe the linearization of the MINLP~\eqref{eq:minlp} in \cref{sec:preliminaries}, which yields the MILP~(\ref{eq:milp}). \cref{sec:column_generation} details the column generation approach used to solve the LP relaxation of the MILP~(\ref{eq:milp}), and \cref{sec:pricing-algorithm} introduces the labeling algorithm for the associated pricing subproblem with a non-additive and non-linear objective function. The lack of additivity prohibits the use of established labeling algorithms \citep{sadykov2021bucket, costa2019exact} and previous attempts to handle ESPPRC variants with non-additive objective functions have proven inefficient \citep{rostami2021branch}. Hence, \cref{sec:pricing-algorithm} also introduces novel dominance and pruning rules, which constitute key methodological contributions and enable an efficient algorithm.

\cref{sec:exact_algo} details the remaining components of the proposed algorithm in \cref{fig:exact_algorithm_overview}, i.e., the MILP heuristic, variable enumeration and re-optimization. 
\cref{sec:heuristic_algo} introduces heuristic procedures that can be used to accelerate the overall algorithm (corridor initialization) as well as stand-alone heuristic methods obtained by leveraging individual components of the complete exact approach.
Finally, all of the above-mentioned components consider formulations featuring the generic logistic function with arbitrary predictors defined in \cref{eq:acceptance_probability}. \cref{sec:concrete-probability-function}
discusses adaptations of our generic algorithms to 
\blue{the logistic acceptance probability function \eqref{eq:acceptance_model} specified in \cref{sec:acceptance-probability}, which is used in our computational study}.

\subsection{Linearization of formulation (2)} \label{sec:preliminaries}

We first show how to linearize formulation \eqref{eq:minlp} under the mild assumption of \textit{separability}, which requires that the probability that an \crowdshipper accepts a bundle only depends on the driver, 
the 
bundle and the 
compensation. Neither the decision whether an operator can offer a certain bundle to a given \crowdshipper nor the compensation value offered for a bundle depend on bundles offered to other drivers or their compensations. As a consequence, 
decisions about offered bundles and compensations can be made independently for each driver / bundle pair.
As the setting considered in this paper does not include any restriction on the operator's bundling decisions (such as a limit on the total number of bundles or a budget constraint on the overall compensation), this separability assumption holds as long as the acceptance probability functions of \crowdshippers do only depend on the compensation, properties of the bundle and the considered driver. 

\citet{cinar2024role} show how to derive optimal compensations for single tasks when assuming separability. \cref{prop:compensation_values,th:compensation_values_logistic}, whose proofs are given in Electronic Companion~\ref{sec:proofs}, generalize these results to  the case of offering bundles instead of individual tasks to \crowdshippers. 



\begin{proposition}\label{prop:compensation_values}	
	\compensationvalues	
\end{proposition}

\cref{th:compensation_values_logistic} shows how to determine the optimal compensation value for a bundle $b\in B$ offered to \crowdshipper $w\in N$ for the logistic acceptance probability function introduced in \cref{def:logistic_prob}. 

\begin{theorem} \label{th:compensation_values_logistic}	
	\optimalcompensation
\end{theorem}

%

\cref{prop:compensation_values} implies that formulation~\eqref{eq:milp} is a linear reformulation of 
\eqref{eq:minlp} 
obtained by calculating optimal compensations $C^*_{wb}$ 
for every pair of occasional driver $w\in N$ and bundle $b\in B$. 

\begin{subequations} \label{eq:milp}
	\begin{align}
		\max \quad   & \sum_{w \in N}\sum_{b \in B} P_w(b,C^*_{wb})\left(\sum_{i \in b}c_i-C^*_{wb} \right)\var \label{eq:milp:objective}  \\
		\mbox{s.t.}\quad & \sum_{w \in N} \sum_{b\in B : i\in b} \var  \leq 1 && i \in M \label{eq:milp:constraint1} \\
		& \sum_{b \in B} \var \leq 1 && w \in N \label{eq:milp:constraint2} \\
		& \var \in \{0,1\} &&  w \in N, b \in B \label{eq:milp:variable1} 
	\end{align}
\end{subequations}

\subsection{Column generation} \label{sec:column_generation}

In this section, we describe the column generation procedure developed for solving the linear relaxation of formulation~\eqref{eq:milp} in case of logistic acceptance probabilities given in \cref{def:logistic_prob}. In each iteration of the column generation process, the 
restricted linear master program (RMP), i.e. the linear relaxation of formulation \eqref{eq:milp} containing a limited set of bundles, is solved and its solution is used to identify bundles with positive reduced costs in the pricing subproblem and add them to the RMP. 
The procedure may be initialized with the empty bundle for each driver
	or with a set of promising bundles that can, e.g., be generated by the 
	corridor heuristic 
	described in \cref{sec:heuristic_algo}.
The procedure terminates if no further bundle with positive reduced costs exist. In this case, the objective function value is an upper bound on the optimal expected cost savings. 
Let $\pi_i\ge 0$, $i\in N$, and $\mu_w\ge 0$, $w\in N$, be the optimal dual variable values 
associated with constraints~\eqref{eq:milp:constraint1} and \eqref{eq:milp:constraint2}, respectively. Then, the reduced cost $\tilde{c}_{w,b}$ of bundle $b\in B$ offered to driver $w\in N$ is defined as 
\begin{equation} \label{eq:sp:obj}
	\tilde{c}_{wb} = P_w(b,C^*_{wb})\left(\sum_{i \in b}c_i-C^*_{wb} \right) - \sum_{i \in b}\pi_i - \mu_w.    
\end{equation}


The following proposition, the proof of which can be found in the \apdx, shows how to calculate \eqref{eq:sp:obj} in case of logistic acceptance probabilities.

\begin{proposition} \label{prop:reduced-cost}
\reducedcost
\end{proposition}

\subsection{Solving the pricing subproblem}\label{sec:pricing-algorithm}

In this section, we describe a labeling algorithm for solving the pricing subproblem that is executed independently for each driver in each iteration of the column generation approach. For each driver, this labeling algorithm either identifies bundles with positive reduced costs that are added to the RMP or proves that no such bundles exist. The process of adding bundles with positive reduced costs and resolving the updated RMP is repeated until there are no driver for which bundles with positive reduced costs exist.
	

Recall that bundles $b\in B$ correspond to ordered sets of tasks (i.e., 
$b=(i_1, i_2 \dots, i_\ell)$ such that $\{i_1, i_2, \dots, i_\ell\}\in 2^M$) and that for each bundle $b$ and \crowdshipper $w\in N$ there exists an optimal (pickup) depot $d_{wb}\in D$. Thus, a bundle $b\in B$ offered to driver $w\in N$ corresponds to a path $(s_w, d_{wb},i_1, i_2, \dots, i_\ell,e_w)$ for driver $w$ starting at the driver's current location $s_w$, visiting the optimal depot $d_{wb}$ and all tasks of bundle $b$, before ending at the drivers destination $e_w$. Therefore, we can solve the pricing subproblem by identifying an elementary 
path with maximum reduced costs \eqref{eq:sp:obj} for each \crowdshipper $w\in N$. Since the objective function~\eqref{eq:sp:obj} is non-additive and non-linear, we propose a tailored labeling algorithm based on 
novel dominance and label pruning rules whose main components we explain in the following.

In the labeling algorithm, we associate a label $L_w^p=(b^p,d^p,q^p,X^p,\bar{C}^p,R^p,\tilde{c}^\blue{p})$ to each path $p=(s_w,d,i_1, \dots i_\ell,e_w)$ corresponding to a bundle $(i_1, \dots, i_\ell)$ and pickup location $d\in D$ offered to driver $w\in N$ . 
Thereby, $b^p=(i_1, \dots, i_\ell)$ is the offered bundle, $d^p\in D$ the pickup location, $X^p=(X^p_1, \dots, X^p_{n})$ is the vector consisting of the non-negative values of the $n$ bundle-dependent predictors of driver $w$, cf.~\cref{def:logistic_prob}. Furthermore, $q^p=\sum_{i\in b^p} q_i$ is the used capacity, $\bar{C}^p=\sum_{i\in b^p} c_i$ is the total 3PL cost of the bundle, $R^p\subseteq T(d^p)$ is the set of 
reachable tasks, i.e. tasks that can be delivered from depot $d^p$ by occasional driver $w$ after delivering all tasks included in bundle $b^p$; and $\tilde{c}^p$ is the reduced cost of bundle $b^p$ with pickup location $d^p$ offered to occasional driver $w\in N$.
In the following, we will avoid case distinctions by using notation $\mathcal{T}^p$ to refer to the last location of path $p$ before a driver reaches its end location, i.e., $\mathcal{T}^p=i_\ell$ if $b^p \ne \emptyset$, while $\mathcal{T}^p=d^p$ otherwise. An example illustrating the following initialization and extension steps is provided in \cref{sec:label_ex}.

\paragraph{Initialization}
For each driver $w\in N$, we initially create one label for each depot $d\in D$ featuring an empty bundle, i.e., one label $L_w^p=(\emptyset, d, 0, \bar{X}^p, 0, T(d), -\mu_w)$ for each $d\in D$ such that $p=(s_w,d,\emptyset,e_w)$ and where $\bar{X}^p$ is vector of initial values of the bundle-dependent predictors.
\paragraph{Extension} 
Adding a task $t\in R^p$ to a label $L_w^p$ associated with bundle $b^p=(i_1, \dots i_\ell)$ generates a label $L_w^{f}$ such that $b^f=(i_1, \dots, i_\ell, t)$, $d^f=d^p$, $q^f=q^p + q_t$, 
$X_i^{f} = X_i^p + \mathcal{U}_i(\mathcal{T}^p,t)$ for $i=1, \dots, n$,
$\bar{C}^{f} =\bar{C}^{p} + c_t$, $R^{f}=R^p\setminus (\{t\}\cup \{ i\in R^p :  q^f + q_i > Q^w\})$, and $\tilde{c}^f=P_w(b^f,C^*_{w{b^f}})\left(\bar{C}^f-C^*_{w{b^f}} \right) - \sum_{i \in b^f}\pi_i - \mu_w$. 
%
Thereby, non-negative values $X^p_i$ are assumed to be additive, cf.\ \cref{def:logistic_prob}. Thus, we 
define the generic function
$\mathcal{U}_i : (M\cup D) \times M \mapsto \mathbb{R}_{\ge0}$ 
that takes as input
the last location $\mathcal{T}^p$ in path $p$, 
and the new task $t$ and returns 
the new value of the bundle-dependent predictor with index $i=1, \dots, n$. 

\paragraph{Dominance and label pruning}
The non-additivity and non-linearity of the reduced cost
prohibits the use of classical dominance rules 
from the literature which are known to have a huge impact on the performance of column-generation based solution methods using labeling algorithms for solving the pricing subproblem. \cref{prop:dominancerules,prop:rcbasedpruning}, whose proofs are given in the \apdx, propose
new dominance and 
label pruning rules to reduce the number of generated labels.

\begin{proposition}{\textbf{Dominance}} \label{prop:dominancerules}
\dominance
\end{proposition}

Let ${U}_i^\mathrm{max}$ and ${U}_i^\mathrm{min}$ be an upper and a lower bound on the change of the value of the predictor $X^p_{i}$ in all potential, future extensions of the considered label calculated as the maximum and minimum feasible additions, i.e. $U^\mathrm{max}_i=\max_{s\in \mathcal{T}^p\cup R^p, t\in R^p} \mathcal{U}_i(s,t)$, 
$U^\mathrm{min}_i=\min_{s\in \mathcal{T}^p\cup R^p, t\in R^p} \mathcal{U}_i(s,t)$. Then, the following proposition holds.

\begin{proposition} \label{prop:rcbasedpruning}
\rcbasedpruning
\end{proposition}

\cref{prop:rcbasedpruning} implies \cref{corr:rcbasedpruning} in which an upper bound on the reduced cost of any extension of a label is computed by over-estimating the values of bundle-dependent predictors, the 3PL cost, and under-estimating the aggregated dual values. Subsequently, \cref{corr:conditionforpruning} provides the detailed rule that will be used to prune labels for which no extension can have positive reduced cost.

\begin{corollary} \label{corr:rcbasedpruning}
\rcbasedpruningcorr
\end{corollary}

\begin{corollary}[\textbf{Reduced cost-based label pruning}] \label{corr:conditionforpruning}
\conditionforpruning
\end{corollary}

A practical method for estimating values of $k_{max}$ 
is discussed in the introduction of \cref{sec:results}.

\subsection{Lower bounding and enumeration} \label{sec:exact_algo}

The column generation procedure described in the previous two subsections solves the continuous relaxation of formulation~\eqref{eq:milp} and, therefore, provides an upper bound \blue{$\bar z$} of the optimal solution value. Subsequently, we first compute a feasible solution providing a lower bound \blue{$\underline{z}$} by solving formulation~\eqref{eq:milp} with a restricted set of bundle variables equal to those added to the RMP during the column generation process. 
Afterwards, we enumerate all columns (i.e., bundle variables) not-yet included in the RMP whose reduced cost are within the optimality gap, i.e., greater than or equal to the difference between the obtained lower and upper bound. Resolving the variant of formulation~\eqref{eq:milp} restricted to all variables added during the column generation process and those found in the enumeration step yields an optimal solution.

\citet{baldacci2008exact} show how to successfully use variable enumeration and re-optimization instead of branch-price-and-cut for solving the capacitated vehicle routing problem. As discussed in the survey of \citet{costa2019exact}, this idea has also been successfully used in exact solution methods for other routing problems. More recently, \citet{paradiso2020exact} and \citet{yang2023exact} apply this methodology within exact solution frameworks for multitrip vehicle routing problems with time windows. The existing literature indicates that such approaches can yield substantial speed-ups over traditional branch-price-and-cut  methods when the LP relaxation yield tight bounds and high-quality heuristic solutions are available too.  Since this is the case for our approach (see, \cref{sec:results}) we also adopt variable enumeration and re-optimization.
%

During the enumeration of columns, the dominance rule of \cref{prop:dominancerules} is not used, as a dominated bundle whose reduced cost is within the optimality gap may be part of the optimal solution. Instead, we use the simpler dominance rule given in \cref{corr:dominance_simple} that considers two labels corresponding to paths that visit the same set of tasks.

\begin{corollary}[Route-based dominance]	\label{corr:dominance_simple}
\dominancesimple
\end{corollary}

\cref{corr:dominance_simple} holds since the four conditions of \cref{prop:dominancerules} trivially hold for any two labels satisfying its conditions. 
 This simpler dominance rule can be applied in the enumeration phase since a solution offering a bundle corresponding to (an extension of) label $L_w^f$ can always be improved by offering a bundle corresponding to (the same extension of) label $L_w^p$ to the same driver $w\in N$.

\subsection{Heuristic approaches} \label{sec:heuristic_algo}

For large instances, the exact algorithm presented in the preceding subsections requires too much computational time. Therefore, we discuss two possibilities to derive heuristic variants that enhance scalability towards such instances.
The first option is to simply skip the last step (i.e., variable enumeration and re-optimization) of the exact algorithm. Several variants of this heuristic approach which computes both a lower and an upper bound on the optimal solution value (and therefore a measure of solution quality) will be considered in our computational study.
The second type of heuristic approach, referred to as \emph{corridor heuristic}, is described in the following. This heuristic can be used either as a stand-alone method for large-scale instances or to speed up either the exact algorithm or the first heuristic variant.
%
%
%
\paragraph{Corridor Heuristic}
The corridor heuristic is 
a modified version of the method introduced by \citet{mancini_bundle_2022}. It focuses on a subset of tasks located within a specific area defined by an angular region, referred to as a ``corridor'', centered along the direction of the occasional driver's destination $e_w$, with its apex at the driver’s current location $s_w$.

To define the corridor, we consider the vector $\vec{v}_w = e_w - s_w$, which serves as the central axis of the corridor. Given a half-angle $\theta > 0$, the corridor includes all tasks $t \in M$ such that both the delivery location of the task and at least one associated depot lie within an angle of $\theta$ from the central axis.

Formally, the corridor subset $Cor_{w,\theta}$ for driver $w$ and angle $\theta$ is defined as:
$$
Cor_{w,\theta} = \left\{ t \in M : \angle(\vec{v}_w, \vec{x}_t) < \theta \wedge \exists d \in D: t \in T(d) \wedge \angle(\vec{v}_w, \vec{x}_d) < \theta \right\},
$$
where $\angle(\cdot,\cdot)$ denotes the smaller angle between two vectors in the plane, and $\vec{x}_t = t - s_w$ and $\vec{x}_d = d - s_w$ are the vectors from the driver’s current location to the task and to the depot $d$. An illustration of the corridor space is shown in \cref{fig:corridor_space}.

\begin{figure}
	\centering
	\resizebox{0.42\textwidth}{!}{%
	\begin{tikzpicture}[
		scale=0.5,
		->,
		>=stealth,
		node distance=1.5cm and 1.5cm,
		every node/.style={draw, circle, minimum size=0.75cm, inner sep=0pt}
		]
		\node[draw, circle, fill=cyan!60] (loc1) {$s_1$};
		\node[below=of loc1, fill=cyan!60] (loc2) {$s_2$};
		\node[right=of loc1, circle, fill=lightgray] (dep1) {$d_1$};
		\node[below=of dep1, circle, fill=lightgray] (dep2) {$d_2$};
		\node[below right=of dep1, fill=orange!50, yshift=.5cm] (task1) {$t_1$};    
		\node[above right=of task1, fill=orange!50, yshift=-0.25cm] (task2) {$t_2$};
		\node[below right=of task1, fill=orange!50, yshift=0.625cm] (task3) {$t_3$};    
		\node[above=of task1, fill=orange!50, yshift=0.25cm] (task4) {$t_4$};
		\node[below=of task1, fill=orange!50, yshift=0.25cm] (task5) {$t_5$};
		\node[right=of task2, fill=cyan!60, xshift=0.75cm, yshift=0cm] (dest1) {$e_1$};
		\node[right=of task3, fill=cyan!60, xshift=0.75cm, yshift=0cm] (dest2) {$e_2$};
	
		\coordinate (start) at (loc1.center);

\begin{scope}[on background layer]
	\filldraw[
	fill=blue!20, opacity=0.3,
	draw=blue,    thick
	]
	(start) -- ($(start)+(-15:20cm)$)
	arc[start angle=-16,end angle=17.5,radius=20cm]
	-- cycle;

\path (start) -- (dest1)
coordinate[pos=0.475] (midL)
coordinate[pos=0.50] (mid)
coordinate[pos=0.525] (midR);

\draw[dashed, thick,-] (start) -- (midL);
\node[draw=none, fill=none, inner sep=0pt] at (mid) {$\vec v_1$};
\draw[dashed, thick,-] (midR) -- (dest1);

\coordinate (arcBot) at ($(start)+(18.25:20cm)$);

\draw pic [
draw,           
->,            
angle radius=7.5cm
] {angle = dest1--start--arcBot};
\node [draw=none,fill=none, xshift=-1.9cm, yshift=-1.9cm] at (arcBot) {$\theta$};

\end{scope}

%
%
%
%
%
%
	\end{tikzpicture}
}
	\caption{Example corridor space defined around $\vec{v}_1=e_1-s_1$ for occasional driver 1 containing $d_1$, $t_2$, and $t_4$.}
	\label{fig:corridor_space}
\end{figure}

The corridor heuristic starts with generating columns that consider only subsets of tasks in the corridor spaces $Cor_{w,\theta}$ for each occasional driver $w \in N$, where $\theta$ is an externally defined parameter. The search continues until no more such bundles with positive reduced cost remain. While this step yields an optimal solution within the restricted corridors, it does not constitute an upper bound for the full problem unless $\theta = 180^\circ$, i.e., the entire task space is considered.

The corridor heuristic can be employed as a standalone heuristic. In this case, the MILP heuristic is executed by only considering the columns identified in the first step. This creates a valid lower bound for the complete problem without any solution quality measurement due to the lack of a valid upper bound. The bundles identified in this first step can also be employed as the initial column set for solving the full problem. In this case, the exact algorithm (or the heuristic variant that skips the variable enumeration and re-optimization step) is executed with the warm start solution identified by the corridor heuristic. 



\subsection{\texorpdfstring{\blue{Adaptations to the concrete} logistic acceptance probability function}{Adaptations to the concrete logistic acceptance probability function}}\label{sec:concrete-probability-function}

The algorithm presented above uses a generic logistic acceptance model with abstract bundle- and driver-dependent predictors.
\blue{Here, we discuss its adaptation to the concrete logistic acceptance probability function \eqref{eq:acceptance_model} introduced in \cref{sec:acceptance-probability}, which is used in our computational study.}
Using the preceding results, we show how to prune the reachable node set, and consequently reduce the required label extensions considering the specific function.

The following \cref{corr:detourupperbound} establishes an upper bound on the maximum detour for which a label can have positive reduced cost. This result is used in \cref{corr:detourlimit} to derive a rule under which tasks can be eliminated from the set of potential extensions of a label.

\begin{corollary}[Detour upper bound] \label{corr:detourupperbound}
	\detourupperbound
\end{corollary}

\begin{corollary}[Detour-based reduction of successors] \label{corr:detourlimit}	
\detourbasedpruning	
\end{corollary}

\section{Experimental Design} \label{sec:experiments}

In this section, we describe the design of the computational study we conduct to evaluate the performance of the proposed algorithms, their benefits compared to a sequential approach inspired by the recent literature, and analyze the sensitivity towards key problem features. 

\subsection{Benchmark instances}

We generated a set of benchmark instances whose locations of tasks and \crowdshippers are inspired by the dataset introduced by \citet{mancini_bundle_2022}. However, since their instances are not large enough to assess the scalability and performance of our algorithm, 
created
larger ones that better stress-test our approach. Specifically, we generated additional sets of larger instances where the locations are distributed according to the same parameters as the ones used for the instances in \citet{mancini_bundle_2022}.
Similar to their 
instances, we consider a single-depot environment where the \crowdshippers are in-store customers.
Task delivery locations and \crowdshipper destinations are uniformly distributed within a square region defined by \(x, y \in [-5, 5]\), which is assumed to represent a $10\,\text{km} \times 10\,\text{km}$ service zone. Task loads \(q_i\) are uniformly drawn integers from the range \([10, 30]\), while each \crowdshipper{} has a capacity of 100 units. 3PL cost is set at a flat rate of \euro4.95, reflecting the fixed package delivery fee charged by the largest Dutch package delivery service provider \citep{postnl_prijzen_2024}. Unlike the original instances, the store is assumed to be located at the origin. 

We generate instances with varying numbers of tasks $|M|\in \{30, 60, 90, 120\}$ and fractions $p\in \{0.1, 0.2, 0.3, 0.4, 0.5\}$ to determine the numbers of drivers $|N|=p\cdot |M|$.
To ensure comparability between instances, we started by generating ten \emph{full-scale instances} as a base set.
 We chose this design in order to be able to draw conclusions about the impact of different ratios of tasks to \crowdshippers while keeping the distribution of locations of tasks and \crowdshippers as well as the \crowdshippers' behavioral classes comparable across instances. 
Each full-scale instance consists of 
120 tasks and 60 drivers. We further diversified the instances by employing distinct driver class allocation patterns, which assign drivers to one of the three unique behavioral classes that will be described 
in \cref{sec:behavioralclasses}.
We created eight such patterns for each full-scale instance: three single-class patterns (each consisting of drivers exclusively from one behavioral class) and five mixed-class patterns. 
For each mixed-class pattern, we ensure for every instance defined by $|M|$ and $p$ that the driver pool contains equal numbers from each class. Moreover, across different $|M|$ and $p$ levels, the drivers in the corresponding full-scale instances are kept in the same class. Class assignments are randomized, subject to these constraints.
This led to 
80 full-scale instances. 
%
For each of the 80 full-scale instances, reduced instances were derived by selecting the first $|M|$
tasks and
the first $p\cdot |M|$
drivers. 
Our approach 
ensures that each instance contains an equal number of drivers from each class and that the same driver is consistently assigned to the same class across instances. 
This procedure led to 20 unique instances of different dimensions generated for each full-scale instance, resulting in a total number of 1,600 of instances. 

\subsection{Behavioral classes} \label{sec:behavioralclasses}

\begin{figure}
\centering
\includegraphics[width=0.70\textwidth]{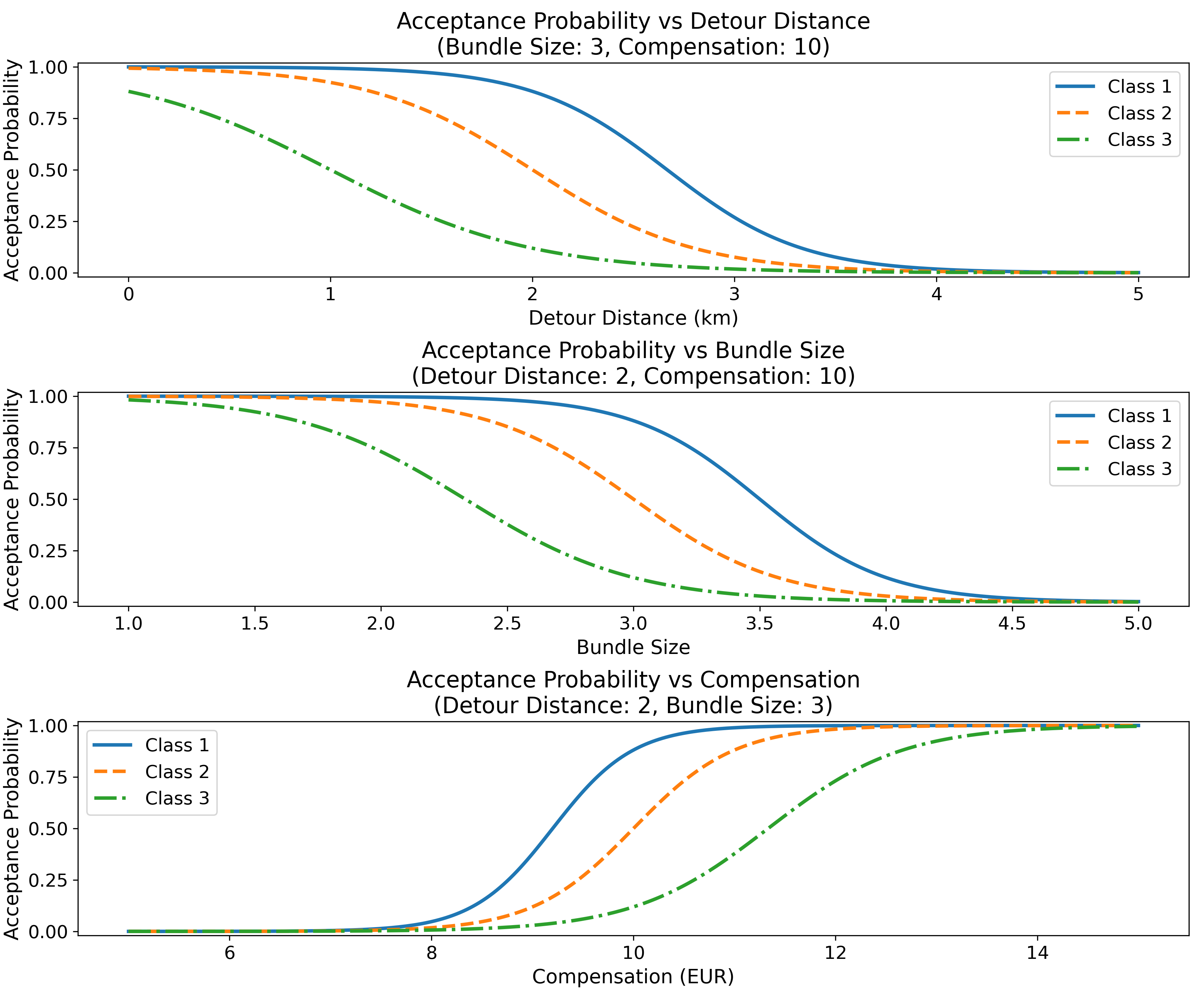}	
\caption{Acceptance probabilities of different \crowdshipper classes towards different offers}
\label{fig:behavior}
\end{figure}

Motivated by the acceptance probability functions introduced by \citet{mohri_modeling_2024}, which identified three latent behavioral classes based on real-world survey data, we defined three distinct \crowdshipper{} classes. \cref{tab:class_coefficients} lists the coefficients (\(\alpha, \beta_1, \beta_2, \gamma\)) for each class, capturing the trade-offs between effort (detour distance, bundle size), and compensation sensitivity. Specifically, Class 1 represents highly effort-averse drivers who are strongly motivated by compensation, Class 2 includes drivers with moderate sensitivity to both effort and compensation, and Class 3 comprises drivers who tolerate effort but are less influenced by compensation. 
\cref{fig:behavior} illustrates how acceptance probabilities vary with detour distance, bundle size, and compensation across the three classes. The figure shows that Class 1 and 3 yield the highest and lowest acceptance probabilities, respectively, towards the same offer. Also, Class 1 exhibits the steepest acceptance probability curves, whereas Class 3 has the least steep curves across varying offer attributes.
\blue{These coefficients are stylized profiles designed for controlled multi-class comparisons rather than direct estimates. As complementary robustness checks, \cref{sec:stress-test} reports the results on our benchmark instances by considering  the perturbed class coefficients based on those in \cref{tab:class_coefficients}, the published estimates from \citet{mohri_modeling_2024}, and the estimates fitted on a real-world dataset from Meituan published for the INFORMS TSL Data-Driven Research Challenge \Citep{wang2025meituan}.}

\begin{table}
	\centering
	\caption{Behavioral Class Coefficients}
	\label{tab:class_coefficients}
	\renewcommand{\arraystretch}{1.2}
	\begin{tabular}{lcccc}
		\toprule
		\textbf{Class} & \(\alpha\) (Intercept) & \(\beta_1\) (Detour) & \(\beta_2\) (Bundle Size) & \(\gamma\) (Compensation) \\ 
		\midrule
		\textbf{Class 1} & -5.0 & -3.0 & -4.0 & 2.5 \\ 
		\textbf{Class 2} & -4.5 & -2.5 & -3.5 & 2.0 \\ 
		\textbf{Class 3} & -4.0 & -2.0 & -3.0 & 1.5 \\ 
		\bottomrule
	\end{tabular}
\end{table}

\subsection{Algorithm variants}

To evaluate the trade-off between computational speed and solution quality as well as the contribution of individual algorithmic features, we evaluated the performance of several algorithm variants. 
\cref{tab:algorithm_variants} summarizes the features of each algorithm variant. 

The \textbf{Heuristic Baseline (H-B)} serves as a baseline method. It applies column generation without the 
dominance and detour limit features introduced in \cref{prop:dominancerules,corr:detourlimit} followed by the MILP heuristic. The \textbf{Heuristic with Dominance (H-D)} extends H-B by enabling dominance while keeping the detour limit off. The \textbf{Heuristic with Dominance and Detour Limit (H-DD)} incorporates both dominance and detour limits. Lastly, the \textbf{Heuristic with Dominance, Detour Limit, and Corridor Search (H-DDC)} integrates initial corridor search alongside dominance and detour limits. All these variants compute both upper and lower bounds, allowing for solution quality assessment.
Next, we introduce the \textbf{Heuristic Corridor Search (H-C)}, which focuses solely on the driver-centric corridor search space during the column generation phase before computing a lower bound and the corresponding solution using the MILP heuristic. 
This approach computes no upper bound. 
%
Finally, the exact variants extend their heuristic counterparts by incorporating the variable enumeration and re-optimization phase to compute the optimal solution. Specifically, \textbf{E-DD} builds upon H-DD, while \textbf{E-DDC} is derived from H-DDC.

\begin{table}
\small
	\centering
	\caption{Algorithm variants and their features}
	\label{tab:algorithm_variants}
	\begin{tabular}{lcccccc}
		\toprule
		\textbf{Variant} & \textbf{Corridor Search} & \textbf{Dominance} & \textbf{Detour Limit} & \textbf{UB} & \textbf{LB} & \textbf{Exact Solution} \\
		\midrule
		H-B   & No  & No  & No  & Yes & Yes & No \\
		H-D   & No  & Yes & No  & Yes & Yes & No \\
		H-DD  & No  & Yes & Yes & Yes & Yes & No \\
		H-DDC & Yes & Yes & Yes & Yes & Yes & No \\ \midrule		
		H-C   & Yes & Yes  & Yes  & No  & Yes & No \\
		\midrule
		E-DD  & No  & Yes & Yes & Yes & Yes & Yes \\
		E-DDC & Yes & Yes & Yes & Yes & Yes & Yes \\
		\bottomrule& 
	\end{tabular}
\end{table}

\section{Computational study} \label{sec:results}

This section presents a comprehensive evaluation of our algorithms. First, we analyze the performance of the algorithmic variants across different instance sizes and behavioral driver classes. 
The impact of algorithmic components is evaluated on mixed-class instances with $p=0.5$, and the best-performing variants are further assessed on single- and mixed-class instances with varying $p$ levels.
We also assess the benefits of jointly optimizing bundle generation, assignment, and compensation decisions	versus an approach where these decisions are made sequentially as in \citet{mancini_bundle_2022, macrina2024bundles}. 
We further conduct a sensitivity analysis to investigate the influence of key parameters, such as the number of tasks, the driver-to-task ratio, and driver classes, on the structure of the solutions. 
Finally, we provide managerial insights focusing on the impact of accounting for driver heterogeneity from the operator's and the drivers' perspectives.

Computations were carried out using the Snellius National Supercomputer, the high-performance computing facility of the Netherlands \citep{SURF_Snellius}. Each experiment was executed on a single core of an AMD EPYC 9004 series processor with a base clock speed of 2.4 GHz. The memory limit 
was set to 16 GB and the time limit to  3 hours. 
The proposed algorithms were implemented in C++ and solved using IBM CPLEX 22.1.1, running on a single thread for both LP and MILP formulations. We used the following settings:
%
In each iteration, we pre-maturely terminate the pricing algorithm for a driver when at least \blue{$\kappa=100$} variables with positive reduced costs have been found for this driver.
Furthermore, we set the half-angle parameter of corridor search to $\theta = 36^\circ$\blue{, as this combination performs the best, see \cref{app:theta-k}.} 
 
 Finally, we determine the maximum number of additional tasks $k_{max}$, that can be added to a label, by sorting all reachable tasks in non-decreasing order of weight. Starting from the lightest task, we then calculate how many of them can fit within the remaining capacity.

\subsection{Algorithmic performance of exact variants}\label{sec:results:exact}

We first analyze and compare the performance of exact variants E-DD considering the dominance and detour limit versus E-DDC which additionally incorporates the initial corridor search. Note that the significant positive impact of dominance and detour limit will be assessed in \cref{sec:results:heuristic} in the context of the heuristic solution variants. 
For runtime analysis, we set the runtime to 3 hours (corresponding to the time limit) for all instances in which we do not reach an optimal solution within the time limit or for which the memory limit was exceeded.

%
%

\begin{table}
	\caption{Numbers of instances solved and average runtimes in seconds (50 instances per configuration).}
	\label{tab:exact_runtimes}
	\centering
	\blue{
	\begin{tabular}{ccccrr}
		\toprule
		\multicolumn{2}{c}{\textbf{Instances}} & \multicolumn{2}{c}{\textbf{Solved Instances}} & \multicolumn{2}{c}{\textbf{Average Runtime}} \\
		\cmidrule(lr){1-2} \cmidrule(lr){3-4} \cmidrule(lr){5-6}
		Tasks & Drivers & \textbf{E-DD} & \textbf{E-DDC} & \textbf{E-DD} & \textbf{E-DDC} \\
		\midrule
		30 & 15 & 50 & 50 & 0.4 & 0.3 \\
		60 & 30 & 50 & 50 & 11.3 & 8.1 \\
		90 & 45 & 49 & 50 & 539.5 & 301.5 \\
		120 & 60 & 42 & 43 & 3101.1 & 3065.4 \\
		\bottomrule
	\end{tabular}
	}
\end{table}

\begin{figure}
	\captionsetup[subfigure]{aboveskip=-0.5pt,belowskip=-0.5pt}
\begin{subfigure}[t]{0.49\textwidth}
	\centering
\includegraphics[width=0.9\textwidth]{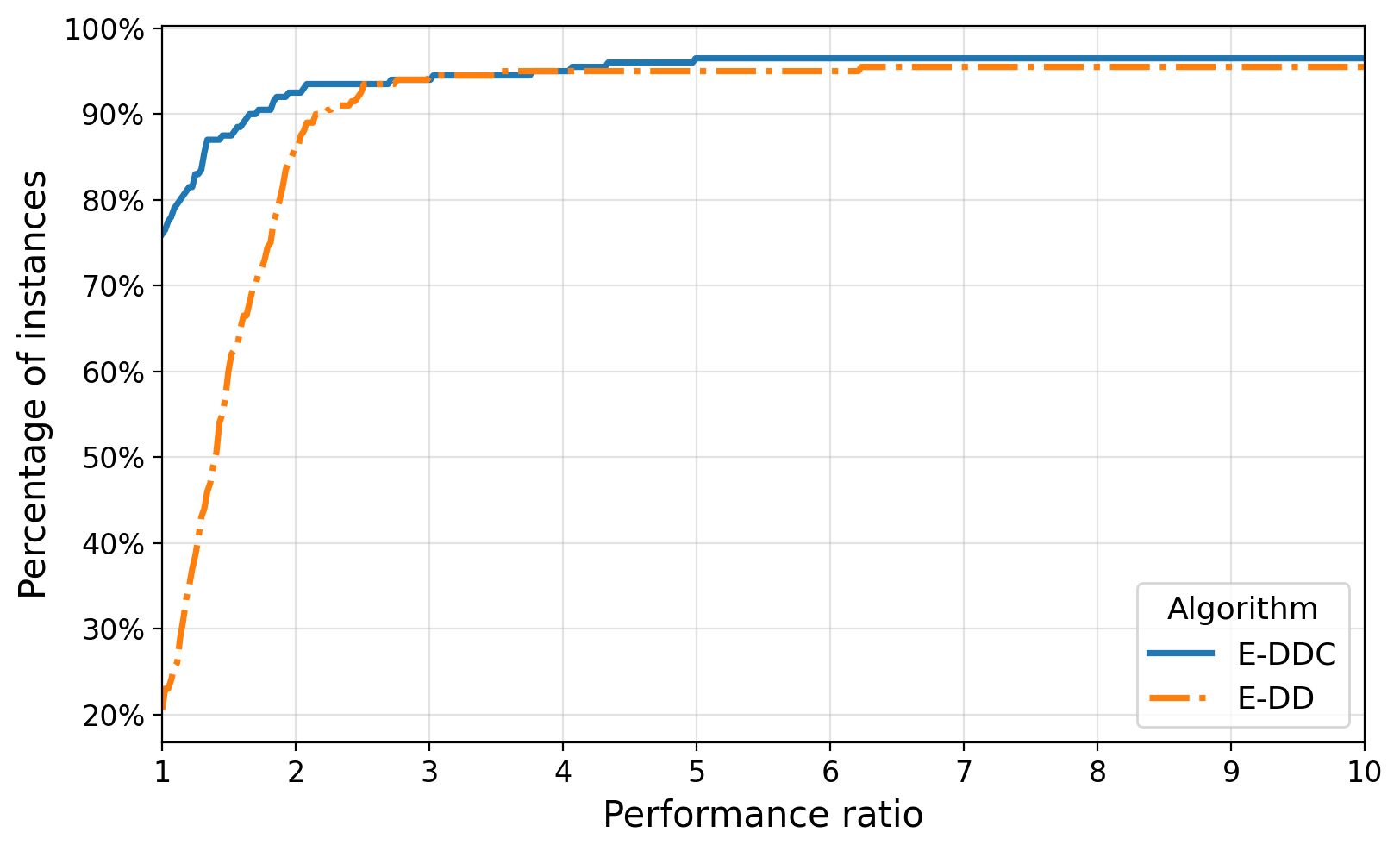}
\caption{\blue{Performance profiles of exact algorithms.}}
\label{fig:exact_perf_prof}
\end{subfigure}	\hfill	
\begin{subfigure}[t]{0.49\textwidth}
	\centering
\includegraphics[width=0.9\textwidth]{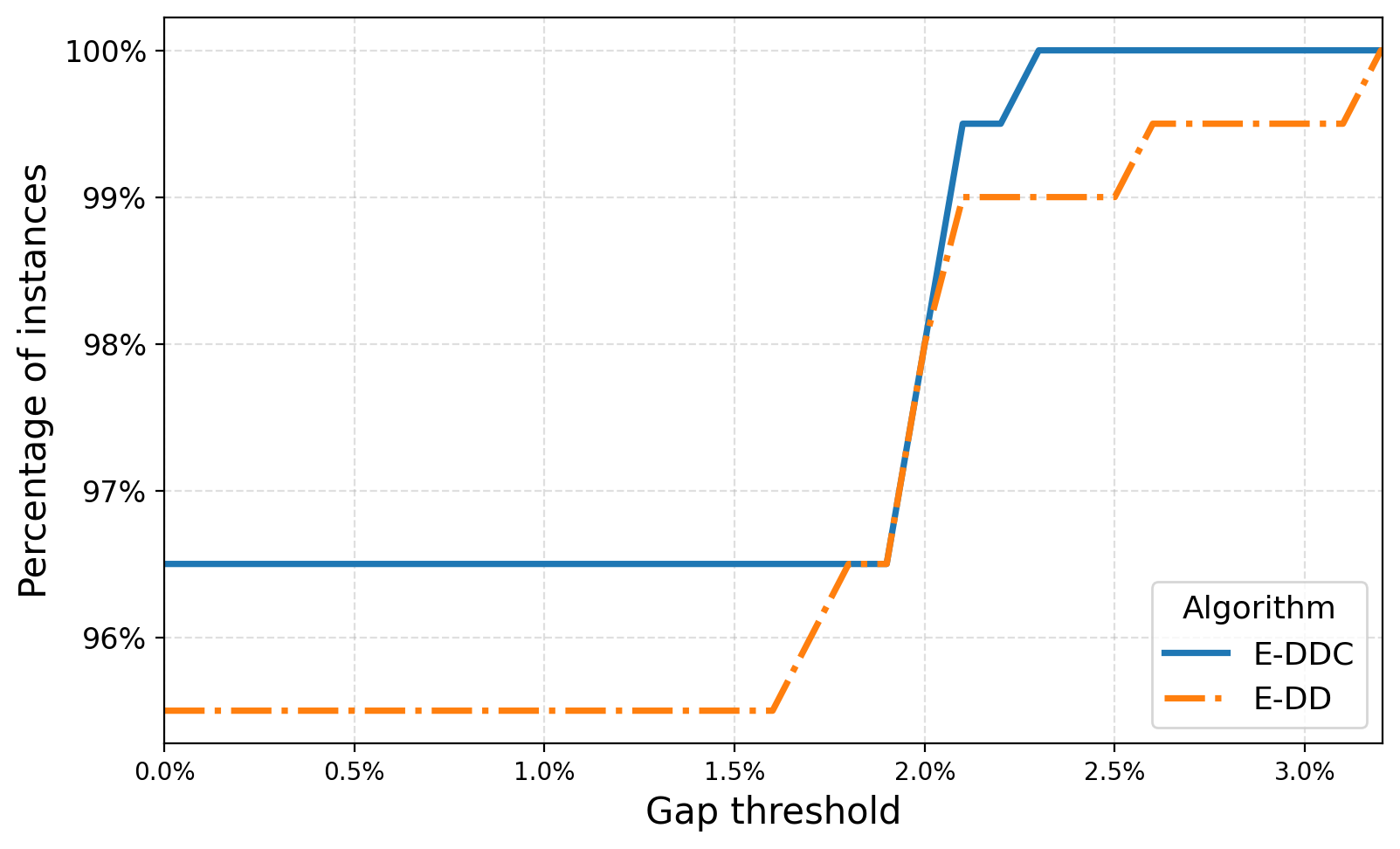}
\caption{\blue{Gaps of exact algorithms.}}
\label{fig:gap_plot}
\end{subfigure}
\caption{\blue{Performance of exact algorithms.}}\label{fig:performance_exact}
\end{figure}

\cref{tab:exact_runtimes} presents the number of instances solved by each exact algorithm and the average runtimes across varying numbers of tasks and drivers for instances with the five mixed-class patterns. 
We first observe that the instances with up to 60 tasks and 30 drivers can be solved quickly by both variants. \blue{For every configuration, E-DDC solves at least as many instances as E-DD and requires a shorter average runtime, so the initial corridor search is beneficial throughout. Its impact is most pronounced for the instances with 90 tasks and 45 drivers, where it solves all instances to optimality.}
Since not all \blue{largest} instances \blue{with 120 tasks} are solved, one must be careful when drawing conclusions from \blue{their} average values only. However, the performance plot given in \cref{fig:exact_perf_prof} shows that E-DDC clearly outperforms E-DD in an instance-by-instance comparison. The figure's horizontal axis represents the performance ratio, calculated by dividing the runtime of each method by the runtime of the best-performing method, while the vertical axis represents the percentage of instances solved within a given performance ratio.
E-DDC solves more than \blue{75}\% of all instances faster than E-DD and achieves a performance ratio of 2 for \blue{more than 90}\% of instances.
\cref{fig:gap_plot} shows the optimality gaps of
instances that were not solved to optimality. We observe that the maximum optimality gaps of E-DD and E-DDC are very small \blue{and that E-DDC dominates E-DD over the entire range of gap thresholds. Moreover, the maximum gap of E-DDC (2.3\%) is clearly smaller than the one of E-DD (3.2\%).}
We conclude that both exact variants of our algorithm perform well and can solve large instances in a reasonable amount of time. Incorporating the corridor search feature \blue{consistently enhances the algorithmic performance.}
Further insights into the performance of different components of our exact algorithms and the impact of the number of tasks, ratio between drivers per task, and the driver classes on the performance can be found in \cref{app:runtime-results}.
These additional results show that the number of tasks and driver percentages are the primary determinants of runtime while the impact of driver class allocation is smaller.


\subsection{Algorithmic performance of heuristic variants}\label{sec:results:heuristic}

In this section, we analyze the performance of the heuristic variants and the contribution of the main algorithmic features, i.e., dominance, detour limit, and corridor search. In this analysis, we consider the four heuristic variants (H-B, H-D, H-DD, and H-DDC) that produce lower and upper bounds as well as variant H-C which focuses on providing a lower bound (only) by rapidly computing a high-quality feasible solution.


\begin{table}
	\caption{Average runtimes in seconds for heuristic variants (50 instances per configuration).}
	\label{tab:heur_performance_comparison}
	\centering
	\blue{
	\begin{tabular}{ccrrrrr}
		\toprule
		Tasks & Drivers & \textbf{H-B} & \textbf{H-D} & \textbf{H-DD} & \textbf{H-DDC} & \textbf{H-C} \\
		\midrule
		30 & 15 & 1959.4 & 3.6 & 0.3 & 0.2 & 0.0 \\
		60 & 30 & - & 772.7 & 7.4 & 3.9 & 0.7 \\
		90 & 45 & - & - & 53.5 & 22.7 & 5.0 \\
		120 & 60 & - & - & 218.0 & 98.1 & 30.5 \\
		\bottomrule
	\end{tabular}
	}
\end{table}

 
\cref{tab:heur_performance_comparison} presents the average runtimes of the considered heuristics across varying instance sizes for instances with the mixed-class driver patterns.
H-B
\blue{terminates within the time limit only for the smallest instances, and even there it requires more than 30 minutes on average,}
as it does not consider the dominance rule or detour limit and, therefore, essentially
needs to
enumerate all feasible bundles in each column generation iteration. The positive effect of adding the dominance rule and (on top of it) the detour limit can be clearly observed from the average runtimes of H-D and H-DD, respectively. While H-D can handle instances up to 60 tasks, H-DD consistently produces feasible solutions within the time limit and its average runtime does not exceed \blue{four} minutes even for the largest instances with 120 tasks. The results obtained for H-DDC show that initial corridor search can further reduce the runtime\blue{, roughly halving the average runtime of H-DD on the instances with 90 and 120 tasks}. Overall, these results clearly demonstrate the beneficial contributions of each algorithmic feature (dominance, detour limit, and corridor search). Finally, we observe that the runtimes of the H-C heuristic are much shorter than those of the other variants.

\begin{figure}
	\centering
	\includegraphics[width=.45\textwidth]{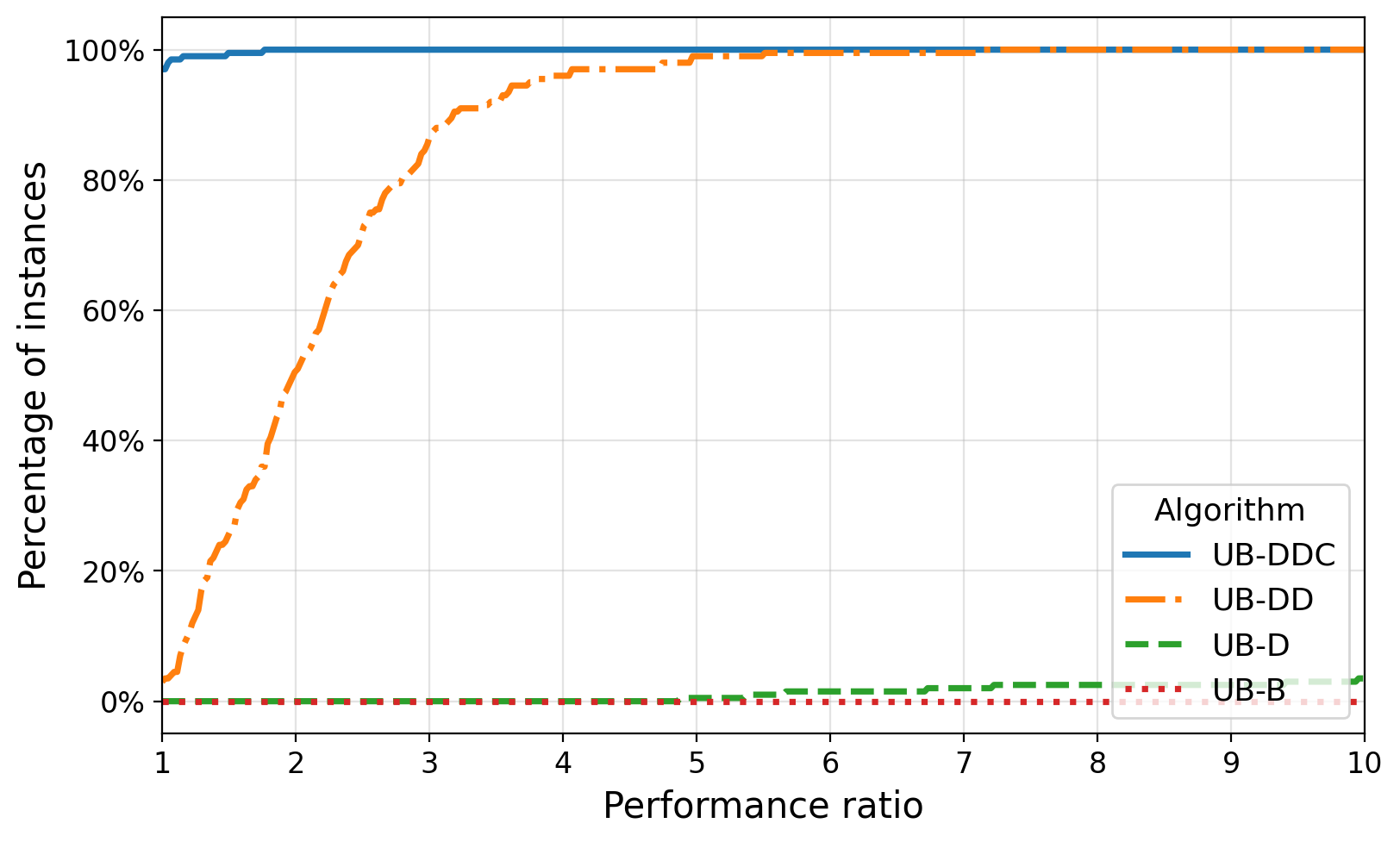}
	\caption{\blue{Performance profiles of the column generation steps of H-B, H-D, H-DD, and H-DDC.}}
	\label{fig:heuristic_perf_prof}
\end{figure}

Further insights for the variants computing upper and lower bounds are obtained from the performance profile of the column generation step given in \cref{fig:heuristic_perf_prof}.
We observe that the column generation step of H-DDC is faster than the one of the other \blue{three} methods in
approximately \blue{97}\% of the instances and takes at most twice as much time than the one of the method for which this step is the fastest on \blue{all} instances. The column generation step of H-DD also has a reasonable performance with a performance ratio of at most four for around \blue{97}\% of the instances\blue{, whereas H-D and H-B perform poorly compared to the others}. 
These results 
confirm the positive contributions of each algorithm feature to reduce the overall runtime.


We now analyze the quality of the solutions produced by the heuristics. \cref{tab:colgen_gaps} shows average gaps ($\mathrm{gap}_\mathrm{H}$) of different heuristics and of their solution values to the optimal solutions ($\mathrm{gap}_\mathrm{OPT}$). 
The reported gap of a heuristic method is the relative gap between its upper bound $\text{UB}_{\text{CG}}$ from the column generation phase (step 1) and the lower bound $\text{LB}_{\text{MIP}}$ obtained from the MILP heuristic (step 2), i.e., $\mathrm{gap}_\mathrm{H}=\frac{\text{UB}_{\text{CG}} - \text{LB}_{\text{MIP}}}{\text{UB}_{\text{CG}}}$.  Since the 
H-C does not compute an upper bound, we report $\mathrm{gap}_\mathrm{H}$  
only for \blue{H-B,} H-D, H-DD, and H-DDC.
\blue{This gap coincides with the root-node optimality gap of a branch-and-price scheme, as the root-node solves the linear relaxation over the whole problem before any branching, and $\mathrm{gap}_\mathrm{H}$ is precisely the resulting gap between the root-node upper and lower bounds. The complete average root-node optimality gaps for H-DDC are reported in \cref{sec:root-node-gaps}.}
Similarly, the gap to the optimal solution $\mathrm{gap}_\mathrm{OPT}= \frac{\text{LB}_{\text{OPT}} - \text{LB}_{\text{MIP}}}{\text{LB}_{\text{OPT}}}$
is the relative gap of a methods lower bound $\text{LB}_{\text{MIP}}$ to the optimal solution $\text{LB}_{\text{OPT}}$ of that instance, which is available in the vast majority of cases as shown in the third column of  \cref{tab:colgen_gaps}. 



\begin{table}
	\centering
	\caption{Average \% gaps of heuristic methods and their gaps to optimal solutions.}
	\label{tab:colgen_gaps}
	\blue{
	\begin{tabular}{@{} r r r r r r r r r r r r @{}}
		\toprule
		\multirow{2}{*}{\textbf{Tasks}} & \multirow{2}{*}{\textbf{Drivers}} & \multirow{2}{*}{\textbf{Solved}} & \multicolumn{4}{c}{\textbf{$\mathrm{gap}_\mathrm{H}$ (\%)}} & \multicolumn{5}{c}{\textbf{$\mathrm{gap}_\mathrm{OPT}$ (\%)}} \\
		\cmidrule(lr){4-7}\cmidrule(lr){8-12}
		& & & \textbf{H-B} & \textbf{H-D} & \textbf{H-DD} & \textbf{H-DDC} & \textbf{H-B} & \textbf{H-D} & \textbf{H-DD} & \textbf{H-DDC} & \textbf{H-C} \\
		\midrule
		30 & 15 & 50 & 1.47 & 1.48 & 1.44 & 1.50 & 0.15 & 0.16 & 0.11 & 0.17 & 10.37 \\
		60 & 30 & 50 & -- & 1.79 & 1.80 & 1.80 & -- & 0.33 & 0.34 & 0.34 & 3.16 \\
		90 & 45 & 50 & -- & -- & 1.65 & 1.68 & -- & -- & 0.34 & 0.38 & 1.39 \\
		120 & 60 & 43 & -- & -- & 1.36 & 1.38 & -- & -- & 0.27 & 0.29 & 0.76 \\
		\bottomrule
	\end{tabular}
	}
\end{table}


We observe that the average \blue{root-node} gaps of all \blue{four} variants producing upper and lower bounds are consistently \blue{at most} 1.8\% when neglecting the cases in which \blue{H-B and} H-D reach 
the time limit. 
This does not only indicate that high-quality solutions are found by the MILP heuristic step for all considered instances sizes, but also shows that the column generation bounds are very tight. The quality of the solutions obtained by H-DD and H-DDC is further highlighted by the extremely small average gaps to the optimal solutions which do not exceed \blue{0.38}\%. We conclude that using H-DD and H-DDC instead of their exact counterparts can be attractive as the loss in solution quality is very small while the required computation time is usually one order of magnitude smaller.


The performance of H-C varies with instance size, whereas the other variants remain stable. Its average gap to the optimal solutions is above 10\% for the smallest instances, but drops below 1\% for the largest instances. 
To better understand this variance, we examine its behavior for different driver to task ratios and compare it to H-DDC.
\cref{tab:corridor_heuristic_gaps} clearly shows that the performance of H-C is highly dependent on number of available tasks and drivers and improves as these numbers increase. H-DDC, however, maintains robust performance whose average gaps are consistently below 0.5\%. This variability in H-C’s performance can be attributed to its reliance on exploring driver-centric corridors. When more tasks and drivers are available, a larger number of bundles can be created within the corridors, leading to better performance. Conversely, with fewer tasks and drivers, the limited number of available bundles restricts its effectiveness. \blue{In contrast, H-DDC constructs bundles by exploring the entire task space rather than only driver-centric corridors, so its optimality gap remains largely unaffected by task or driver density.}

\begin{table}
	\centering
	\caption{Numbers of instances for which an optimal solution is known and average \% gaps to optimal solutions of H-C and H-DDC for different driver to task ratios.} 
	\blue{
	\label{tab:corridor_heuristic_gaps}
	\begin{tabular}{c ccccc ccccc ccccc}
		\toprule
		\multirow{2}{*}{} 
		& \multicolumn{5}{c}{Solved Instances} 
		& \multicolumn{5}{c}{H-C} 
		& \multicolumn{5}{c}{H-DDC} \\
		\cmidrule(lr){2-6} \cmidrule(lr){7-11} \cmidrule(lr){12-16}
		& \multicolumn{5}{c}{Drivers per Task} 
		& \multicolumn{5}{c}{Drivers per Task} 
		& \multicolumn{5}{c}{Drivers per Task} \\
		\cmidrule(lr){2-6} \cmidrule(lr){7-11} \cmidrule(lr){12-16}
		Tasks 
		& 0.1 & 0.2 & 0.3 & 0.4 & 0.5 
		& 0.1 & 0.2 & 0.3 & 0.4 & 0.5 
		& 0.1 & 0.2 & 0.3 & 0.4 & 0.5 \\
		\midrule
		30   & 80 & 80 & 80 & 80 & 80
		& 33.07 & 25.31 & 22.30 & 14.20 &  9.83
		&  0.08 &  0.06 &  0.13 &  0.12 &  0.14 \\
		60   & 80 & 80 & 80 & 80 & 80
		& 15.97 &  9.17 &  5.09 &  3.85 &  2.81
		&  0.07 &  0.34 &  0.35 &  0.31 &  0.29 \\
		90   & 80 & 80 & 80 & 80 & 80
		& 10.65 &  4.20 &  2.73 &  1.63 &  1.27
		&  0.25 &  0.35 &  0.36 &  0.38 &  0.38 \\
		120  & 80 & 72 & 67 & 62 & 65
		&  5.15 &  2.35 &  1.48 &  0.86 &  0.67
		&  0.19 &  0.31 &  0.36 &  0.30 &  0.31 \\
		\bottomrule
	\end{tabular}
	}
\end{table}

Average runtimes of H-C and H-DDC for different numbers of drivers and behavioral classes are given in \cref{tab:sensitivity_performance:heur}. 
 These results confirm the observations made for the exact variant E-DDC, in particular that the main parameters influencing the overall runtime are the numbers of available tasks and drivers. Again, instances can be solved faster when drivers have the lowest sensitivity to compensation (Class 3), instances of Class 2 seem to be slightly easier than those with Class 1 drivers, and 
\blue{runtimes for mixed-class instances fall between those of the individual classes.}

In conclusion, these results show that the algorithmic enhancements (dominance rule, detour limit, and corridor search) proposed in \cref{sec:methodology} are crucial ingredients for the overall success of the proposed algorithms. Consequently, H-DDC incorporating all these components outperforms all other heuristic variants, and consistently generates high-quality solutions for all considered instances. 
%
Finally, H-C can serve as an effective stand-alone heuristic for quickly producing feasible solutions. Its solution quality is, however, highly dependent on the numbers of available tasks and drives. As its performance increases when these values increase, it is an ideal heuristic for tackling very large scale instances for which the runtimes of the other heuristics become prohibitive.

\subsection{Comparison to a sequential approach}

In this section, we compare the benefits of the integrated approach introduced in this article to
a method in which the decisions on bundle generation, compensation determination, and bundle offering are taken sequentially as in \citet{mancini_bundle_2022, macrina2024bundles}. The developed sequential method consists of three main steps. 
In its first step (whose details are given in \apdx) 
a diverse and extensive set of bundles (based on distances between items and their weights) for each depot and driver is generated. In the second step, optimal compensations are computed for each bundle generated in the first step. 
In this step, bundles with an acceptance probability below 50\% are filtered out. 
Similarly, bundles with expected savings less than $0.1$ are discarded.
This step generates on average 5.51 times as many bundles as considered in H-DDC's MILP heuristic step, with the multiplier increasing as the number of tasks grows.
Finally, in the third step a solution is computed by solving the MILP \eqref{eq:milp} using the resulting set of bundles and their optimal compensation values.

\cref{tab:gap_to_best_solution} reports average percentage gaps between the objective values of the sequential method and the best-known integrated solutions, computed as
$\mathrm{gap}_\mathrm{BK} = \frac{\text{LB}_{\text{BK}} - \text{LB}_{\text{Seq}}}{\text{LB}_{\text{BK}}}$.
Here, $\text{LB}_{\text{BK}}$ denotes the objective value from E-DDC for instances solved to optimality, and from H-DDC otherwise, since these provide the best-known values across all instances. $\text{LB}_{\text{Seq}}$ refers to the objective value obtained via the proposed sequential approach. 
We first observe that the integrated approach significantly outperforms the sequential method for all considered settings. Furthermore, these benefits increase with a decreasing number of tasks or drivers per task. These results clearly demonstrate the advantages of jointly optimizing bundle generation, compensation, and assignment decisions, especially when tasks or drivers are scarce.

\begin{table}
	\caption{Average \% gaps between sequential and best-known integrated solutions}
	\label{tab:gap_to_best_solution}
	\centering
	
	\begin{tabular}{c rrrrr}
	\toprule
	& \multicolumn{5}{c}{Drivers per Task}  \\
	\cmidrule(lr){2-6}
	Tasks & 0.1 & 0.2 & 0.3 & 0.4 & 0.5 \\
	\midrule
	30 & 13.83 & 7.86 & 6.30 & 4.82 & 3.38 \\
	60 & 5.53 & 5.60 & 4.07 & 3.35 & 2.93 \\
	90 & 6.66 & 4.26 & 3.21 & 2.55 & 2.25 \\
	120 & 5.55 & 4.08 & 3.55 & 2.50 & 2.24 \\
	\bottomrule
\end{tabular}
\end{table}

\subsection{Sensitivity analysis}\label{sec:sensitivity}

In this section, we analyze the four key criteria (acceptance probability, compensation, detour, and bundle size) of the offers made to occasional drivers using best-known solutions. In the following paragraphs, we examine the sensitivity of these factors to variations in the number of tasks, the driver to task ratio, the behavioral classes of drivers, and the distribution patterns of these classes. The additional results provided in \cref{sec:md} of the electronic companion show that all general trends confirmed in the following subsections remain valid for instances with multiple depots, spatially clustered tasks, and drivers with heterogeneous capacities.

\subsubsection{Sensitivity to numbers of tasks} \label{sec:sensitivity:tasks} 

\begin{figure}
	\captionsetup[subfigure]{aboveskip=-0.5pt,belowskip=-0.5pt}
	\begin{subfigure}[t]{0.49\textwidth}
		\centering
		\includegraphics[width=0.9\textwidth]{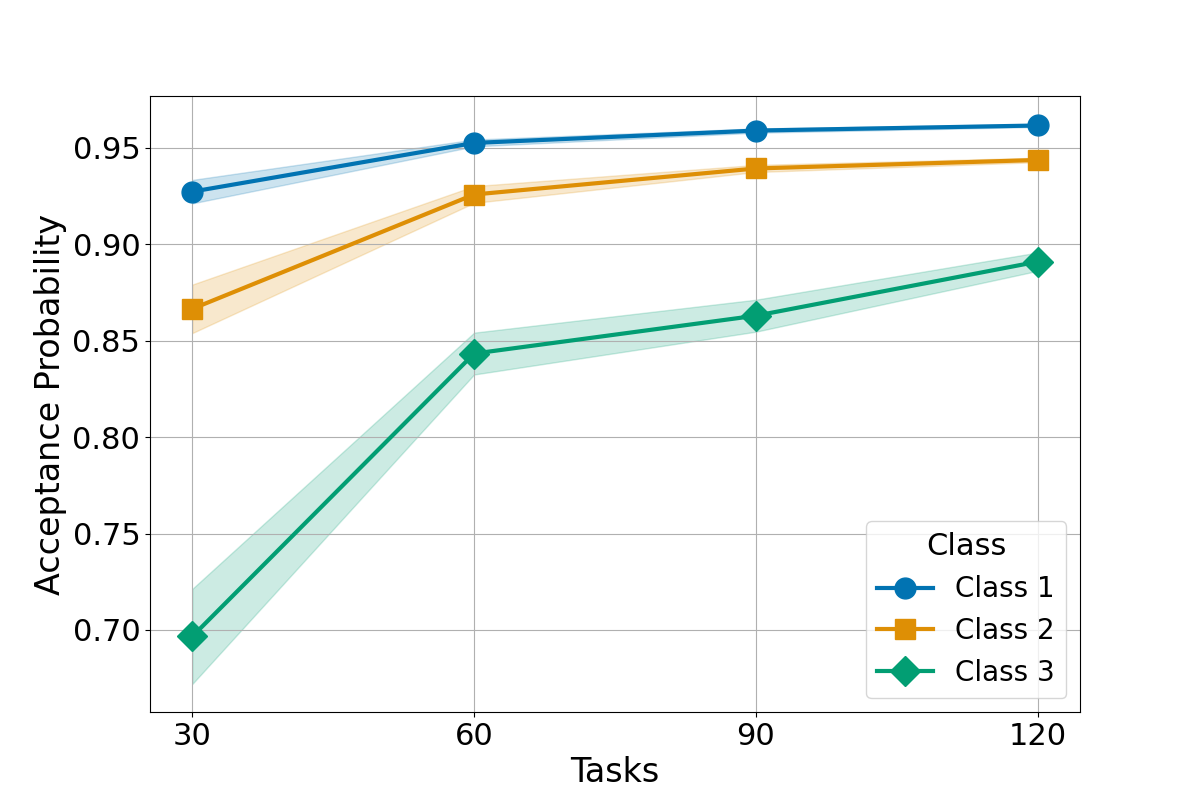}
		\caption{Acceptance probability.} 
	\label{fig:task_acceptance}
\end{subfigure}	\hfill	
\begin{subfigure}[t]{0.49\textwidth}
	\centering
	\includegraphics[width=0.9\textwidth]{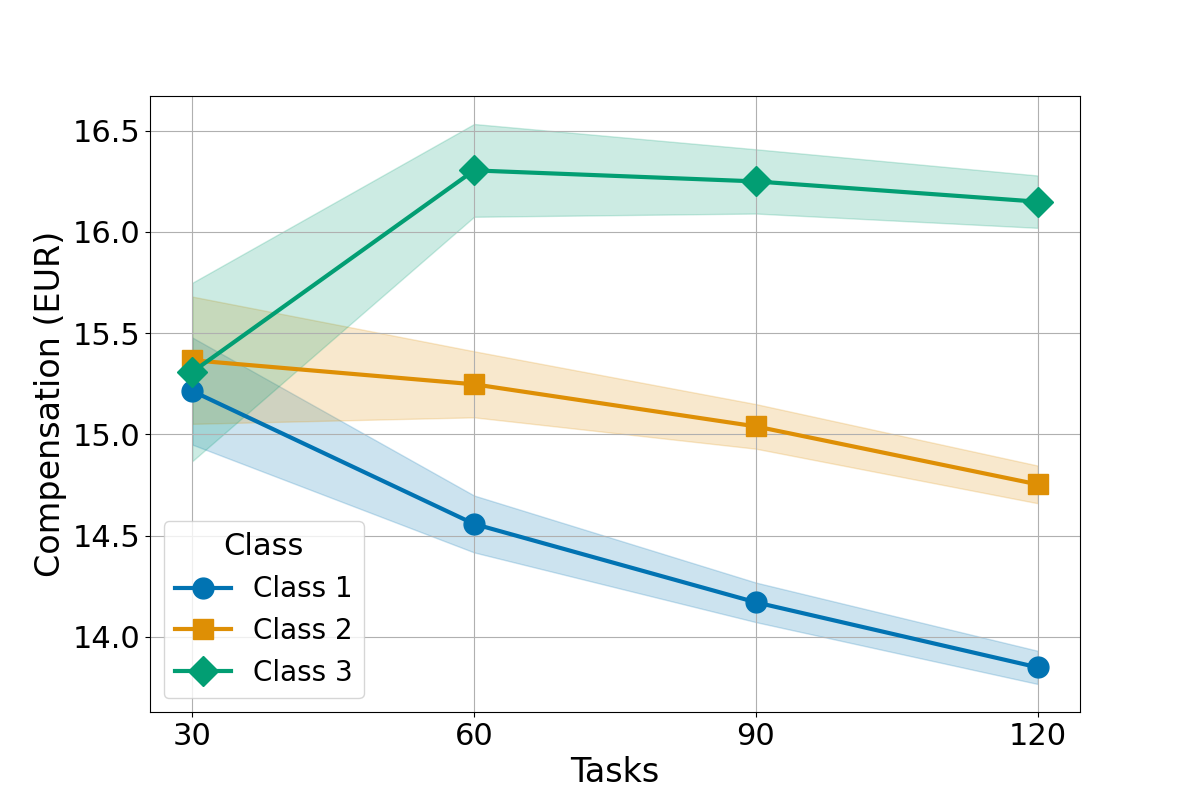}
	\caption{Compensation.} 
\label{fig:task_compensation}
\end{subfigure}
\begin{subfigure}[t]{0.49\textwidth}
	\centering
\includegraphics[width=0.9\textwidth]{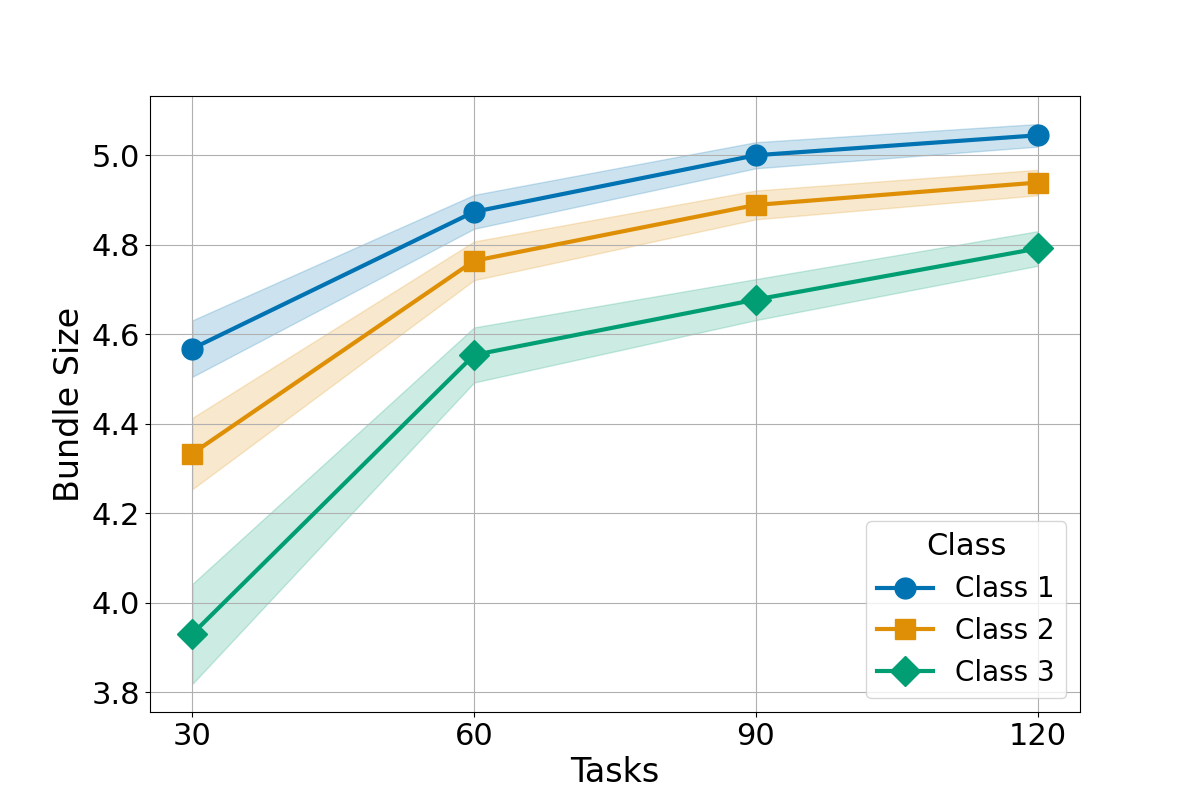}
\caption{Bundle size.} 
\label{fig:task_bundle}
\end{subfigure}	\hfill	
\begin{subfigure}[t]{0.49\textwidth}
	\centering
\includegraphics[width=0.9\textwidth]{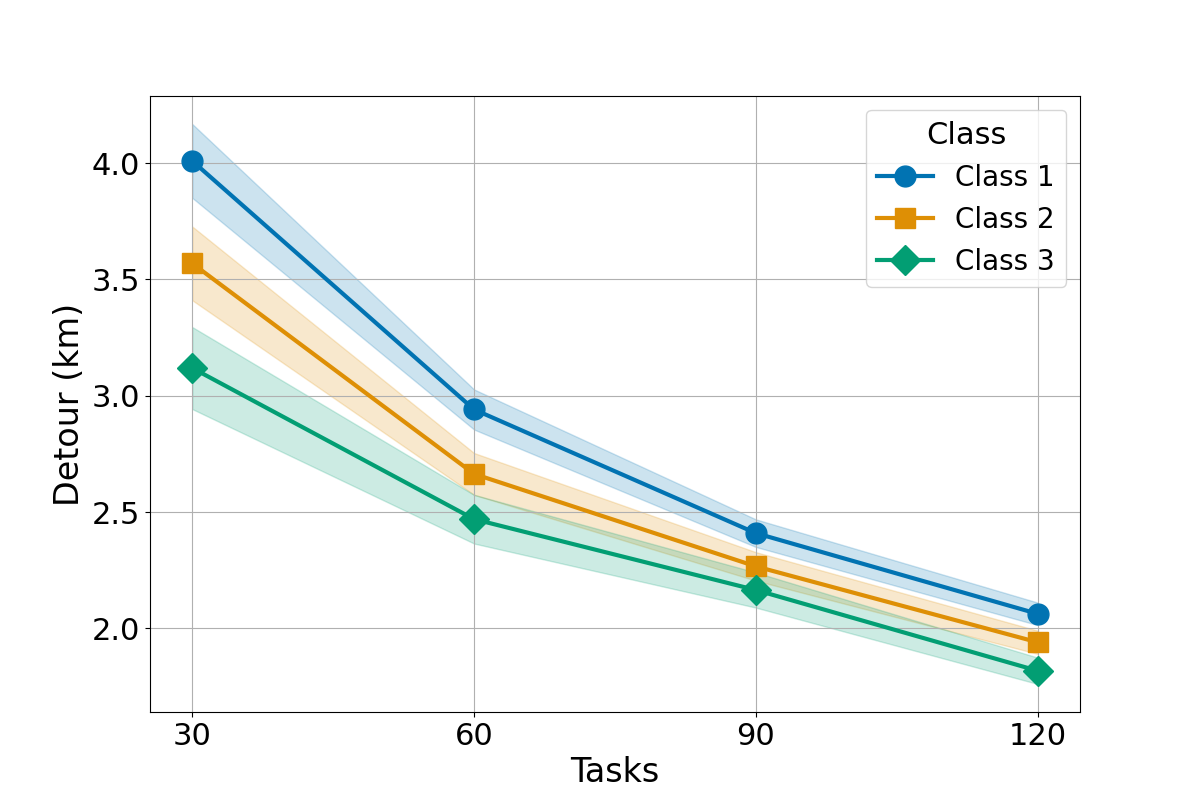}
\caption{Detour.} 
\label{fig:task_detour}
\end{subfigure}
\caption{Mean values and their 95\% confidence intervals of different metrics for different numbers of tasks and behavioral classes.}
\label{fig:task_all_metrics}
\end{figure}

\cref{fig:task_all_metrics} shows the impact of the number of tasks on the four key metrics using mean values and their 95\% confidence intervals. We depict the impact on each of the three considered behavioral classes. As can be seen, for the most part, the trends are similar for the three behavioral classes. We will first discuss this general trend before discussing the difference between class behavior in \cref{sec:results:sensitivity:behavioral}.

We observe that for all driver classes offers with higher acceptance probabilities are made when the number of tasks increases. More tasks results in a larger number of possible bundles some of which may better align with drivers' preferences, thereby contributing to this trend. It can be clearly seen that despite a small increase in bundle size of up to 15\%, the mean detour drops significantly with a larger number of tasks by up to 50\%. As bundles become more attractive, acceptance probabilities increase despite a lower mean compensation. One notable exception occurs for Class 3 drivers when the number of tasks increases from 30 to 60. Here, mean compensation actually needs to increase to also lead to an increasing acceptance probability. This may be an artifact of the general trend stemming from a rather limited number possibilities to create attractive bundles for Class 3 drivers when considering only 30 tasks. Finally, we observe that as the number of tasks increases, the confidence-interval bands narrow across all metrics, indicating that the resulting offers become more stable. Moreover, the bands are consistently the widest for Class~3 and the narrowest for Class~1, suggesting that the relatively homogeneous offers are made to Class~1 drivers, while offers with a broader range of characteristics are made to Class~3 drivers.

\subsubsection{Sensitivity to driver to task ratio} \label{sec:sensitivity:drivers}

\begin{figure}
		\captionsetup[subfigure]{aboveskip=-0.5pt,belowskip=-0.5pt}
	\begin{subfigure}[t]{0.49\textwidth}
		\centering
		\includegraphics[width=0.9\textwidth]{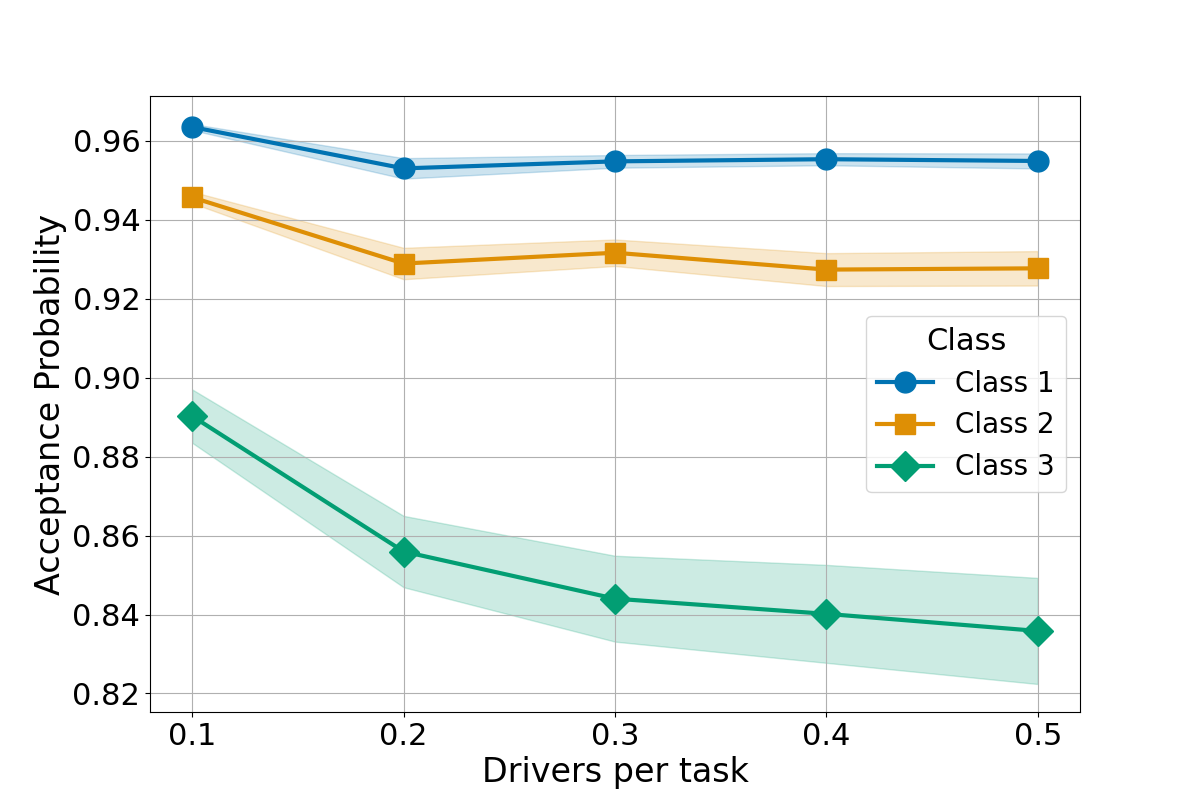}
		\caption{Acceptance probability} 
	\label{fig:driver_acceptance}
\end{subfigure}	\hfill	
\begin{subfigure}[t]{0.49\textwidth}
	\centering
	\includegraphics[width=0.9\textwidth]{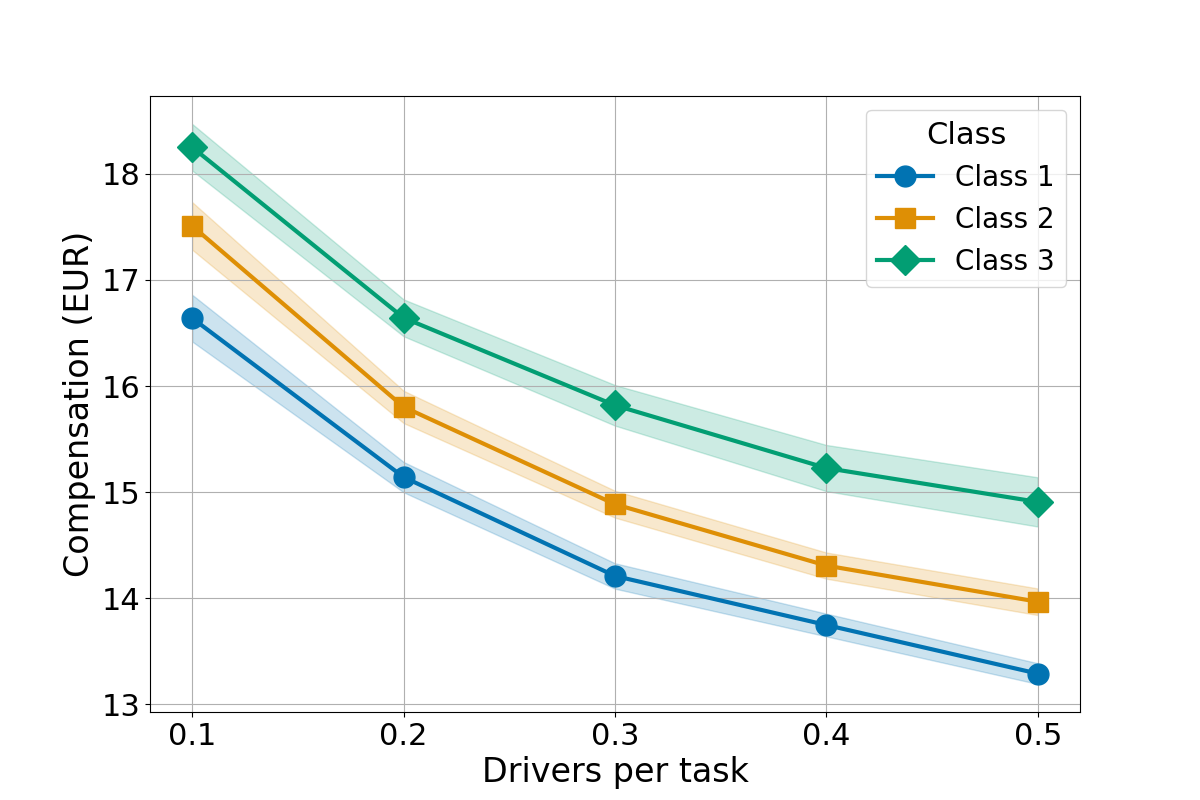}
	\caption{Compensation} 
\label{fig:driver_compensation}
\end{subfigure}
\begin{subfigure}[t]{0.49\textwidth}
	\centering
\includegraphics[width=0.9\textwidth]{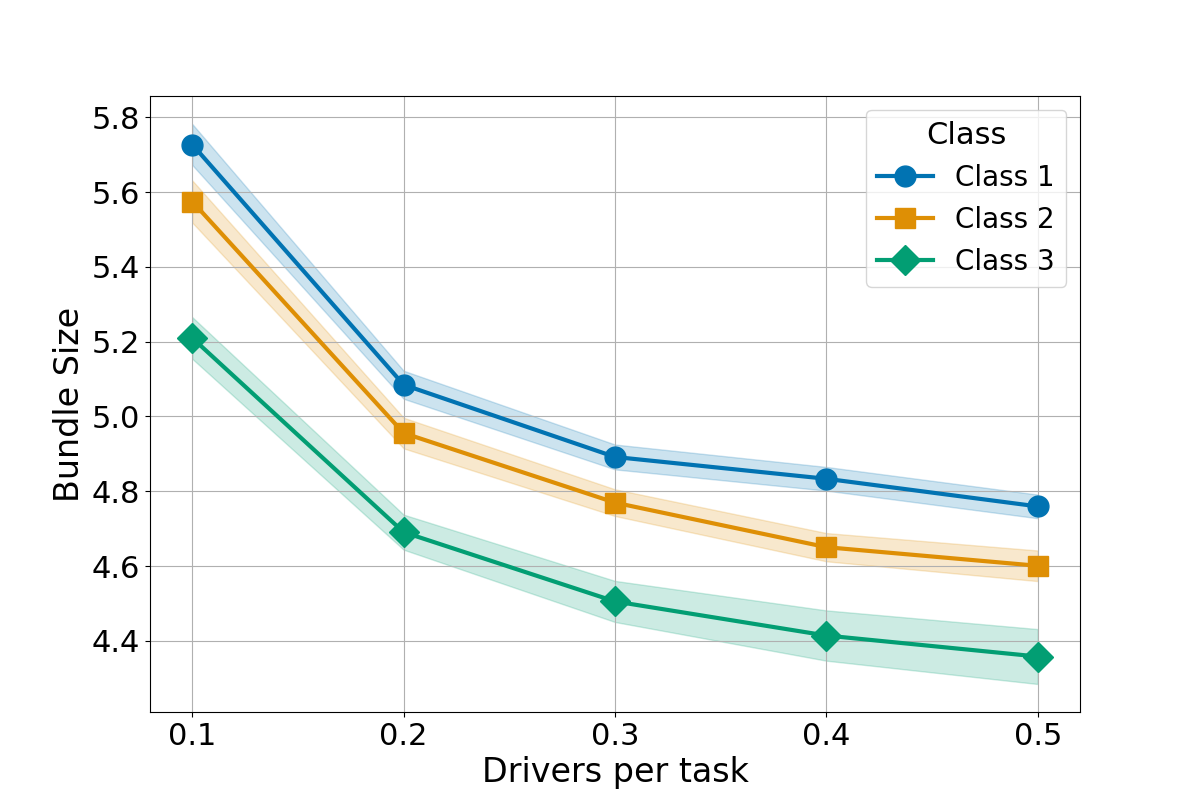}
\caption{Bundle size} 
\label{fig:driver_bundle}
\end{subfigure}	\hfill	
\begin{subfigure}[t]{0.49\textwidth}
	\centering
\includegraphics[width=0.9\textwidth]{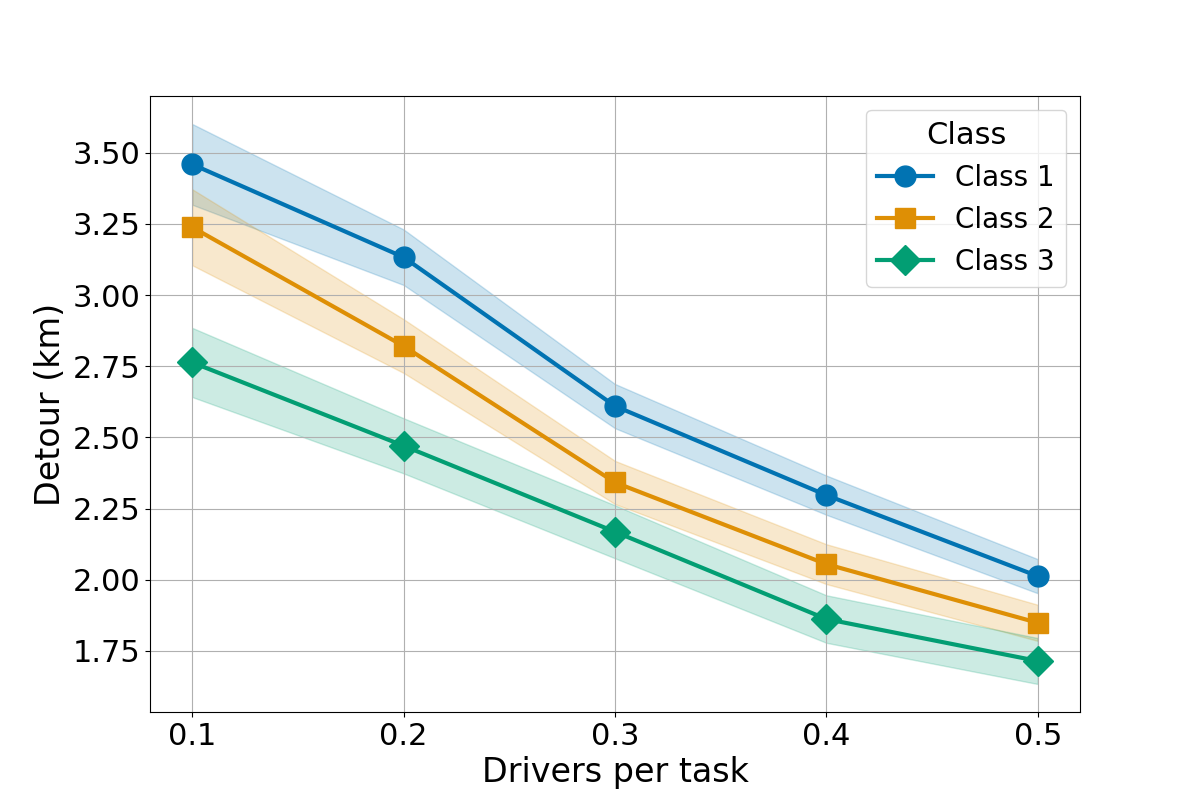}
\caption{Detour} 
\label{fig:driver_detour}
\end{subfigure}

\caption{Mean values and their 95\% confidence intervals of different metrics for different driver to task ratios and behavioral classes.} 
\label{fig:driver_all_metrics}
\end{figure}

\cref{fig:driver_all_metrics} shows that acceptance probabilities remain relatively stable when the driver-to-task ratio exceeds 0.1. 
The only two exceptions are the steeper drop from 0.1 to 0.2, 
and in the case of Class 3. The former may be linked to the median bundle size of 5 which suggests that no efficient distribution of tasks to bundles exist when the driver to task ratio is 0.1.
%
The steeper decrease in the case of Class 3 is mainly caused by the lower base acceptance probability of the class in combination with an overall lower mean compensation paid as more drivers become available, as 
\cref{fig:driver_compensation} shows. Similar to the case above when increasing the number of tasks, a higher  driver to task ratio 
increases the likelihood of finding drivers better suited for the same bundle, allowing the system to offer lower compensation while still securing acceptance. The figures also clearly demonstrate that as the number of tasks and the ratio between drivers and tasks 
increase, offers with lower detours are consistently made to drivers across all classes.  As we will see in \cref{sec:results:sensitivity:behavioral}, in the case of mixed-class instances, bundles can be offered to other classes more easily and for a lower compensation. Thus, Class 3 receives less attractive 
bundles, leading to an even lower acceptance probability.
Finally, we observe that the confidence bands for acceptance probability and (to a lesser extent) bundle size and compensation widen more markedly for Class~3 than for the other classes 
This pattern reinforces that Class~3 is deprioritized and is utilized opportunistically, receiving a more diverse set of offers than the other classes. In contrast, the detour bands for Classes~1 and~2 shrink as the driver-per-task ratio increases, further reflecting this prioritization structure.

\subsubsection{Sensitivity to behavioral classes of drivers}\label{sec:results:sensitivity:behavioral}

The results of the previous sections show that offers with the highest acceptance probabilities are consistently made to Class 1, while Class 3 receives those with the lowest. Despite having the lowest mean acceptance probabilities, Class 3 receives offers with the highest mean compensations, while the exact opposite is observed for Class 1. This contrast highlights the influence of sensitivity to compensation across classes, with Class 1 being more responsive to compensation, while Class 3 requires higher compensation to accept offers. Class 3 receives offers with the smallest bundles and shortest detours, while Class 1 receives those with the largest bundles and longest detours. This contrasts with the behavioral class characteristics presented in \cref{tab:class_coefficients}, where Class 3 has the lowest sensitivity to detour and bundle size, while Class 1 has the highest. This discrepancy suggests that compensation plays the most influential role in offer allocation, as Class 1, despite receiving offers with the highest detours and lowest compensation, still exhibits the highest acceptance probabilities.

	\begin{table}
	\caption{Mean values of different metrics for different behavioral classes.} 
\label{tab:mean_metrics}
\resizebox{\textwidth}{!}{%
	\begin{tabular}{lrrrrrr}
		\toprule
		& \multicolumn{2}{c}{Class 1} & \multicolumn{2}{c}{Class 2} & \multicolumn{2}{c}{Class 3} \\ \cmidrule(lr){2-3} \cmidrule(lr){4-5} \cmidrule(lr){6-7}
		& Mixed-Class & Single-Class & Mixed-Class & Single-Class & Mixed-Class & Single-Class \\
		\midrule
		Acceptance Probability & 0.96 & 0.95 & 0.93 & 0.93 & 0.82 & 0.88 \\
		Bundle Size & 5.01 & 4.82 & 4.82 & 4.84 & 4.40 & 4.82 \\
		Compensation & 14.39 & 13.84 & 15.02 & 14.98 & 15.64 & 16.57 \\
		Detour & 2.59 & 2.41 & 2.37 & 2.31 & 2.26 & 2.13 \\
		\bottomrule
	\end{tabular}
}
\end{table}

\cref{tab:mean_metrics} shows the impact of the class configuration on the four key metrics. For each class, the mean values are compared between mixed-class instances, where the class constitutes one third of the available drivers, and single-class instances, where all drivers are from that class.
The main trends among the classes discussed above are still evident. 
Additionally, we observe that the discrepancy in acceptance probabilities across different classes is smaller in single-class settings but becomes more pronounced in mixed-class instances. This suggests that when drivers exhibit homogeneous behavior, they receive similar offers in terms of acceptance probability, with small influence from inter-class behavioral differences. In contrast, in heterogeneous driver settings, inter-class discrepancies increase. Drivers who are more sensitive to compensation tend to receive offers with higher acceptance probabilities, while those who are less sensitive to compensation are offered bundles with lower acceptance probabilities. This is highlighted by the observation that drivers of Class 1 appears to be preferred for bundle allocation in the mixed-class setting, as compensation and bundle size increase compared to the single-class setting, while they decrease for Classes 2 and 3. Due to the slightly higher compensation in the mixed-class setting, offers to Class 1 exhibit slightly higher acceptance probabilities compared to the single-class setting. Conversely, acceptance probabilities for Class 3 are lower in single-class instances compared to the mixed-class setting, and a similar yet less pronounced trend observed for Class 2.

\subsection{Managerial insights} \label{sec:managerial_insights}

In this section, we investigate the impact of explicitly accounting for driver heterogeneity from both the operator's and the drivers' perspectives. To this end, we benchmark our method of personalized compensation offer design 
(henceforth denoted as POD) against two alternatives that (partially) neglect driver heterogeneity. One-Class Optimization (OSO) ignores heterogeneity by assuming that all drivers belong to a single representative behavioral class and designs offers accordingly. Calibrated Compensations (CC) builds on the bundle assignments of OSO but re-calibrates compensations according to \cref{th:compensation_values_logistic} using the driver’s actual behavioral parameters.

\cref{tab:perf_by_class_calc} reports objective values as a percentage of the POD objective and average results per offer for the previously identified key criteria (acceptance probability, compensation, bundle size, and detour) separately for each driver class (C1, C2, C3). The driver class assumed by OSO is indicated in the first column of \cref{tab:perf_by_class_calc}. \blue{The main takeaway is that POD achieves the highest operator objective, while CC recovers most of this performance through class-specific repricing and OSO can incur substantial losses when its assumed representative class does not match drivers' actual behavior.}

\begin{table}[htbp]
	\centering
	\small
	\caption{Metrics by class and calculation mode}
	\label{tab:perf_by_class_calc}
	\begin{tabular}{llc ccc ccc ccc ccc}
		\toprule
		& & & \multicolumn{3}{c}{Acceptance (\%)} & \multicolumn{3}{c}{Compensation} & \multicolumn{3}{c}{Bundle size} & \multicolumn{3}{c}{Detour} \\
		\cmidrule(lr){4-6}\cmidrule(lr){7-9}\cmidrule(lr){10-12}\cmidrule(lr){13-15}
		Class & Method & Obj. (\%) & C1 & C2 & C3 & C1 & C2 & C3 & C1 & C2 & C3 & C1 & C2 & C3 \\
		\midrule
		Mixed & POD & 100.00 
		& 95.36 & 92.14 & 73.08 
		& 15.09 & 15.45 & 14.16 
		& 5.06 & 4.82 & 3.91
		& 3.12 & 2.76 & 2.19 \\
		\midrule
			\multirow{2}{*}{Class 1} & OSO & 66.11
		& 94.79 & 55.21 & 7.94
		& 14.21 & 14.31 & 14.34
		& 4.80 & 4.82 & 4.80
		& 2.76 & 2.81 & 2.89 \\
		& CC & 96.47
		& 94.79 & 92.29 & 85.41
		& 14.21 & 15.52 & 17.40
		& 4.80 & 4.82 & 4.80
		& 2.76 & 2.81 & 2.89 \\
		\midrule
			\multirow{2}{*}{Class 2} & OSO & 81.24
		& 99.13 & 92.18 & 32.27
		& 15.29 & 15.37 & 15.36
		& 4.80 & 4.82 & 4.79
		& 2.64 & 2.68 & 2.75 \\
		& CC & 96.60
		& 94.81 & 92.18 & 85.57
		& 14.08 & 15.37 & 17.21
		& 4.80 & 4.82 & 4.79
		& 2.64 & 2.68 & 2.75 \\
		\midrule
			\multirow{2}{*}{Class 3} & OSO & 80.71
		& 98.92 & 97.23 & 84.93
		& 16.76 & 16.80 & 16.76
		& 4.76 & 4.77 & 4.73
		& 2.42 & 2.47 & 2.51 \\
		& CC & 96.34
		& 94.39 & 91.20 & 84.93
		& 13.75 & 14.99 & 16.76
		& 4.76 & 4.77 & 4.73
		& 2.42 & 2.47 & 2.51 \\
		\bottomrule
	\end{tabular}

\end{table}

The results shown in \cref{tab:perf_by_class_calc} 
clearly demonstrate the benefits of POD which explicitly accounts for driver heterogeneity. 
It achieves the highest objective values 
while ensuring high acceptance for drivers from the first two classes and maintaining meaningful engagement of drivers from Class~3 at moderate compensation levels. In contrast, OSO is highly sensitive to the assumed representative class, leading to substantial performance losses under misclassification. When OSO is calibrated to Class~1, the objective drops to 66.11\% of POD, driven by low acceptance probabilities for Classes~2 and~3 (55.21\% and 7.94\%). While the results of OSO improve significantly when calibrated to Class~2 or Class~3, POD remains superior.

Even though CC uses the bundling and assignment decisions made by OSO it achieves 96.3-96.6\% of POD’s objective due to re-optimizing compensation based on the recipient’s true driver class. This large difference between CC and OSO highlights the value of class-specific compensation optimization. The remaining gap to POD shows, however, that the potential savings cannot be fully recovered through repricing alone which confirms the value of making bundling, assignment, and pricing decisions in an integrated manner while considering driver classes. The results clearly indicate that integrating personalized compensation offers with bundling decisions significantly increases the operator's profitability.


\begin{table}[htbp]
	\centering
	\small
	\setlength{\tabcolsep}{5pt}
	\renewcommand{\arraystretch}{1.15}
	\caption{\blue{Sensitivity of cross-class utility ranges to the utility-proxy parameters. CC and OSO results are reported by calibration class.}}
	\label{tab:welfare_money_metric_sensitivity}
	\blue{
	\begin{tabular}{ccccccccc}
		\toprule
		& & & \multicolumn{3}{c}{\textbf{CC calibration}} & \multicolumn{3}{c}{\textbf{OSO calibration}} \\
		\cmidrule(lr){4-6}\cmidrule(lr){7-9}
		$\boldsymbol{\lambda}_{\boldsymbol{\Delta}}$ & $\boldsymbol{\lambda}_b$ & \textbf{POD} & \textbf{Class 1} & \textbf{Class 2} & \textbf{Class 3} & \textbf{Class 1} & \textbf{Class 2} & \textbf{Class 3} \\
		\midrule
		0.15 & 0.75 & 2.56 & 1.75 & 1.74 & 1.64 & 8.85 & 7.53 & 1.79 \\
		0.23 & 0.75 & 2.48 & 1.76 & 1.75 & 1.65 & 8.66 & 7.39 & 1.77 \\
		0.30 & 0.75 & 2.42 & 1.77 & 1.76 & 1.66 & 8.49 & 7.27 & 1.75 \\
		\midrule
		0.15 & 1.25 & 1.77 & 1.98 & 1.96 & 1.87 & 6.76 & 5.93 & 1.44 \\
		0.23 & 1.25 & 1.69 & 1.99 & 1.98 & 1.89 & 6.57 & 5.79 & 1.42 \\
		0.30 & 1.25 & 1.63 & 2.00 & 1.99 & 1.90 & 6.40 & 5.67 & 1.40 \\
		\midrule
		0.15 & 2.00 & 0.68 & 2.31 & 2.30 & 2.23 & 3.64 & 3.52 & 0.92 \\
		0.23 & 2.00 & 0.71 & 2.33 & 2.32 & 2.24 & 3.44 & 3.38 & 0.90 \\
		0.30 & 2.00 & 0.74 & 2.34 & 2.33 & 2.25 & 3.28 & 3.26 & 0.88 \\
		\bottomrule
	\end{tabular}
	}
\end{table}

While a better understanding of driver behavior is beneficial for the operator, \blue{we also find that POD generally yields a more balanced welfare distribution as the monetary value assigned to driver effort increases.} To evaluate these benefits, we define a \emph{money-metric utility proxy} that approximates driver welfare as the \emph{expected net benefit} of receiving an offer. For this, we use the class-level averages reported in \cref{tab:perf_by_class_calc} (namely $p_{jk}$, $c_{jk}$, $|s_{jk}|$, and $\Delta_{jk}$ denoting the acceptance probability, compensation, bundle size, and detour for class \(j\) under method \(k\), respectively) to compute the utility proxy $U_{jk}=p_{jk}(c_{jk}-\lambda_{\Delta}\Delta_{jk}-\lambda_b |s_{jk}|)$ for each driver class $j$ and method $k$. Here, the term $(c_{jk} - \lambda_{\Delta}\Delta_{jk} - \lambda_{b}|s_{jk}|)$ represents the net monetary benefit obtained by subtracting the monetary value of effort (detour and bundle-related handling) from the compensation. Multiplying the net benefit by the acceptance probability $p_{jk}$ yields the expected net benefit, since compensation and effort are realized only if the driver accepts the offer. As parameters, we use (i) $\blue{\lambda_{\Delta}\in \{0.15,0.23,0.30\}}$ as a proxy for the cost of each km of detour (consistent with the Dutch benchmark reimbursement rate \blue{\euro$0.23$} per km \citep{businessgovnl_travel_allowance_2026}), and (ii) $\lambda_{b}\in\{0.75,1.25,2.00\}$ as the proxy for the money-equivalent cost for the time of handling a single task (\blue{consistent with three to eight} minutes 
	valued at the Dutch minimum wage of \euro$14.71$ \citep{governmentnl_minimum_wage_amounts_2026}).



\blue{Table~\ref{tab:welfare_money_metric_sensitivity} reports the sensitivity of the cross-class utility ranges to all combinations of $\lambda_{\Delta}$ and $\lambda_b$. When delivery effort is assigned a low monetary value ($\lambda_b=0.75$), CC and OSO calibrated to Class~3 produce smaller ranges than POD. This reveals a trade-off between POD's individualized bundle and compensation design and the more consistent offer structures of the benchmarks. OSO calibrated to Class~3 provides similar offers with consistently high compensation, whereas CC retains the OSO bundle and assignment structure but personalizes compensation according to the recipient's class. These alternatives can produce a more balanced welfare distribution when effort costs are valued less, but they achieve lower platform benefits than POD. Thus, the results reveal a trade-off between platform performance and the distribution of driver welfare rather than uniform dominance by one method. As the monetary value assigned to the delivery effort increases, the POD range generally decreases, although this improvement is driven primarily by the bundle-handling proxy $\lambda_b$, while changes in $\lambda_{\Delta}$ have only a limited (or adverse) effect. At $\lambda_b=1.25$, the POD ranges are consistently lower than CC ranges but still lower than OSO calibrated to Class~3. Finally, POD achieves the smallest cross-class range when $\lambda_b=2.00$, indicating that personalized offer and compensation design not only yields the highest platform benefits but also achieves the most balanced welfare distribution when drivers assign a high value to their effort.}

\section{Conclusions and future work} \label{sec:conclusion}


In this paper, we proposed and analyzed a novel algorithmic solution method for crowdsourced delivery which simultaneously optimizes task bundling, bundle assignment, and compensation decisions while considering bundle- and compensation-dependent individual acceptance probabilities of occasional drivers. 
By \blue{jointly optimizing} these interrelated decisions \blue{within an exact framework}, our work bridges a critical gap in the literature, which previously addressed them \blue{individually, partially or sequentially}.
The key challenge of optimizing compensation decisions under
uncertain (acceptance) behavior 
is highly relevant for the broader gig economy with crowdsourced delivery being one particular example. Indeed, the optimal determination
of worker-specific incentives has been demonstrated as an effective means to address these challenges \citep{allon2023impact}. Thus, our approach offers a strong foundation for diverse applications across the gig economy.

We formulated the considered problem as a mixed-integer nonlinear program with an exponential number of variables to account for all possible bundles, and solved it exactly using a column generation-based algorithm. Our formulation incorporates generic logistic functions that predict occasional drivers' acceptance probabilities based on driver- and bundle-specific predictors. We showed that the pricing subproblem can be solved by identifying an elementary path with resource constraints and maximum reduced costs. As this path problem features a non-linear and non-additive objective function, we developed tailored dominance and reduced cost-based pruning rules since classical dominance and pruning rules rely on the additivity property of the objective function. 
We also developed a corridor-based heuristic that can be used either as a standalone approach or to initialize other algorithms. We defined a logistic acceptance probability function that features some of the most significant predictors for occasional drivers' choice (detour distance, bundle size, and compensation). 
Furthermore, we showed how to align our generic dominance and pruning rules with externally derived logistic functions, including theoretical detour limits for task pruning.

Our computational experiments, conducted on a large instance library varying in task numbers, occasional driver ratios, and driver behavioral classes (with differing levels of effort averseness and compensation sensitivity), included instances with up to 120 tasks and 60 occasional drivers. We evaluated several heuristic and exact algorithm variants that integrated the proposed components, and the results show all of these 
components improve the
performance. In particular, the exact and heuristic variants (E-DDC and H-DDC) that incorporate dominance rules, theoretical detour limits, and corridor-based initialization deliver the best overall performance. 
E-DDC exactly solves all instances with up to 60 tasks and 30 drivers, 
most instances up to 90 tasks and 45 drivers, and a substantial amount of instances up to 120 tasks and 60 drivers. For the instances that cannot be solved to optimality, the resulting optimality gaps remain very small.
H-DDC consistently achieved dual gaps within 2\% and primal gaps within 0.5\%, with an average runtime of 400 seconds on the largest instances. Additionally, a standalone corridor-based heuristic (H-C) demonstrates an average runtime of 80 seconds for the largest instances and delivers higher quality solutions for instances with the highest task numbers and occasional driver ratios. Consequently, these algorithm variants offer practitioners an effective balance between exploration and exploitation, delivering exact solutions when needed while also ensuring fast runtimes for large-scale or time-sensitive scenarios. The results further highlight the benefits of our integrated method over a benchmark method inspired from the literature, where bundle generation, compensation, and assignment decisions are made sequentially. Finally, the sensitivity analysis revealed that the method effectively leverages the distinct characteristics of driver classes to design and offer tailored bundles for each profile.


These findings offer valuable insights into crowdsourced delivery systems characterized by information asymmetry that favors the operator. Future research could focus on addressing the potential disadvantages faced by occasional drivers and promoting fairness in offer allocation. Expanding the model to include full-time gig workers who depend on delivery tasks as their primary income source also offers a valuable avenue for future research. \blue{Such future research may be devoted to multiple pickup locations, preparation and depot-waiting times, strategic interactions among market participants commonly observed on the aforementioned gig economy platforms, and the associated routing, behavioral, and game-theoretic extensions.}

Another promising avenue for future research 
	involves embedding 
	the model and approach studied in this article 
	into a multi-period framework.  
\blue{Our setting could be relevant to such dynamic operational environments in several ways. Since it admits exact solutions, it could serve as a benchmark for multiperiod models that consider all available drivers and delivery tasks. Also, depending on the operator’s decision-making frequency, suitable heuristic or exact variants could be employed as myopic policies that optimize assignments among the tasks and drivers available at each decision epoch. By incorporating estimates of the future value of the tasks and/or drivers, the static problem could balance immediate performance with future opportunities. Additionally, although the current setting allows at most one offer per driver and includes each task in at most one offer, a dynamic implementation could reconsider drivers who reject offers and tasks that are rejected or postponed at subsequent decision epochs. Thus, while the restriction would continue to apply within each decision epoch, drivers and tasks could receive multiple offers over time, effectively relaxing it over the planning horizon.}

\blue{Finally, the effectiveness and realism of such dynamic frameworks could be further enhanced by updating acceptance probability estimates online based on observed driver responses, and by endogenizing customer demand and driver participation to capture how current decisions affect future task arrivals and driver availability. To effectively use the static problem as a decision tool in environments that require faster solution methods, future research could focus on improving computational performance by developing new dominance and pruning rules that further restrict the label search space, as well as more efficient heuristic variants.}


\ifTR

\bibliographystyle{plainnat}

\else

\bibliographystyle{elsarticle-harv}
\fi

\bibliography{ref}







\ifTR

\else

\ifBlind
\ECRRHSecondLine{\bf Authors’ names blinded for peer review}%
\ECLRHSecondLine{\bf Authors’ names blinded for peer review}%
\fi

%

\appendix
\section{Proofs} \label{sec:proofs}


\repTH{\cref{th:compensation_values_logistic}}{\optimalcompensation}
{	
	By inserting the definition of the logistic acceptance probability function \eqref{eq:acceptance_probability} into the formula for optimal compensations derived in \cref{prop:compensation_values} we obtain
\begin{equation*}
	C^*_{wb} = \argmax_{C \ge 0} ~ \frac{\bar{C}-C}{1+e^{-(\alpha^w + B^w X_{wb} + D^w Y_w + \gamma^w C)}}
\end{equation*}
whose first derivative with respect to $C$ is equal to
\begin{equation*}
	\frac{e^{-(\alpha^w + B^w X_{wb} + D^w Y_w + \gamma^w  C)}\left( -e^{(\alpha^w + B^w X_{wb} + D^w Y_w + \gamma^w C)}-1+\gamma^w(\bar{C}-C) \right)}{\left(1+e^{-(\alpha^w + B^w X_{wb} + D^w Y_w + \gamma^w  C)}\right)^2}.
\end{equation*}

The latter is equal to zero if and only if the function 
$$-e^{(\alpha^w + B^w X_{wb} + D^w Y_w + \gamma^w C)}-1+\gamma^w(\bar{C}-C)$$ 
equals zero at $C = C^*_{wb}$. Since $\gamma^w > 0$ (cf.\ \cref{def:logistic_prob}), this function is strictly monotonically decreasing in $C$.
The function is always strictly negative at $C=\bar C$. Also, it is non-positive at $C=0$ when $\gamma_w\bar C \le 1+e^{\alpha^w + B^w X_{wb} + D^w Y_{w}}$. So, there can be at most one root in $(0,\bar C_b)$ if $\gamma_w\bar C > 1+e^{\alpha^w + B^w X_{wb} + D^w Y_{w}}$.

By the following manipulations,
\begin{align*}
	-e^{(\alpha^w + B^w X_{wb} + D^w Y_w + \gamma^w C^*_{wb})}-1+\gamma^w(\bar{C}-C^*_{wb})  & = 0 & \Leftrightarrow \\
	-1+\gamma^w(\bar{C}-C^*_{wb})  & = e^{(\alpha^w + B^w X_{wb} \blue{+} D^w Y_w + \gamma^w C^*_{wb})} & \Leftrightarrow \\
	\left(-1+\gamma^w(\bar{C}-C^*_{wb})\right)e^{-(\alpha^w + B^w X_{wb} + D^w Y_w + \gamma^w  C^*_{wb})} & =1 & \Leftrightarrow \\
	\left(-1+\gamma^w(\bar{C}-C^*_{wb})\right)e^{-1+\gamma^w(\bar{C}-C^*_{wb})} & =e^{\alpha^w + B^w X_{wb} + D^w Y_w + \gamma^w\bar{C}-1}
\end{align*}
we obtain an equation of the form $\eta e^\eta=\phi$, where $\eta=-1+\gamma^w(\bar{C}-C^*_{wb})$ and $\phi=e^{\alpha^w + B^w X_{wb} + D^w Y_w + \gamma^w\bar{C}-1} >0$. Hence, this equation 
can be solved using the principal real branch of Lambert $W$ function $W(\phi)=\eta$, which is single-valued for $\phi>0$. The result follows since
\begin{align*}
&	W(e^{\alpha^w + B^w X_{wb} + D^w Y_w + \gamma^w\bar{C}-1})  = -1+\gamma^w(\bar{C}-C^*_{wb})   \Leftrightarrow \\ 
&	C^*_{wb}  =-\frac{W(e^{\alpha^w + B^w X_{wb} + D^w Y_w + \gamma^w\bar{C}-1})-\gamma^w\bar{C}+1}{\gamma^w}. 
\end{align*}	
}

\repTH{\cref{prop:reduced-cost}}{\reducedcost}{
	The reduced cost is given as $P_w(b,C^*_{wb})(\bar{C}-C^*_{wb})-\sum_{i \in b}\pi_i-\mu_w$ in \cref{eq:sp:obj}. We first use 
	$P_w(b,C^*_{wb})=\frac{1}{1+e^{-(\alpha^w + B^w X_{wb} + D^w Y_w + \gamma^wC^*_{wb})}}$ and $C^*_{wb} = - \frac{W(e^{\alpha^w + B^w X_{wb} + D^w Y_w + \gamma^w \bar{C} - 1})-\gamma^w \bar{C} + 1}{\gamma^w}$ 
	to rewrite the first term $P_w(b,C^*_{wb})(\bar{C}-C^*_{wb})$ as follows:
\begin{align*}
	P_w(b,C^*_{wb})(\bar{C}-C^*_{wb}) &=\frac{\bar{C}-\left( - \frac{W(e^{\alpha^w + B^w X_{wb} + D^w Y_w + \gamma^w \bar{C} - 1})-\gamma^w \bar{C} + 1}{\gamma^w} \right)}{1+e^{-\left(\alpha^w + B^w X_{wb} + D^w Y_w + \gamma^w \left( - \frac{W(e^{\alpha^w + B^w X_{wb} + D^w Y_w + \gamma^w \bar{C} - 1})-\gamma^w \bar{C} + 1}{\gamma^w} \right)\right)}} \\
	&=\frac{\frac{W(e^{\alpha^w + B^w X_{wb} + D^w Y_w + \gamma^w \bar{C} - 1}) + 1}{\gamma^w} }{1+e^{-(\alpha^w + B^w X_{wb} + D^w Y_w + \gamma^w \bar{C} - 1  - W(e^{\alpha^w + B^w X_{wb} + D^w Y_w + \gamma^w \bar{C} - 1}))}} \\
	&=\frac{\frac{W(e^{\alpha^w + B^w X_{wb} + D^w Y_w + \gamma^w \bar{C} - 1}) + 1}{\gamma^w} }{1+e^{-(\alpha^w + B^w X_{wb} + D^w Y_w + \gamma^w \bar{C} - 1)}e^{W(e^{\alpha^w + B^w X_{wb} + D^w Y_w + \gamma^w \bar{C} - 1}))}} 
\end{align*}

Substituting $\phi=e^{\alpha^w + B^w X_{wb} + D^w Y_w + \gamma^w \bar{C} - 1}$ and using $W(\phi)=\eta \Leftrightarrow \eta e^\eta=\phi$ we obtain 
\begin{align*}
	P_w(b,C^*_{wb})(\bar{C}-C^*_{wb}) =\frac{\frac{W(\phi)+1}{\gamma^w} }{1+\phi^{-1}e^{W(\phi)}} = \frac{\frac{\eta+1}{\gamma^w} }{1+\frac{e^{-\eta}}{\eta}e^{\eta}} = \frac{\frac{\eta+1}{\gamma^w} }{\frac{\eta+1}{\eta}} = \frac{\eta}{\gamma^w} 
	= \frac{W(e^{\alpha^w + B^w X_{wb} + D^w Y_w + \gamma^w \bar{C} - 1})}{\gamma^w}. 
\end{align*}

Thus, $\tilde{c}_{wb}=\frac{W(e^{\alpha^w + B^w X_{wb} + D^w Y_w + \gamma^w \bar{C} - 1})}{\gamma^w} -\sum_{i \in b}\pi_i-\mu_w$ which concludes the proof.
}

\repTH{\cref{prop:dominancerules}}{\dominance}{

A label $L^p$ \emph{dominates} label $L^f$ if and only if 
\begin{inparaenum}[(a)]
	\item $\tilde{c}^p \ge \tilde{c}^f$ and 
	\item for every feasible extension $L^{f'}$ of $L^f$ corresponding to bundle $(i^f_1, \dots, i^f_{\ell'}, i^f_{\ell'+1}, \dots i^f_{k})$ the same extension $L^{p'}$ of $L^p$ corresponding to bundle $(i^p_1, \dots, i^p_{\ell}, i^f_{\ell'+1}, \dots i^f_{k})$, $k\ge 1$, is feasible, and
	\item the relation $\tilde{c}^{p'} \ge \tilde{c}^{f'}$ holds for any such extension.
\end{inparaenum}
We will show that the Conditions (i)-(iv) imply (a)-(c).

To show that condition (a) holds, we first recall that $W(\cdot)$ and $e^{(\cdot)}$ are strictly increasing functions (see \cref{prop:reduced-cost}) and the terms $\alpha^w$, $D^wY_w$, and $\mu_w$ are identical for both bundles as they are offered to the same occasional driver. Consequently, Conditions $(i)$ and $(ii)$ imply that $\tilde{c}^p \ge \tilde{c}^f$.

Thus, consider extension $L^{f'}$ of $L^f$ corresponding 
to bundle $(i^f_1, \dots, i^f_{\ell'}, i^f_{\ell'+1}, \dots i^f_{k})$, $k\ge 1$. The extension rules concerning the set of reachable nodes and used capacity detailed in \cref{sec:pricing-algorithm} ensure that if Conditions (iii) and (iv) hold for two labels, they also hold for any extension of them that add the 
same tasks in the same sequence for both of them. Thus, label $L^{p'}$ extending $L^p$ by tasks $i^f_{\ell'+1}, 
\dots i^f_{k}$ is feasible, proving that Condition (b) holds. Further, $R^{p'}\supseteq R^{f'}$ as well as $q^{p'}\le q^{f'}$, i.e., conditions (iii) and (iv) hold for the considered extensions.

To show that condition (c) holds, we prove that that Conditions (i) and (ii) are fulfilled for the extension as well. Repeating the argument made to prove condition (a), we thus also obtain $\tilde{c}^{p'} \ge \tilde{c}^{f'}$ for the extension and the proposition follows.

Condition (ii) holds since $\sum_{i\in b^p} \pi_i \le \sum_{i\in b^f} \pi_i \Rightarrow \sum_{i\in b^{p'}} \pi_i \le \sum_{i\in b^{f'}} \pi_i$ as the same tasks are added to both labels. To see that Condition (i) holds for any extension, we observe that the last task within the label and the newly added task are identical in all extension steps made for $L^f$ and $L^p$ to arrive at $L^{f'}$ and $L^{p'}$, respectively. Thus, by definition of the generic update functions $\mathcal{U}_i$, the increase of all bundle-dependent predictors are identical for all extensions made for $L^f$ and $L^p$. Furthermore, $\bar{C}^{p'} - \bar{C}^{f'} = \bar{C}^{p} - \bar{C}^{f}$ since the predetermined 3PL costs (of the same tasks) are added in each extension step. As a consequence of these two observations, $B^\omega X^p + \gamma^w \bar{C}^p \ge B^w X^f + \gamma^w \bar{C}^f \Rightarrow B^\omega X^{p'} + \gamma^w \bar{C}^{p'} \ge B^w X^{f'}  + \gamma^w \bar{C}^{f'}$, i.e., Condition (i) holds. }

\repTH{\cref{corr:detourupperbound}}{\detourupperbound}{
Applying \cref{corr:rcbasedpruning} for the logistic acceptance probability function~\eqref{eq:acceptance_model}, $\hat{c}^k=\frac{W(e^{\alpha^w + \beta_1^w\Delta^k+\beta_2^w\hat{b}^k + \gamma^w\hat{C}^k - 1})}{\gamma^w} - \hat{\pi}^k - \mu_w$ provides an upper bound on the reduced cost of any extension of label $L_w^p$ by exactly $k$ tasks with a detour of $\Delta^k$. Thereby, 
$\hat{b}^k=|b^p|+k$, $\hat{C}^k=\bar{C}^p+kc^{max}$, and $\hat{\pi}^k=\pi^p+k\pi^{min}$.
%
We first observe that $\hat{c}^k\ge 0$ if and only if $W(e^{\alpha^w + \beta_1^w\Delta^k+\beta_2^w\hat{b}^k + \gamma^w\hat{C}^k - 1}) \geq \gamma^w(\hat{\pi}^k + \mu_w )$. 
%
The latter inequality can be further simplified as follows since $W(e^t) \geq \eta \Leftrightarrow e^t \geq \eta e^\eta$: 
\begin{align*}
	e^{\alpha^w + \beta_1^w\Delta^k+\beta_2^w\hat{b}^k + \gamma^w\hat{C}^k - 1} & \geq \gamma^w(\hat{\pi}^k + \mu_w )e^{\gamma^w(\hat{\pi}^k + \mu_w )} \\
	\alpha^w + \beta_1^w\Delta^k+\beta_2^w\hat{b}^k + \gamma^w\hat{C}^k - 1 & \geq \ln(\gamma^w(\hat{\pi}^k + \mu_w )) + \gamma^w(\hat{\pi}^k + \mu_w )
\end{align*}

Since $\beta_1^w\le 0$, we obtain the upper bound of 
\begin{align*}
	\Delta^k \leq \frac{\ln(\gamma^w(\hat{\pi}^k + \mu_w )) + \gamma^w(\hat{\pi}^k + \mu_w - \hat{C}^k) - \alpha^w - \beta_2^w (\hat{b}^k) + 1 }{\beta_1^w} = \bar{\Delta}^k
\end{align*}
on $\Delta^k$ for which the reduced costs can be non-negative.
}

\section{Literature overview} \label{sec:literature_apdx}

\cref{tab:literature} provides an  
	overview of the related literature indicating which papers consider acceptance probability, compensation and bundling. If an article considers more than one of these dimensions, the table also highlights how they are integrated. Specifically, sequential indicates that decisions are made consecutively, whereas joint indicates that they are determined simultaneously. In this literature review, we further discuss these studies and highlight
	their key assumptions and limitations compared to our approach.

\begin{table}[H]
	\centering
			\caption{Comparison of the related literature.}
		\label{tab:literature}
		\begin{tabular}{lcccc}
			\toprule
			Study                        & Acceptance Probability & Compensation &   Bundling   & Integration \\ \midrule
			Gdowska et al. (2018)        &      $\checkmark$      &              &              &             \\
			Santini et al. (2022)        &      $\checkmark$      &              &              &             \\
			Ausseil et al. (2024)        &      $\checkmark$      &              &              &             \\
			Barbosa et al. (2023)        &      $\checkmark$      & $\checkmark$ &              & Sequential  \\
			Hou et al. (2022)            &      $\checkmark$      & $\checkmark$ &              & Sequential  \\
			Çınar et al. (2024)          &      $\checkmark$      & $\checkmark$ &              &    Joint    \\
			Kafle et al. (2017)          &                        &              & $\checkmark$ &             \\
			Mancini and Gansterer (2022) &                        &              & $\checkmark$ &             \\
			Mancini and Gansterer (2024) &                        &              & $\checkmark$ &             \\
			Mancini et al. (2025)        &                        &              & $\checkmark$ &             \\
			Wang et al. (2023)           &                        &              & $\checkmark$ &             \\
			Torres et al. (2022b)        &                        &              & $\checkmark$ &             \\
			Le et al. (2021)             &                        & $\checkmark$ &              &             \\
			Horner et al. (2021)         &      $\checkmark$      &              &              &             \\
			Behrendt et al. (2024)       &      $\checkmark$      & $\checkmark$ &              &    Joint    \\
			Torres et al. (2022a)        &      $\checkmark$      &              & $\checkmark$ &    Joint    \\
			Cerulli et al. (2024)        &                        & $\checkmark$ & $\checkmark$ &    Joint    \\
			Macrina et al. (2024)        &      $\checkmark$      & $\checkmark$ & $\checkmark$ & Sequential  \\ \midrule
			Our work                     &      $\checkmark$      & $\checkmark$ & $\checkmark$ &    Joint    \\ \bottomrule
		\end{tabular}
\end{table}

\section{Model fitting} \label{sec:model_fitting}

	Crowdsourced delivery operators can collect and process a wide range of operational data. Table~\ref{tab:offers_schema} presents an illustrative database for logging historical offers, listing the recorded fields, their data types, and their meanings (with additional details on drivers, tasks, and depots available in separate entity-specific databases). The database can include identifiers for the offer, the targeted driver, the pickup depot, and the ordered sequence of task identifiers in the offered bundle, as well as key decision attributes such as the offered compensation and the required detour. It also records the offer outcome (accepted or rejected) and may be augmented with contextual information relevant to the operating environment, such as timestamps, traffic and weather conditions. 

\begin{table}[htbp]
	\centering	
		\caption{Example offer database.}
		\label{tab:offers_schema}
		\begin{tabular}{llp{8cm}}
			\toprule
			\textbf{Field} & \textbf{Type} & \textbf{Meaning} \\
			\midrule
			offer\_id & int & Unique offer identifier \\
			driver\_id & int & Target driver \\
			pickup\_depot\_id & int & Pickup depot \\
			ordered\_task\_ids & array & Ordered sequence of task ids \\
			compensation & float & Offered payment \\
			detour & float & Total additional distance/time vs. baseline \\
			status & enum &  accepted / rejected  \\
			\addlinespace
			context\_\textit{attr} & various & Optional contextual attributes \\
			\bottomrule
		\end{tabular}
\end{table}

	Such a database is useful for fitting models that predict the acceptance behavior of occasional drivers. In this paper, we consider an estimation approach, where drivers are initially clustered based on their historical offer-response patterns (e.g., sensitivity to offer-related factors) before separate logistic regression models are fitted for each cluster to estimate class-specific acceptance probabilities. This two-stage approach is attractive because it yields interpretable behavioral parameters for each segment and facilitates comparisons across segments. 
	Alternatives include latent-class formulations, which jointly infer unobserved behavioral segments and estimate segment-specific acceptance functions, as illustrated by \citet{mohri_modeling_2024}. However, their probabilistic class assignment and estimation structure make cross-class comparisons less interpretable than our approach.
	
	Concerning the use of probit models over logit models, we observe that they often yield similar fitted probabilities (see Section 2.3 of \cite{cramer2003}) and that their use would complicate the analytical developments made in the remainder of this paper. 
	Our choice for using a logistic regression model is therefore based on (i) the prevalence of logistic regression/binary logit models in transportation research, 
	(ii) their analytical tractability, and (iii) the empirical similarity of fitted probabilities for logit and probit models.

\section{Example illustrating label initialization and extension} \label{sec:label_ex}

	
	Consider one driver $w$ with capacity $Q_w=10$ starting from depot $d$ and three tasks $\{1,2,3\}$ that can be delivered from depot $d$, i.e., $T(d)=\{1,2,3\}$. 
	Task loads and 3PL costs are given as $(q_1,c_1)=(4,6)$, $(q_2,c_2)=(5,5)$, and $(q_3,c_3)=(6,4)$ and we consider two additive bundle-dependent predictors $X^p=(\Delta^p,\ |b^p|)$, where $\Delta^p$ is the detour and $|b^p|$ the bundle size. In the following, we will detail the label initiations and the main components of the label extensions resulting from first adding task 1 and then task 2. Thereby, we assume that the detour of task $1$ is equal to $0.8$ while the additional detour of visiting task 2 after task 1 is equal to 0.7, i.e., $U_\Delta(d,1)=0.8$ and $U_\Delta(1,2)=0.7$. 
	%
	In the \emph{initialization} step a label $L_w^p=(b^p,d,q^p,X^p,\bar{C}^p,R^p,-\mu_w)$ for path $p=(s_w=d,d,e_w)$ with an empty set of tasks $b^p=\emptyset$, $q^p=0$, $X^p=(0,0)$, $\bar{C}^p=0$, and $R^p=\{1,2,3\}$ is generated and $\mathcal{T}^p=d$ since $b^p=\emptyset$.
	%
	The \emph{first extension including task $1$} yields a label $L^f_w=(b^f,d,0,X^f,\bar{C}^f,R^f,\tilde{c}^f)$ for path $f=(s_w=d,d,1,e_w)$ with $b^f=(1)$, $q^f=0+q_1=4$, $\bar C^f=0+c_1=6$, $X^f=(0,0)+(0.8,1)=(0.8,1)$, and $R^f=R^p\setminus(\{1\}\cup\{i\in R^p: q^f+q_i>Q_w\})=\{2,3\}$ for which $\mathcal{T}^f=1$ holds. Further \emph{extension by task 2} results in label $L^g=(b^g,d,0,X^g,\bar{C}^g,R^g,\tilde{c}^g)$ for path $g=(s_w=d,d,1,2,e_w)$ with $b^{g}=(1,2)$, $q^{g}=4+q_2=9$, $\bar C^{g}=6+c_2=11$,
	$X^{g}=(0.8,1)+(0.7,1)=(1.5,2)$, and
	$R^{g}=R^f\setminus(\{2\}\cup\{i\in R^f: q^g+q_i>Q_w\})=\emptyset$
	since $9+q_3>10$.

\section{Bundle generation the in sequential method}\label{app:sequential_bundles} 

We generate candidate bundles using a heuristic algorithm that 
that executes the following steps for each depot $d\in D$ and driver $w\in N$.
The algorithm first selects a set of initial (seed) tasks which includes all tasks for which the detour or driver $w\in N$ is at most $\Delta_{\mathrm{max}}$ which is set to $\Delta_{\mathrm{max}} = 5$, as 91.51\% of the bundles in all best-known solutions have a detour below this threshold. 
Parameter $R$ (which is set to $R=250$ to ensure sufficient diversity and coverage) defined the number of (randomized) extensions made for each initial seed task. In each of these $R$ extensions, a new bundle initially containing only the respective seed task is created. Then, a candidate task pool is computed that consists of all other tasks that (i) are reachable from the considered depot, (ii) do not exceed the driver’s remaining capacity, and (iii) can be added to the bundle without exceeding the detour limit $\Delta_{\mathrm{max}}$. The tasks in this pool are in non-decreasing order based on their additional detour, and one of the top $K$ (which we set to $K=5$) candidates is randomly selected to extend the bundle. This extension process continues until no more tasks can be added to the bundle, i.e., until no more candidate tasks remain, the capacity or detour limit is reached. All partial bundles generated during one of the $R$ runs are collected and, finally, only the bundle with the lowest detour is retained among those that contain the same set of tasks.

\section{\blue{Impact of the half-angle $\theta$ and pricing cutoff $\kappa$ on E-DDC}}\label{app:theta-k}
\blue{We assess the joint effect of the corridor half-angle $\theta$ and the pricing cutoff $\kappa$.
We test $\theta\in\{30^\circ,36^\circ,45^\circ,60^\circ\}$ and $\kappa\in\{50,100,200\}$ on the 50 largest instances with 120 tasks and 60 drivers. \cref{tab:theta-k-eddc} reports the average runtime of E-DDC and the number of instances solved to optimality for each parameter pair. 
}

\begin{table}[h]
	\vspace{-1em}
	\centering
	\caption{\blue{Average E-DDC runtimes in seconds and numbers of instances solved to optimality (out of 50) for each combination of the corridor half-angle $\theta$ and pricing cutoff $\kappa$.}}
	\label{tab:theta-k-eddc}
	\blue{
	\begin{tabular}{crrrrrr}
		\toprule
		\multirow{2}{*}{$\theta$}
		& \multicolumn{2}{c}{$\kappa=50$} & \multicolumn{2}{c}{$\kappa=100$} & \multicolumn{2}{c}{$\kappa=200$} \\
		\cmidrule(lr){2-3}\cmidrule(lr){4-5}\cmidrule(lr){6-7}
		& Time & Solved & Time & Solved & Time & Solved \\
		\midrule
		$30^\circ$ & 3469.8 & 41 & 3463.3 & 40 & 3170.6 & 43 \\
		$36^\circ$ & 3328.8 & 41 & 3142.1 & 43 & 3551.0 & 39 \\
		$45^\circ$ & 3444.8 & 40 & 3405.3 & 42 & 3571.4 & 41 \\
		$60^\circ$ & 3228.6 & 43 & 3297.4 & 42 & 3031.0 & 42 \\
		\bottomrule
	\end{tabular}
	}
		\vspace{-1em}
\end{table}

\blue{The combination $\theta=36^\circ$ and $\kappa=100$ solves 43 of the 50 instances to optimality, tying with $(\theta,\kappa)=(60^\circ,50)$ and $(30^\circ,200)$ for the highest number of exact solutions. Among these three combinations, it has the lowest average runtime (3142.1 seconds, compared with 3228.6 and 3170.6 seconds, respectively). We therefore prioritize exact-solve coverage and use average runtime as a tie-breaker, which supports the choice $\theta=36^\circ$ and $\kappa=100$ used in the computational experiments. While this combination yields the best overall performance, the results are robust across all parameter combinations. }

\blue{
\section{Root-node optimality gaps} \label{sec:root-node-gaps}

\cref{tab:colgen-optimality-gap} reports the average root-node optimality gap of H-DDC, broken down by the number of tasks and the driver-to-task ratio. As discussed in \cref{sec:results}, this gap coincides with the column generation gap $\mathrm{gap}_\mathrm{H}$, since the root node solves the linear relaxation over the whole problem before any branching. The gaps remain quite small across the entire grid with individual cells range from 0.93\% to 1.94\%, and the overall average is only 1.52\%. Neither the number of tasks nor the driver-to-task ratio has a consistent effect on the gap, which confirms that H-DDC consistently attains tight root-node bounds regardless of instance size or driver density.
}
\begin{table}[h]
\vspace{-1em}
\centering
\caption{\blue{Average root-node optimality gap (\%) by number of tasks and driver-to-task ratio.}}
\label{tab:colgen-optimality-gap}
\blue{
\begin{tabular}{lrrrrrr}
\hline
\multirow{2}{*}{Tasks} & \multicolumn{5}{c}{Drivers per Task} & \multirow{2}{*}{Average} \\
\cline{2-6}
 & 0.1 & 0.2 & 0.3 & 0.4 & 0.5 & \\
\hline
30 & 1.16 & 1.65 & 1.53 & 1.41 & 1.31 & 1.41 \\
60 & 0.93 & 1.86 & 1.94 & 1.87 & 1.75 & 1.67 \\
90 & 1.17 & 1.63 & 1.71 & 1.74 & 1.71 & 1.59 \\
120 & 1.02 & 1.53 & 1.59 & 1.47 & 1.42 & 1.41 \\
\hline
\textbf{Average} & \textbf{1.07} & \textbf{1.67} & \textbf{1.69} & \textbf{1.62} & \textbf{1.55} & \textbf{1.52} \\
\hline
\end{tabular}
}
\vspace{-2em}
\end{table}

\section{Detailed runtime analysis}\label{app:runtime-results}

\subsection{Runtime consumption per component of exact algorithms} \cref{tab:runtime_components} reports the relative runtime consumption per components of E-DD and E-DDC. For larger instances with 90 and 120 tasks, these results show that time is mostly spent in the variable enumeration \& re-optimization step. This indicates that the final step of our exact algorithm is the main bottleneck when solving large instances. The main reason for this behavior likely is that the dominance and pruning rules introduced in \cref{sec:pricing-algorithm} cannot be applied during enumeration as some labels that are dominated under the usual criteria with positive reduced costs may correspond to bundles used in an optimal solution.


\begin{table}[h]
	\centering
	\caption{Relative runtime consumption (in \%) per components of E-DD and E-DDC.}
	\label{tab:runtime_components}
	\resizebox{\linewidth}{!}{%
	\blue{
		\begin{tabular}{cccccccc}
			\toprule
			& \multicolumn{3}{c}{E-DD} & \multicolumn{4}{c}{E-DDC} \\ \cmidrule(lr){2-4} \cmidrule(lr){5-8}
			Tasks & \makecell{Column\\Generation} & \makecell{MILP\\Heuristic} & \makecell{Enumeration \&\\Re-optimization} & \makecell{Corridor\\Initialization} & \makecell{Column\\Generation} & \makecell{MILP\\Heuristic} & \makecell{Enumeration \&\\Re-optimization} \\ \midrule
			30 & 57.13 & 20.23 & 22.64 & 11.00 & 49.93 & 15.15 & 23.93 \\
			60 & 56.99 & 11.61 & 31.41 & 10.00 & 35.86 & 9.12 & 45.02 \\
			90 & 36.92 & 3.77 & 59.31 & 5.06 & 14.64 & 3.36 & 76.94 \\
			120 & 18.97 & 2.06 & 78.97 & 3.78 & 5.11 & 1.69 & 89.42 \\
			\bottomrule
		\end{tabular}
	}
	}
\end{table}

\subsection{Impact of numbers of tasks and drivers per task on E-DDC}
Further insights into the performance of E-DDC can be obtained from \cref{fig:runtime_profiles_by_M} which shows the percentage of instances solved (y-axis) by the runtimes indicated on x-axis across varying task numbers $|M|$ and drivers per task ratios. For $|M|=30$, the exact algorithm solves all instances within a few seconds for all drivers per task ratios, with $0.1$ being the fastest and higher values only mildly increasing runtimes. For $|M|=60$, runtimes increase to the tens-of-seconds range, yet the curves still reach 100\% for all drivers per task ratios, again with $0.1$ clearly fastest. For $|M|=90$, the distribution spreads into the minutes-to-hours, while \blue{all} ratios 
still achieve complete coverage.
Finally, for $|M|=120$, $0.1$ still reaches 100\% quickly, but for the ratios above $0.1$ the curves level off below 100\% by the 3-hour limit, showing that the exact algorithm solves only a subset of instances at this scale.

\begin{figure}[H]
	\centering
	\begin{subfigure}[b]{0.48\textwidth}
		\centering
		\includegraphics[width=\linewidth]{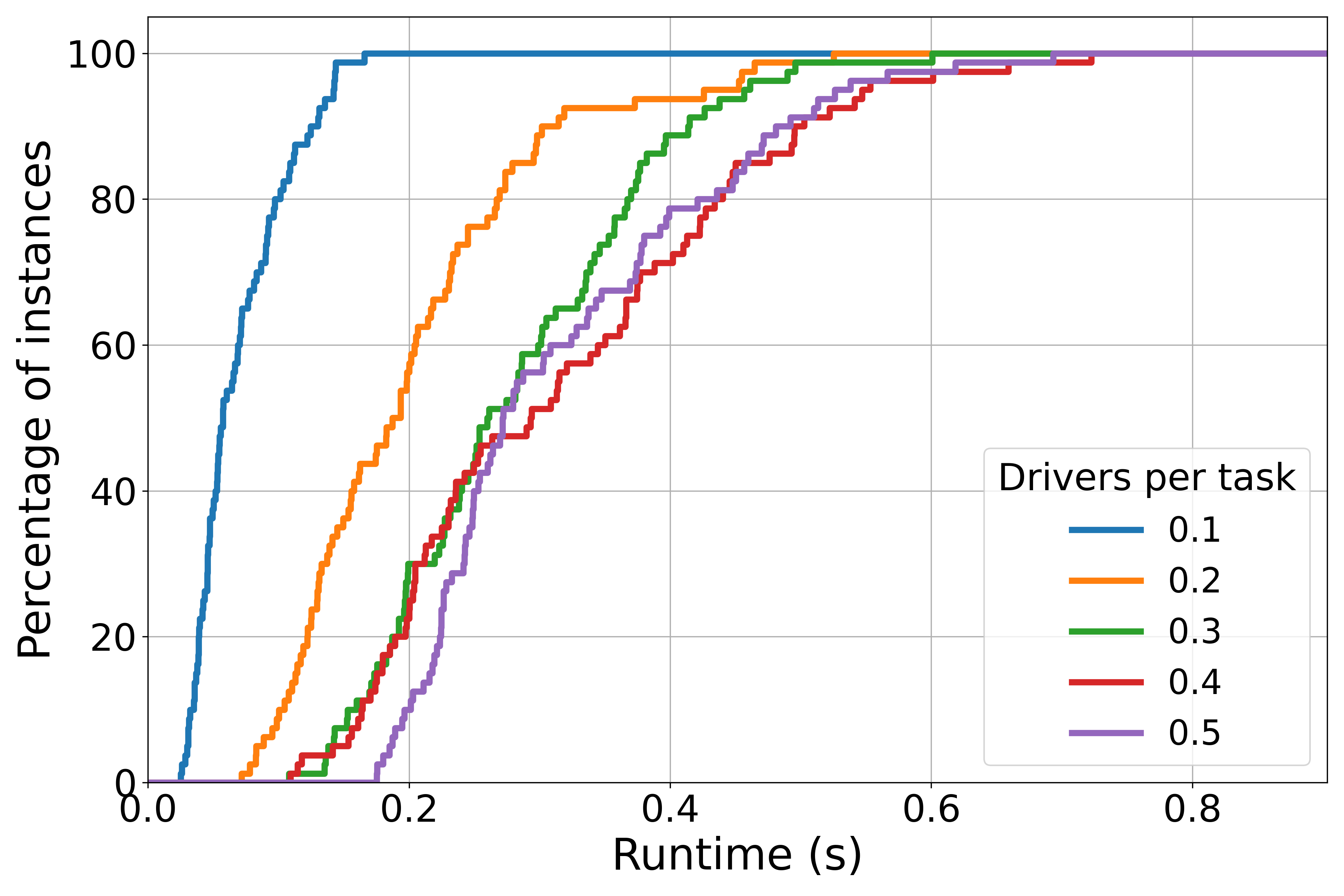}
		\caption{$|M|=30$}
		\label{fig:runtime_profile_m30}
	\end{subfigure}\hfill
	\begin{subfigure}[b]{0.48\textwidth}
		\centering
		\includegraphics[width=\linewidth]{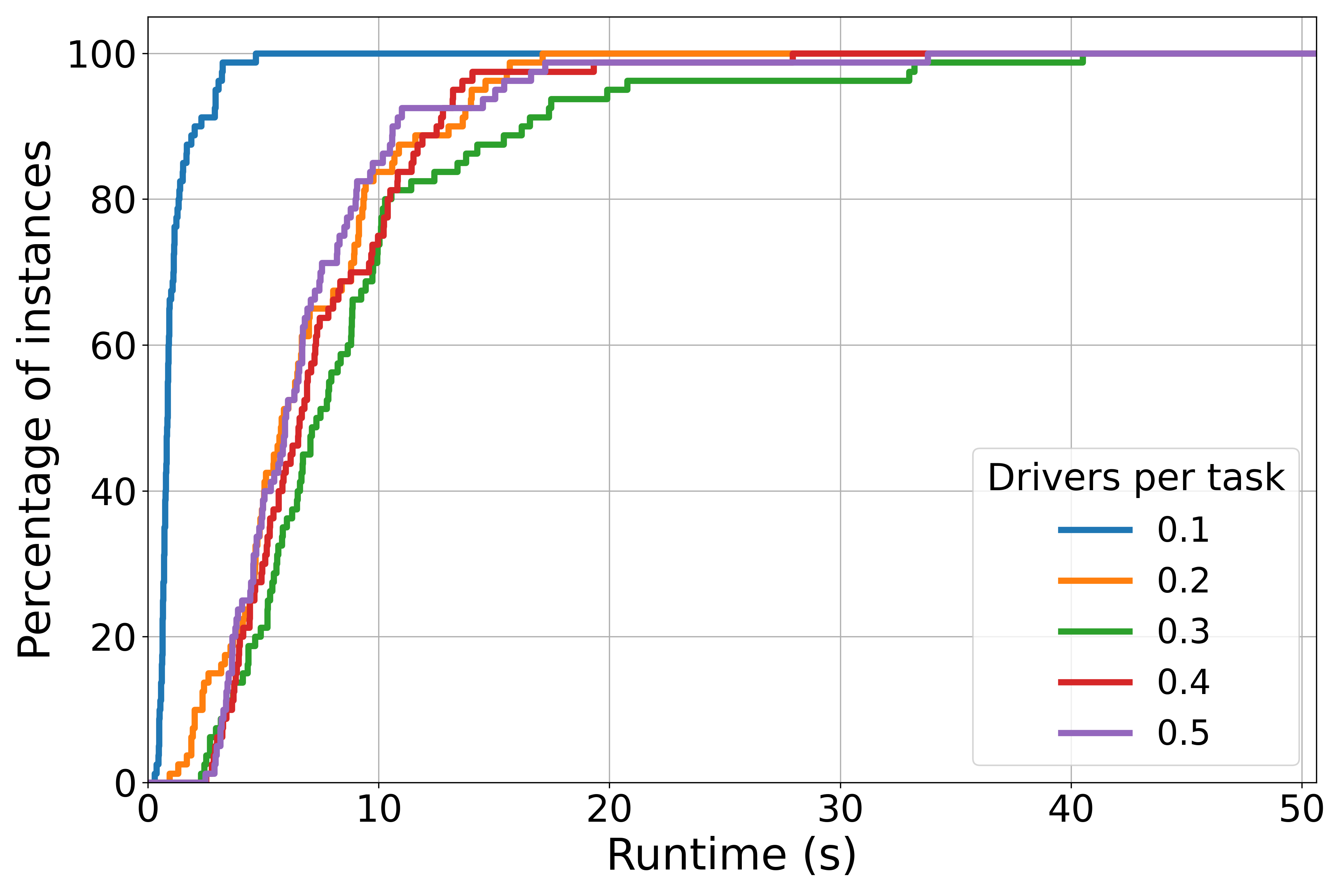}
		\caption{$|M|=60$}
		\label{fig:runtime_profile_m60}
	\end{subfigure}
	
	\vspace{0.6em}
	
	\begin{subfigure}[b]{0.48\textwidth}
		\centering
		\includegraphics[width=\linewidth]{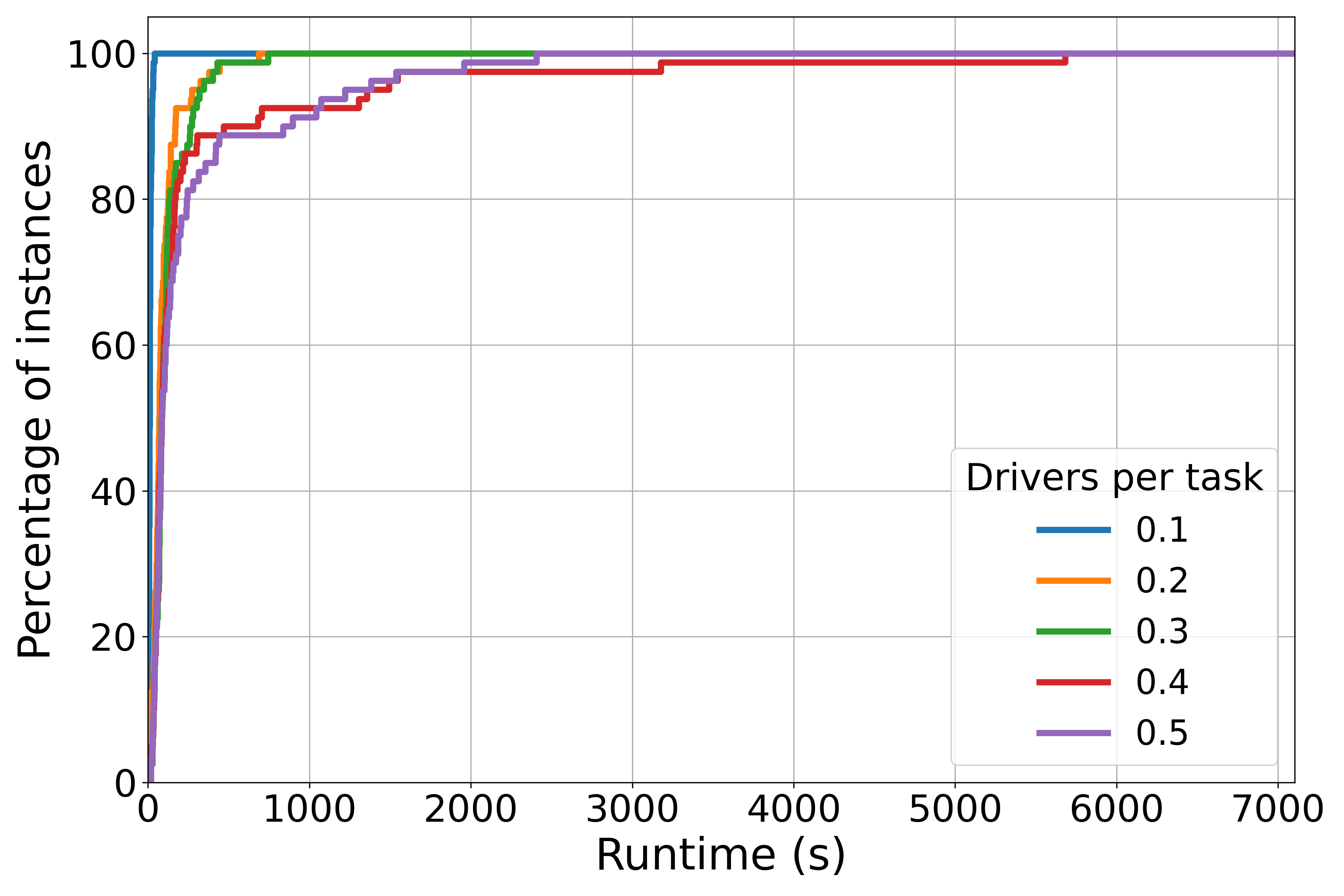}
		\caption{$|M|=90$}
		\label{fig:runtime_profile_m90}
	\end{subfigure}\hfill
	\begin{subfigure}[b]{0.48\textwidth}
		\centering
		\includegraphics[width=\linewidth]{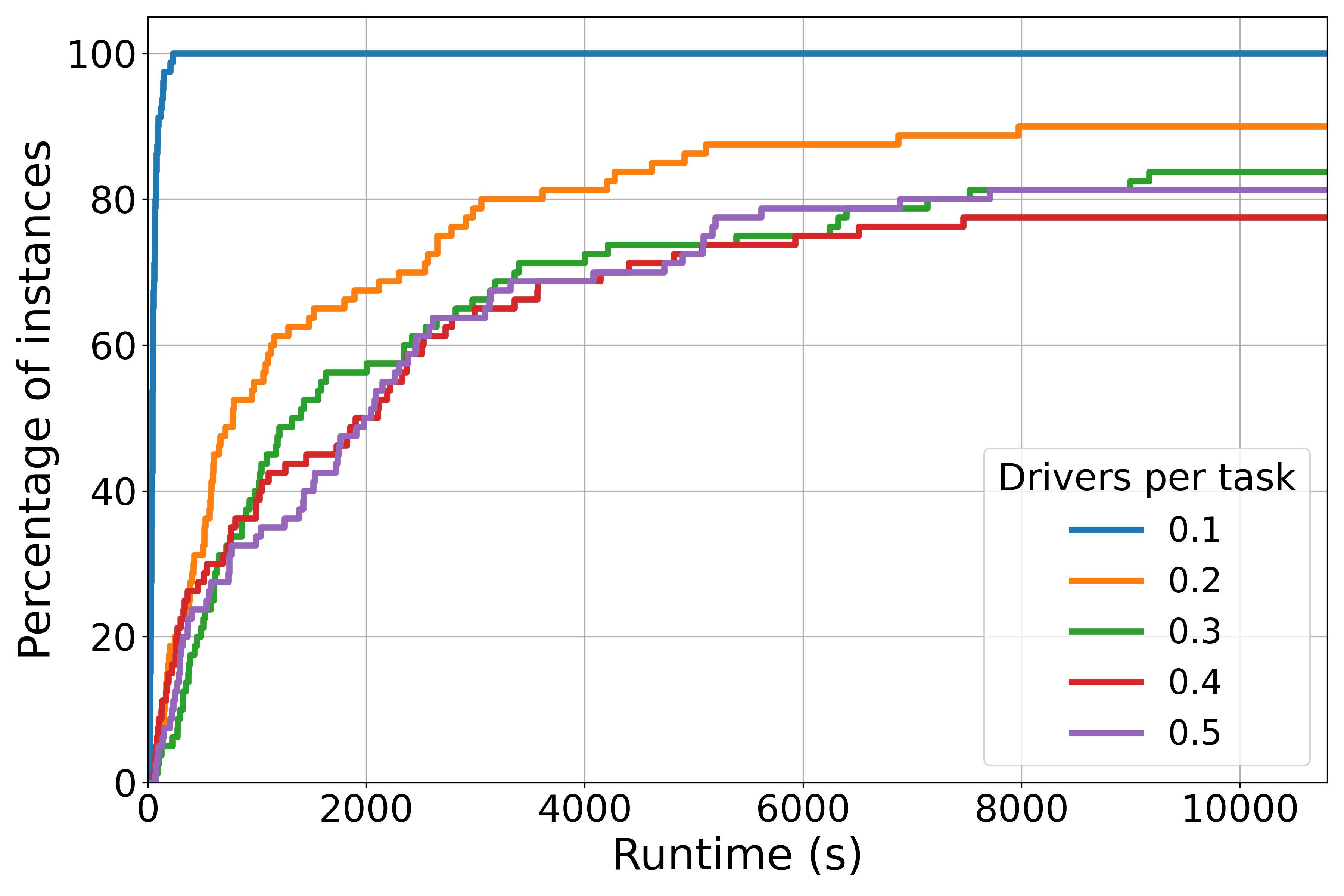}
		\caption{$|M|=120$}
		\label{fig:runtime_profile_m120}
	\end{subfigure}
	
	\caption{\blue{Runtime profiles of the exact method across varying task numbers $|M|$ and driver-to-task ratios.}}
	\label{fig:runtime_profiles_by_M}
\end{figure}

\subsection{Impact of behavioral classes on E-DDC}

Next, we analyze the impact of behavioral classes on the performance of E-DDC. \cref{tab:exact_drivers_classes} details numbers of solved instances and average runtimes for various numbers of tasks, ratios between drivers and tasks, and assignments of drivers to the behavioral classes discussed in \cref{sec:behavioralclasses}. 
%
These results confirm that E-DDC can rapidly solve all instances up to 60 tasks. They also show that the instances with a driver to task ratio of $p=0.1$ are easier to solve than those with a larger value of $p$. This trend is even more pronounced for larger instances. E-DDC also consistently solves \blue{all} instances with 90 tasks \blue{within 6000 seconds, with an average runtime of approximately 155 seconds.}
The proposed exact algorithm can solve all instances with 120 tasks and 12 drivers ($p=0.1$), but reaches its
limit on instances with 120 tasks and a higher value of $p$. While it can only solve 
\blue{more than 75\% of}
these instances, we note that the maximum optimality gaps are still very small ($\le 2.5\%$). 

\begin{table}[h]
	\caption{Numbers of instances solved by E-DDC and average runtimes in seconds for various numbers of drivers and behavioral class configurations (10 instances per class 1,2,3; and 50 instances for the mixed case).}
	\label{tab:exact_drivers_classes}
		\centering
		\blue{
			\begin{tabular}{ccrrrrrrrr}
				\toprule
				\multicolumn{2}{c}{\textbf{Instances}} & \multicolumn{4}{c}{\textbf{Solved Instances}} &  \multicolumn{4}{c}{\textbf{Average Runtime}} \\
				\cmidrule(lr){1-2} \cmidrule(lr){3-6} \cmidrule(lr){7-10}
				Tasks & Drivers per Task  & Class 1 & Class 2 & Class 3 & Mixed & Class 1 & Class 2  & Class 3 & Mixed \\
				\midrule
				\multirow[t]{5}{*}{30}  & 0.1        &           10 &           10 &           10 &         50 & 0.1 & 0.1 & 0.0 & 0.1 \\
				& 0.2        &           10 &           10 &           10 &         50 & 0.3 & 0.2 & 0.1 & 0.2 \\
				& 0.3        &           10 &           10 &           10 &         50 & 0.4 & 0.3 & 0.2 & 0.3 \\
				& 0.4        &           10 &           10 &           10 &         50 & 0.4 & 0.3 & 0.2 & 0.3 \\
				& 0.5        &           10 &           10 &           10 &         50 & 0.4 & 0.3 & 0.3 & 0.3 \\
				\midrule
				\multirow[t]{5}{*}{60}  & 0.1        &           10 &           10 &           10 &         50 & 1.6 & 1.1 & 0.7 & 1.1 \\
				& 0.2        &           10 &           10 &           10 &         50 & 11.5 & 7.4 & 3.2 & 6.4 \\
				& 0.3        &           10 &           10 &           10 &         50  & 13.2 & 10.6 & 5.5 & 8.6 \\
				& 0.4        &           10 &           10 &           10 &         50 & 10.2 & 6.4 & 3.7 & 8.0 \\
				& 0.5        &           10 &           10 &           10 &         50 & 6.2 & 5.7 & 3.4 & 8.1 \\
				\midrule
				\multirow[t]{5}{*}{90}  & 0.1        &           10 &           10 &           10 &         50 & 18.8 & 12.7 & 7.4 & 12.0 \\
				& 0.2        &           10 &           10 &           10 &         50 & 163.7 & 88.7 & 63.7 & 92.7 \\
				& 0.3        &           10 &           10 &           10 &         50 & 86.8 & 64.9 & 43.8 & 152.2 \\
				& 0.4        &           10 &           10 &           10 &         50 & 238.3 & 216.8 & 69.2 & 357.5 \\
				& 0.5        &           10 &           10 &           10 &         50 & 195.3 & 226.4 & 116.8 & 301.5 \\
				\midrule
				\multirow[t]{5}{*}{120} & 0.1        &           10 &           10 &           10 &         50  & 74.1 & 46.6 & 22.3 & 55.5 \\
				& 0.2        &            9 &            9 &            9 &         45 & 3560.3 & 2066.2 & 1389.2 & 2317.0 \\
				& 0.3        &            8 &            7 &           10 &         42 & 4237.4 & 4950.0 & 1814.0 & 3223.4 \\
				& 0.4        &            7 &            7 &            7 &         41 & 4403.6 & 4323.5 & 4027.1 & 3384.1 \\
				& 0.5        &            7 &            8 &            7 &         43 & 4437.1 & 4592.7 & 4139.6 & 3065.4 \\
				\bottomrule
			\end{tabular}
		}
	\end{table}

While the impact of the considered behavioral class on the numbers of solved instances seems limited, runtimes do indicate that instances that only feature Class 3 drivers are generally solved faster, likely because these drivers have the lowest sensitivity to compensation. As a result, fewer bundles need to be evaluated, reducing average runtimes. A similar, though less pronounced, trend is observed between Class 1 and Class 2 drivers. 
\blue{Moreover, runtimes for mixed-class instances either fall between those of the individual classes or are even lower. As indicated by the sensitivity analysis, this may be because, within a heterogeneous driver pool, some driver classes are prioritized over others, thereby reducing the effort spent on deprioritized classes and accelerating solution times.}
Overall, while driver class allocation does influence performance, the results suggest that the number of tasks and driver percentages are the primary determinants of runtime.

\subsection{Impact of behavioral classes and numbers of tasks and drivers per task on H-C and H-DDC}
\cref{tab:sensitivity_performance:heur} shows average runtimes of H-C and H-DDC for different numbers of drivers and behavioral classes. These results confirm the observations made for the exact variant E-DDC in \cref{sec:results}, in particular that the main parameters influencing the overall runtime are the numbers of available tasks and drivers. Again, instances can be solved faster when drivers have the lowest sensitivity to compensation (Class 3), instances of Class 2 seem to be slightly easier than those with Class 1 drivers, and \blue{mixed-class instances exhibit an intermediate level of complexity.}

\begin{table}[h]
	\caption{Average runtimes (s) of H-C and H-DDC for different numbers of drivers and behavioral classes (10 instances per configuration for classes 1,2,3 and 50 instances for the mixed case)}
	\label{tab:sensitivity_performance:heur}
	\blue{
		\begin{tabular}{ccrrrrrrrr}
			\toprule
			\multicolumn{2}{c}{\textbf{Instances}} & \multicolumn{4}{c}{\textbf{H-C}} & \multicolumn{4}{c}{\textbf{H-DDC}} \\ 
			\cmidrule(lr){1-2} \cmidrule(lr){3-6} \cmidrule(lr){7-10} 
			Tasks & Drivers per Task & Class 1 & Class 2 & Class 3 & Mixed & Class 1 & Class 2 & Class 3 & Mixed \\
			\midrule
			\multirow[t]{5}{*}{30} & 0.1 & 0.0 & 0.0 & 0.0 & 0.0 & 0.1 & 0.1 & 0.0 & 0.1 \\ 
			& 0.2 & 0.0 & 0.0 & 0.0 & 0.0 & 0.2 & 0.2 & 0.1 & 0.2\\ 
			& 0.3 & 0.0 & 0.0 & 0.0 & 0.0 & 0.3 & 0.2 & 0.1 & 0.2\\ 
			& 0.4 & 0.1 & 0.0 & 0.0 & 0.0 & 0.3 & 0.2 & 0.1 & 0.2\\ 
			& 0.5 & 0.0 & 0.0 & 0.1 & 0.0 & 0.3 & 0.2 & 0.2 & 0.2\\ 
			\midrule
			\multirow[t]{5}{*}{60} & 0.1 & 0.1 & 0.1 & 0.1 & 0.1 & 1.3 & 0.8 & 0.5 & 0.8 \\ 
			& 0.2 & 0.3 & 0.2 & 0.2 & 0.2 & 5.7 & 3.7 & 1.8 & 3.5 \\ 
			& 0.3 & 0.5 & 0.4 & 0.3 & 0.4 & 6.1 & 4.2 & 2.2 & 4.5 \\ 
			& 0.4 & 0.6 & 0.5 & 0.4 & 0.5 & 4.6 & 3.3 & 2.0 & 4.1 \\ 
			& 0.5 & 0.7 & 0.6 & 0.5 & 0.7 & 3.3 & 3.1 & 2.0 & 3.9 \\ 
			\midrule
			\multirow[t]{5}{*}{90} & 0.1 & 0.7 & 0.6 & 0.4 & 0.6 & 12.3 & 8.6 & 4.2 & 7.6 \\ 
			& 0.2 & 2.9 & 2.3 & 1.6 & 2.2 & 51.4 & 29.5 & 13.2 & 28.3 \\ 
			& 0.3 & 4.2 & 4.2 & 3.0 & 3.8 & 31.8 & 25.3 & 14.0 & 29.2 \\ 
			& 0.4 & 4.8 & 4.7 & 3.9 & 4.2 & 22.7 & 21.4 & 14.2 & 25.1 \\ 
			& 0.5 & 6.1 & 5.7 & 5.0 & 5.0 & 19.5 & 18.5 & 13.5 & 22.7 \\ 
			\midrule
			\multirow[t]{5}{*}{120} & 0.1 & 4.2 & 3.3 & 2.0 & 3.2 & 38.5 & 24.7 & 13.5 & 28.8 \\ 
			& 0.2 & 16.2 & 13.0 & 8.5 & 12.2 & 201.2 & 127.0 & 51.7 & 121.3 \\ 
			& 0.3 & 21.5 & 23.4 & 34.2 & 22.7 & 166.5 & 136.9 & 78.6 & 126.2 \\ 
			& 0.4 & 23.3 & 27.0 & 22.7 & 25.4 & 98.4 & 133.6 & 89.1 & 111.1 \\ 
			& 0.5 & 74.2 & 31.6 & 29.1 & 30.5 & 161.2 & 90.4 & 107.5 & 98.1 \\ 
			\bottomrule
		\end{tabular}
	}
\end{table}

\section{Sensitivity analysis in a multi-depot setting}\label{sec:md}

In this section, we evaluate our algorithm on instances featuring multiple depots, spatially clustered tasks, and drivers with \blue{homogenous} capacities \blue{(\cref{sec:md_hom}) and with heterogenous capacities (\cref{sec:md_het})} that are created by extending the original set of instances. Its main goal is to understand whether the conclusions drawn in \cref{sec:sensitivity} remain valid for instances with modified structure. As in \cref{sec:sensitivity}, we therefore focus on the sensitivity of four key criteria (acceptance probability, compensation, detour, and bundle size) of the offers made to occasional drivers using best-known solutions.
In the new set of multi-depot instances, 
all drivers are assumed to start at the origin and 
the service area is partitioned into four zones corresponding to the four quadrants. A depot is located at the center of each quadrant, at coordinates $(\pm 2.5, \pm 2.5)$. Tasks are assigned to their corresponding zones and can only be served by the depot in their zone. No depot constraints are enforced on the drivers, hence the operator may construct bundles originating from any depot. Finally, each driver is assigned a heterogeneous capacity \blue{in \cref{sec:md_het}}, drawn uniformly at random between 50 and 150.

\subsection{\blue{Multi-depot setting with homogenous capacity}} \label{sec:md_hom}

\cref{fig:md_hom_task_all_metrics} illustrates how the number of tasks affects the four key metrics. Overall, the qualitative patterns mirror those in the single-depot setting, i.e., are in line with the trends in \cref{fig:task_all_metrics}. One noticeable trend is that, as the number of tasks increases, compensation rises more sharply for Class~3 up to 90 tasks, with a similar (though milder) increase emerging for Class~2. Combined with the multi-depot setting’s consistently lower acceptance probabilities and slightly smaller bundle sizes, alongside higher compensations and larger detours across all classes, this suggests that the multi-depot setting amplifies the patterns observed in the single-depot case. 
	
A similar observation holds for the general trends in \cref{fig:md_hom_driver_all_metrics}, which reports the impact of the driver-to-task ratio in the multi-depot setting, in comparison to its single-depot counterpart in \cref{fig:driver_all_metrics}. Some noticeable differences are that Class~3 acceptance probabilities increase from 0.2 to 0.4, and that this class exhibits slightly different patterns in compensation, bundle size, and detour particularly around these values. This suggests that, because the multi-depot setting offers fewer feasible bundling opportunities, the optimizer adjusts Class~3 offers more actively as needed. This interpretation is further supported by the comparatively more stable confidence intervals for Class~3 across supply levels relative to its single-depot counterpart. We also observe a shift in class prioritization in terms of compensation (especially at driver-per-task ratios of 0.4 and 0.5) indicating that the model increases payments to more compensation sensitive classes to sustain high acceptance rates.

\begin{figure}[H]
	\captionsetup[subfigure]{aboveskip=-0.5pt,belowskip=-0.5pt}
	\begin{subfigure}[t]{0.49\textwidth}
		\centering
		\includegraphics[width=0.9\textwidth]{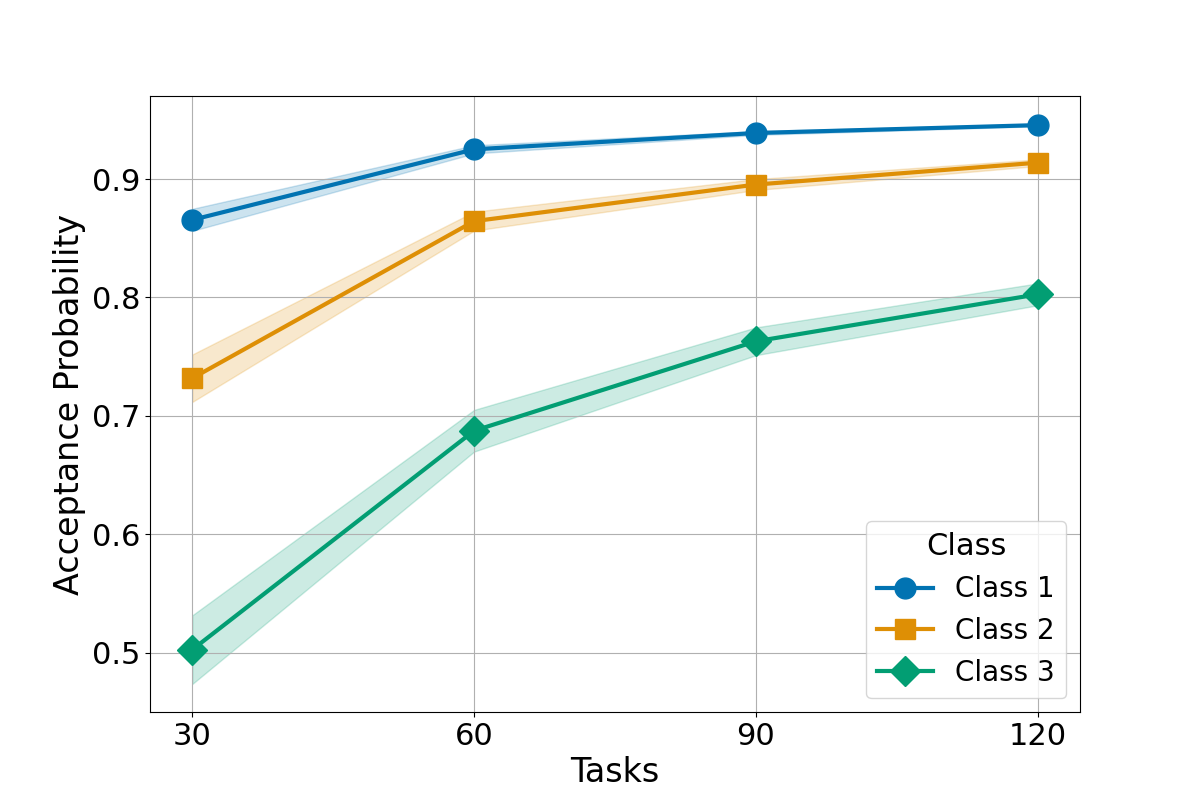}
		\caption{Acceptance probability.} 
	\label{fig:md_hom_task_acceptance}
\end{subfigure}	\hfill	
\begin{subfigure}[t]{0.49\textwidth}
	\centering
	\includegraphics[width=0.9\textwidth]{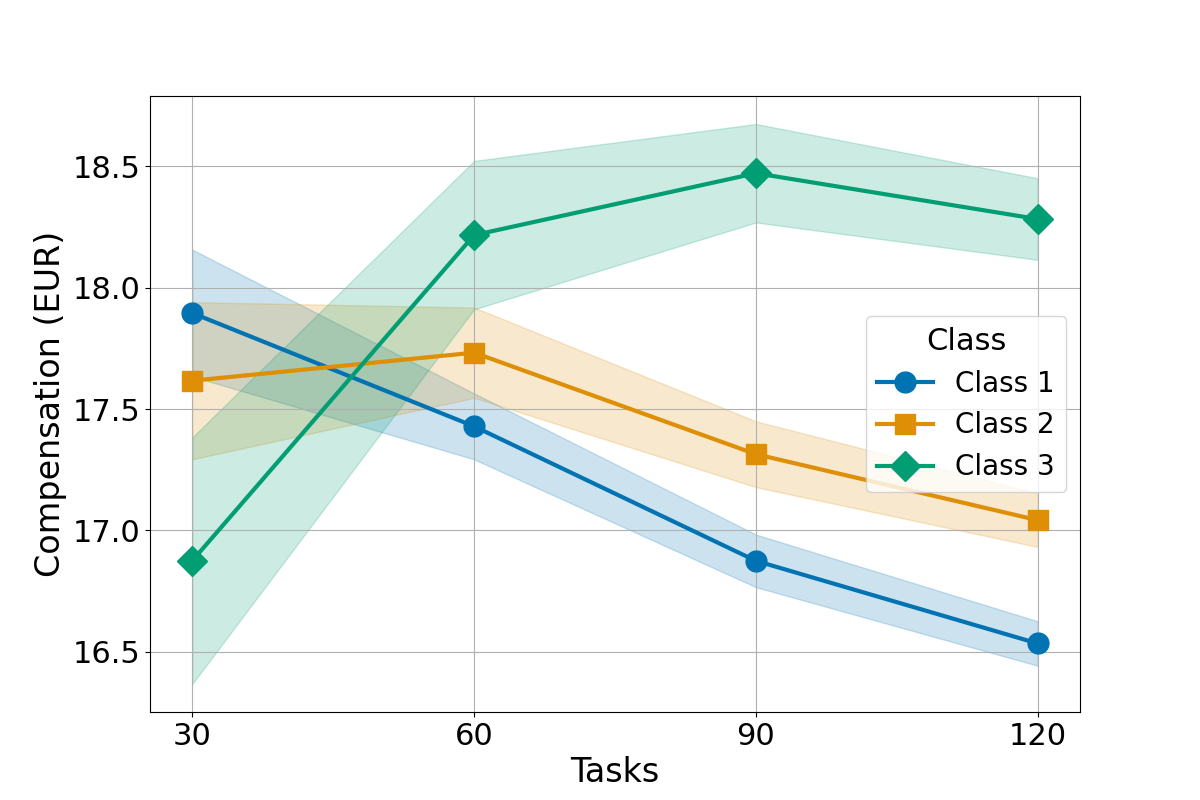}
	\caption{Compensation.} 
\label{fig:md_hom_task_compensation}
\end{subfigure}
\begin{subfigure}[t]{0.49\textwidth}
\centering
\includegraphics[width=0.9\textwidth]{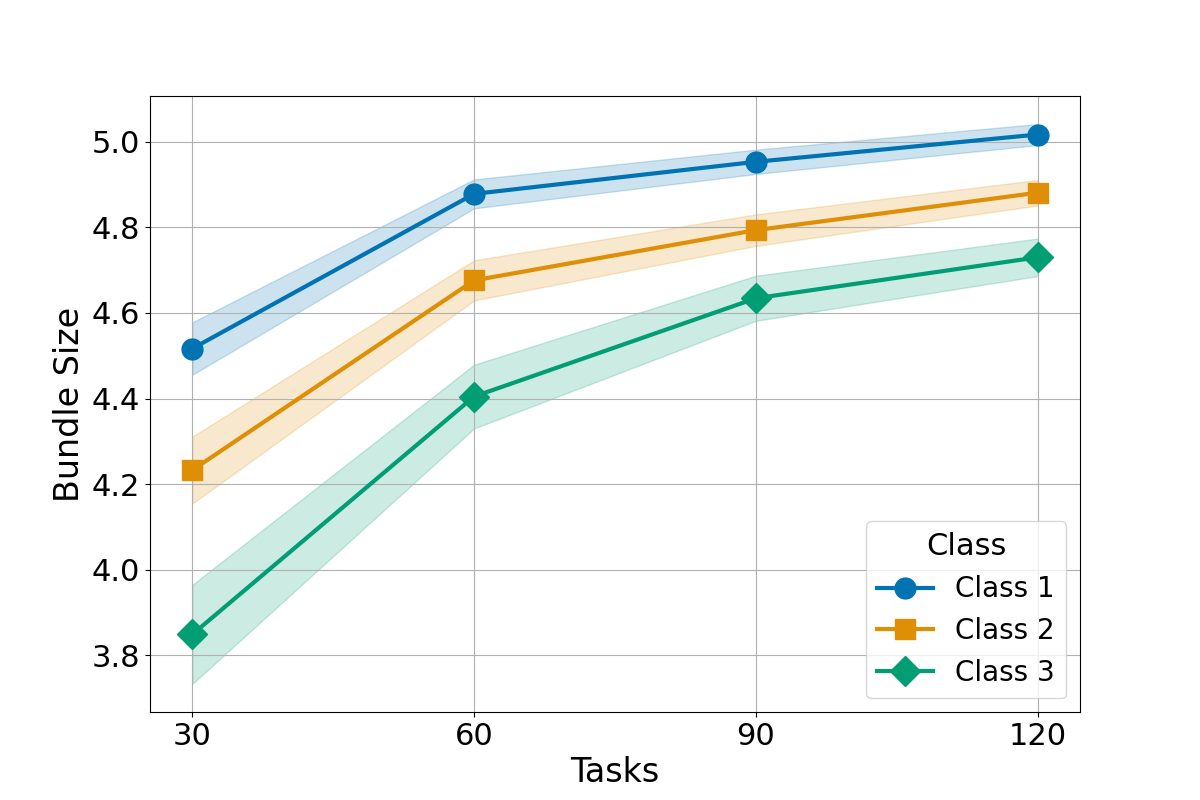}
\caption{Bundle size.} 
\label{fig:md_hom_task_bundle}
\end{subfigure}	\hfill	
\begin{subfigure}[t]{0.49\textwidth}
\centering
\includegraphics[width=0.9\textwidth]{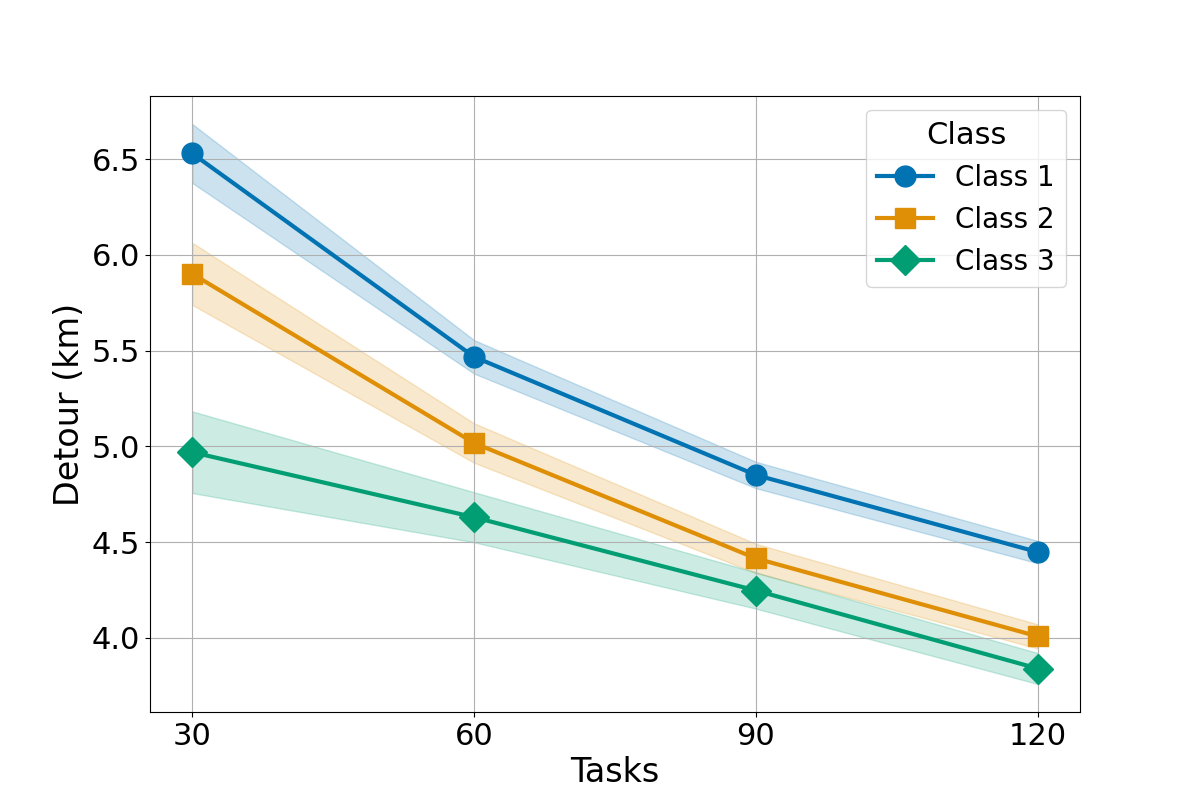}
\caption{Detour.} 
\label{fig:md_hom_task_detour}
\end{subfigure}
\caption{Mean values and their 95\% confidence intervals of different metrics for different numbers of tasks and behavioral classes in a multi-depot setting \blue{with homogenous capacities}.}
\label{fig:md_hom_task_all_metrics}
\vspace{-1em}
\end{figure}

\begin{figure}[H]
\captionsetup[subfigure]{aboveskip=-0.5pt,belowskip=-0.5pt}
\begin{subfigure}[t]{0.49\textwidth}
\centering
\includegraphics[width=0.9\textwidth]{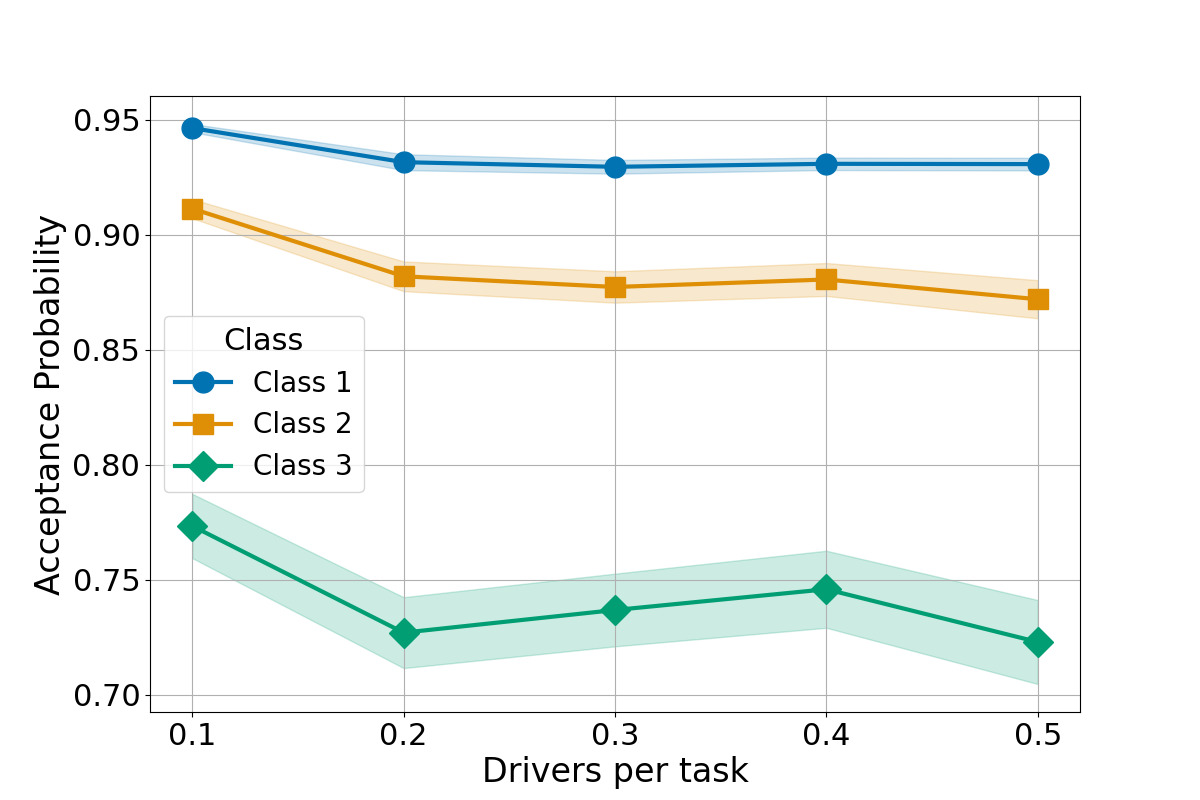}
\caption{Acceptance probability} 
\label{fig:md_hom_driver_acceptance}
\end{subfigure}	\hfill	
\begin{subfigure}[t]{0.49\textwidth}
\centering
\includegraphics[width=0.9\textwidth]{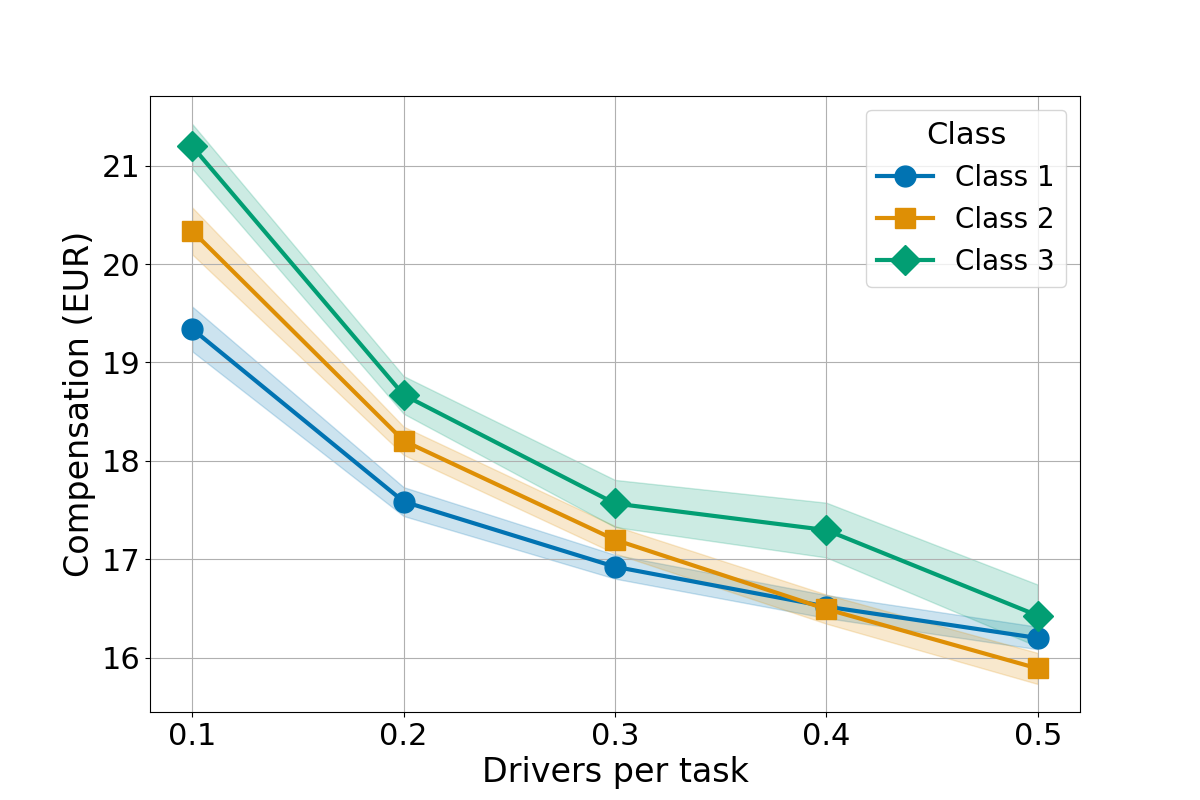}
\caption{Compensation} 
\label{fig:md_hom_driver_compensation}
\end{subfigure}
\begin{subfigure}[t]{0.49\textwidth}
\centering
\includegraphics[width=0.9\textwidth]{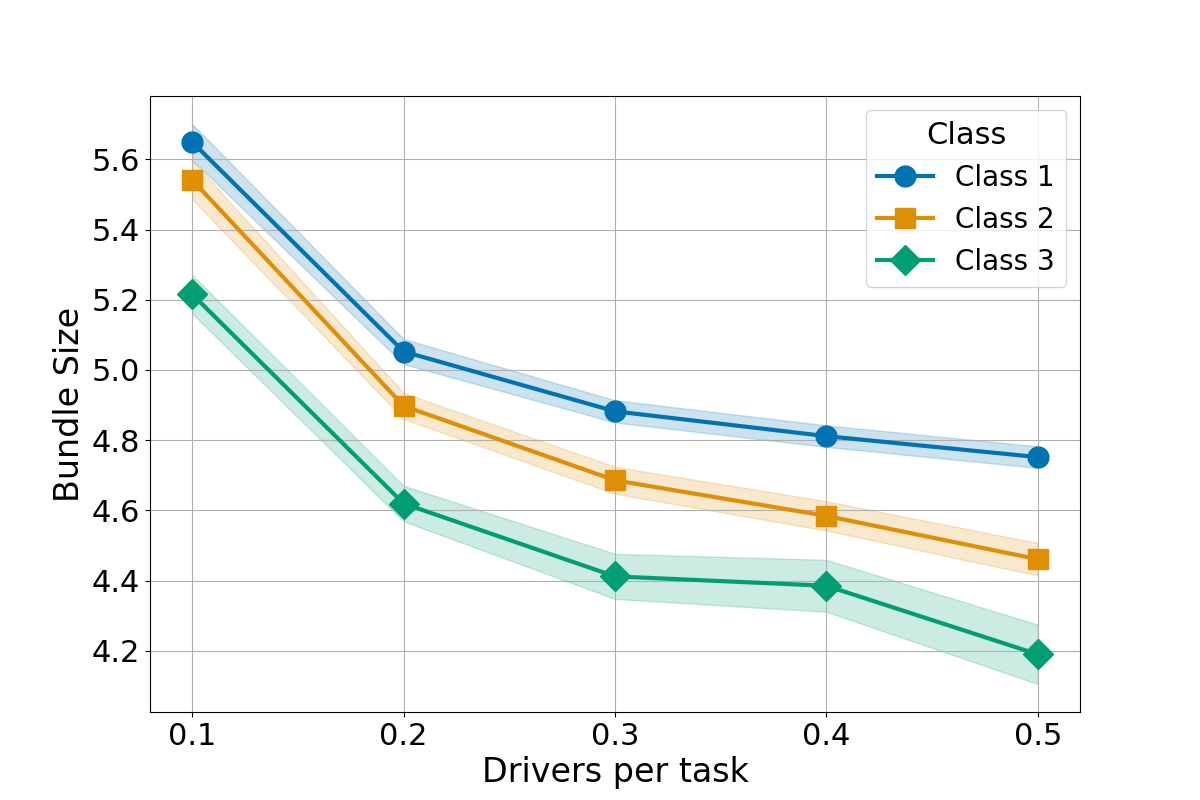}
\caption{Bundle size} 
\label{fig:md_hom_driver_bundle}
\end{subfigure}	\hfill	
\begin{subfigure}[t]{0.49\textwidth}
\centering
\includegraphics[width=0.9\textwidth]{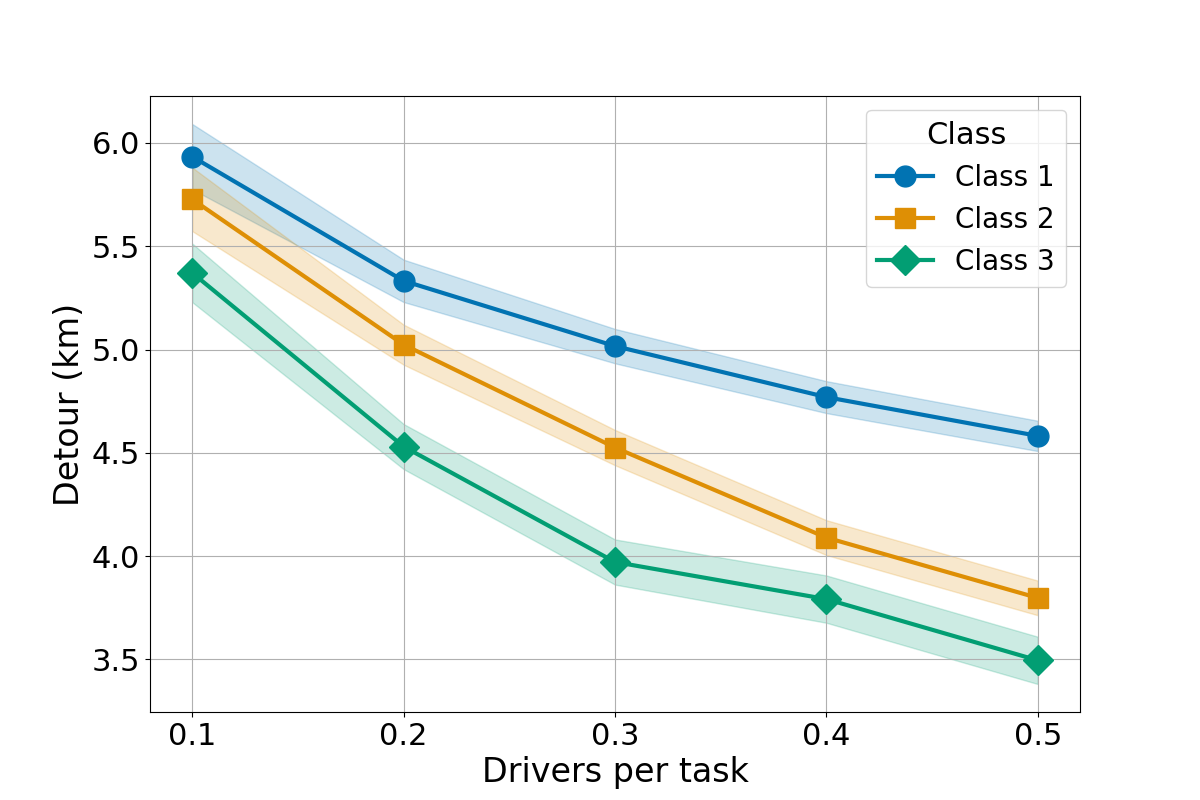}
\caption{Detour} 
\label{fig:md_hom_driver_detour}
\end{subfigure}

\caption{Mean values and their 95\% confidence intervals of different metrics for different driver to task ratios and behavioral classes in a multi-depot setting \blue{with homogenous capacities}.} 
\label{fig:md_hom_driver_all_metrics}
\end{figure}

Finally, we find that the general trends observed when comparing single-class and mixed-class cases in the single-depot setting (\cref{tab:mean_metrics}) also carry over to the multi-depot setting (\cref{tab:md_hom_mean_metrics}). Moreover, the same systematic metric shifts discussed above (lower acceptance probabilities and smaller bundle sizes, together with higher compensations and larger detours) remain present here as well. A plausible explanation for these shifts is that drivers must first travel to a depot that is not co-located with their initial positions, which mechanically increases detour. In addition, the smaller task pool available at each depot can reduce bundling opportunities, thereby increasing the effective effort per served task. The model responds to these frictions by raising compensations, but acceptance probabilities still remain below those observed in the single-depot setting.

\begin{table}[h]
\caption{Mean values of different metrics for different behavioral classes in a multi-depot setting \blue{with homogenous capacities}.} 
\label{tab:md_hom_mean_metrics}
\resizebox{\textwidth}{!}{%
	\begin{tabular}{lrrrrrr}
\toprule
& \multicolumn{2}{c}{Class 1} & \multicolumn{2}{c}{Class 2} & \multicolumn{2}{c}{Class 3} \\ \cmidrule(lr){2-3} \cmidrule(lr){4-5} \cmidrule(lr){6-7}
& Mixed-Class & Single-Class & Mixed-Class & Single-Class & Mixed-Class & Single-Class \\
\midrule
Acceptance Probability & 0.93 & 0.93 & 0.88 & 0.89 & 0.69 & 0.79 \\
Bundle Size & 4.99 & 4.76 & 4.73 & 4.80 & 4.23 & 4.82 \\
Compensation & 17.22 & 16.26 & 17.23 & 17.46 & 17.06 & 19.14 \\
Detour & 5.11 & 4.63 & 4.48 & 4.55 & 4.02 & 4.40 \\
\bottomrule
\end{tabular}
}
\end{table}

\subsection{\blue{Multi-depot setting with heterogeneous capacities}} \label{sec:md_het}

\blue{
\cref{fig:md_het_task_all_metrics} illustrates how the number of tasks affects the four key metrics when driver capacities are heterogeneously distributed. Overall, the qualitative patterns are in line with the trends in \cref{fig:md_hom_task_all_metrics}. Acceptance probabilities and bundle sizes increase in the number of tasks, while detours decrease for all classes, and bundle sizes, compensations, and detours are uniformly higher in nominal terms compared to their multi-depot homogeneous counterparts. Three differences are noticeable. First, Class~1's compensation now increases up to 60 tasks before decreasing, whereas it declines monotonically under homogeneous capacities; Class~2 exhibits the same trend in both settings. Second, Class~3's compensation keeps increasing up to 120 tasks, at a decreasing rate, whereas it turns downward after 90 tasks under homogeneous capacities. Third, the gap in acceptance probability between Class~3 and the two other classes narrows, and this narrowing is driven mainly by lower acceptance probabilities for Classes~1 and~2 rather than by higher ones for Class~3. In the same line, the bundle sizes of Classes~1 and~2 become statistically indistinguishable at 120 tasks, whereas they remain clearly separated under homogeneous capacities. These observations could be attributed to the fact that the model concentrates tasks on the drivers with larger capacities. Serving more tasks per bundle requires higher compensation to keep the additional effort acceptable to the drivers, and it leaves less room to differentiate bundle sizes across behavioral classes.

	\begin{figure}[H]
	\captionsetup[subfigure]{aboveskip=-0.5pt,belowskip=-0.5pt}
	\begin{subfigure}[t]{0.49\textwidth}
		\centering
		\includegraphics[width=0.9\textwidth]{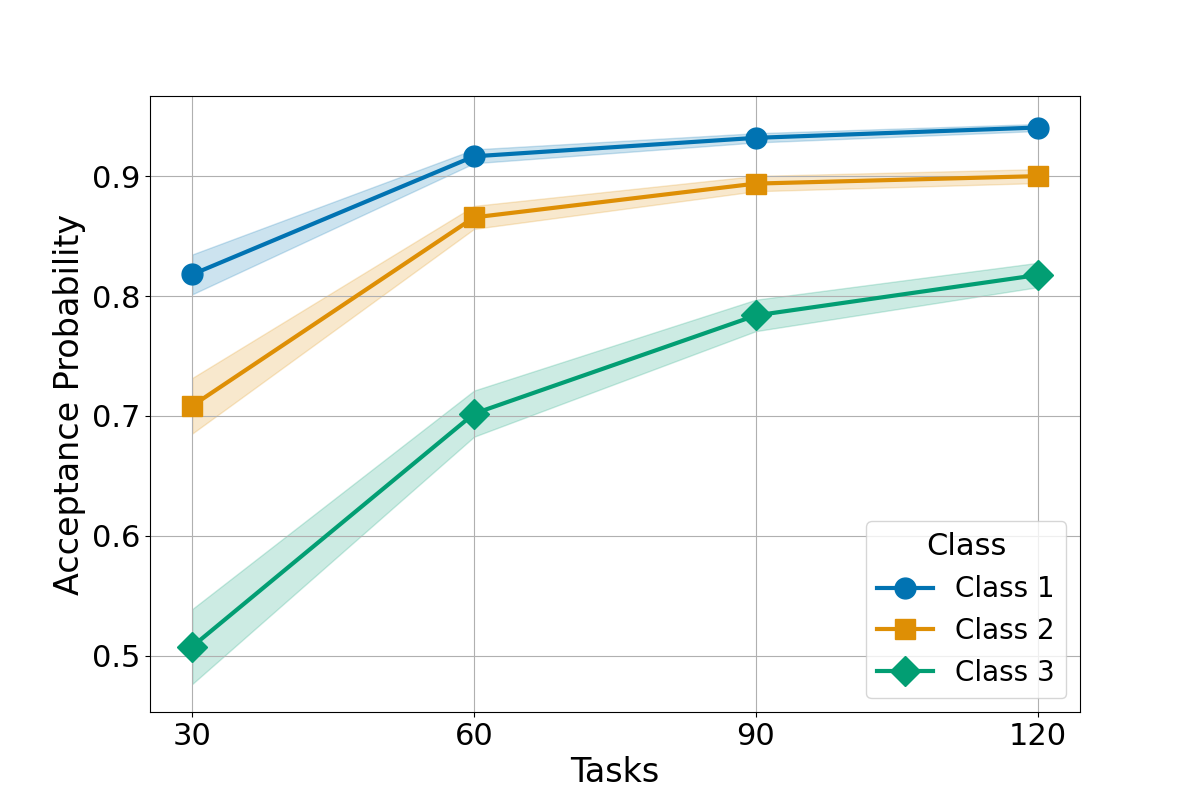}
		\caption{Acceptance probability.} 
	\label{fig:md_het_task_acceptance}
\end{subfigure}	\hfill	
\begin{subfigure}[t]{0.49\textwidth}
	\centering
	\includegraphics[width=0.9\textwidth]{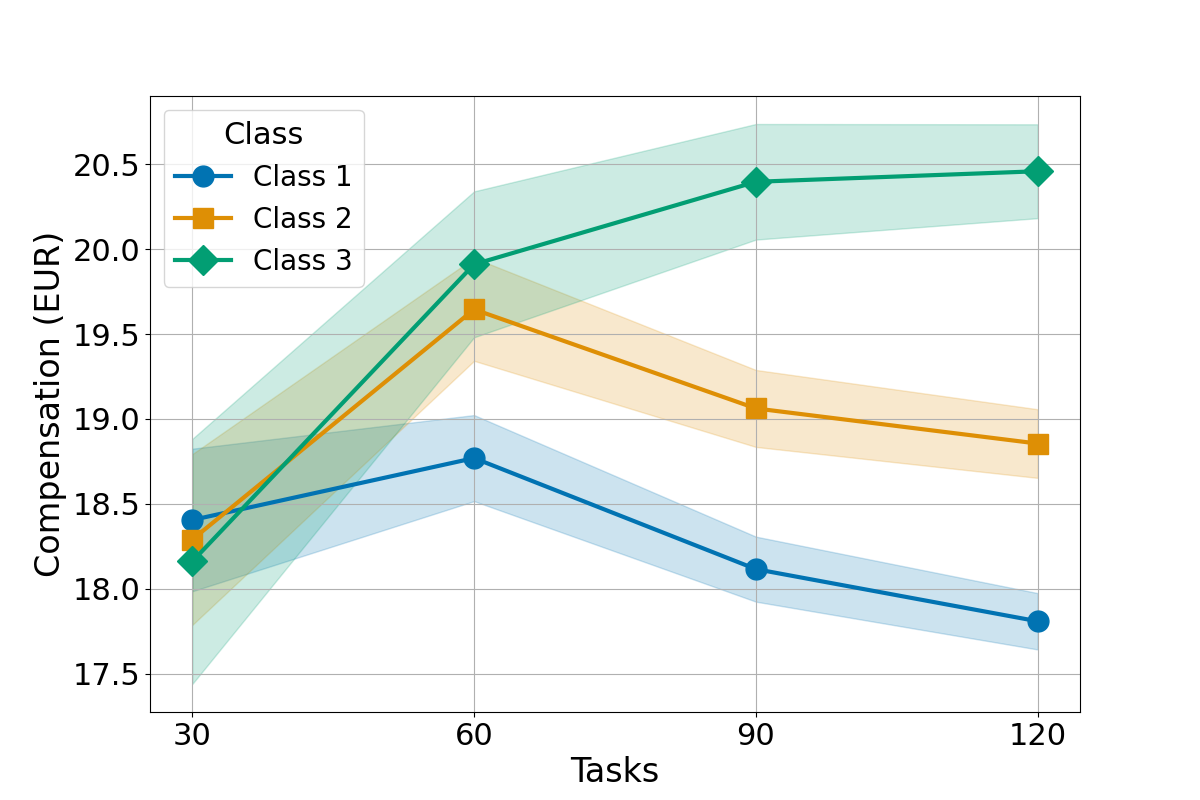}
	\caption{Compensation.} 
\label{fig:md_het_task_compensation}
\end{subfigure}
\begin{subfigure}[t]{0.49\textwidth}
\centering
\includegraphics[width=0.9\textwidth]{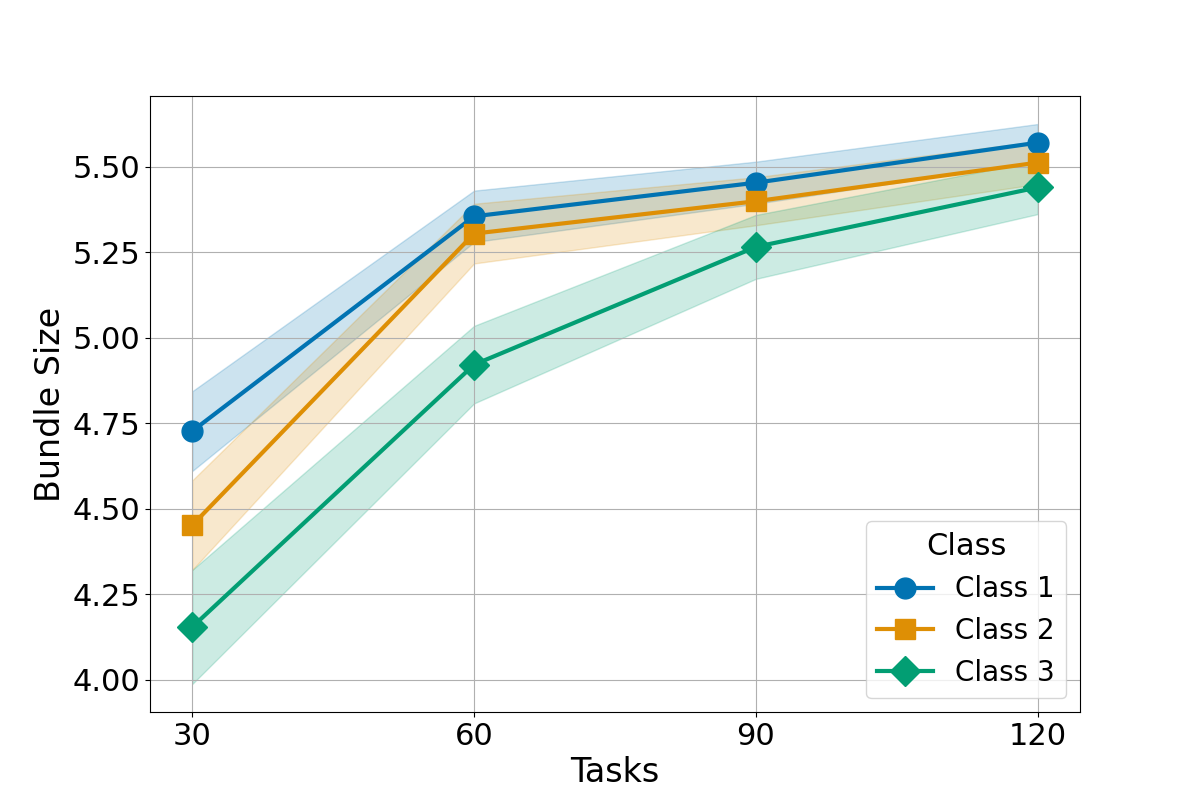}
\caption{Bundle size.} 
\label{fig:md_het_task_bundle}
\end{subfigure}	\hfill	
\begin{subfigure}[t]{0.49\textwidth}
\centering
\includegraphics[width=0.9\textwidth]{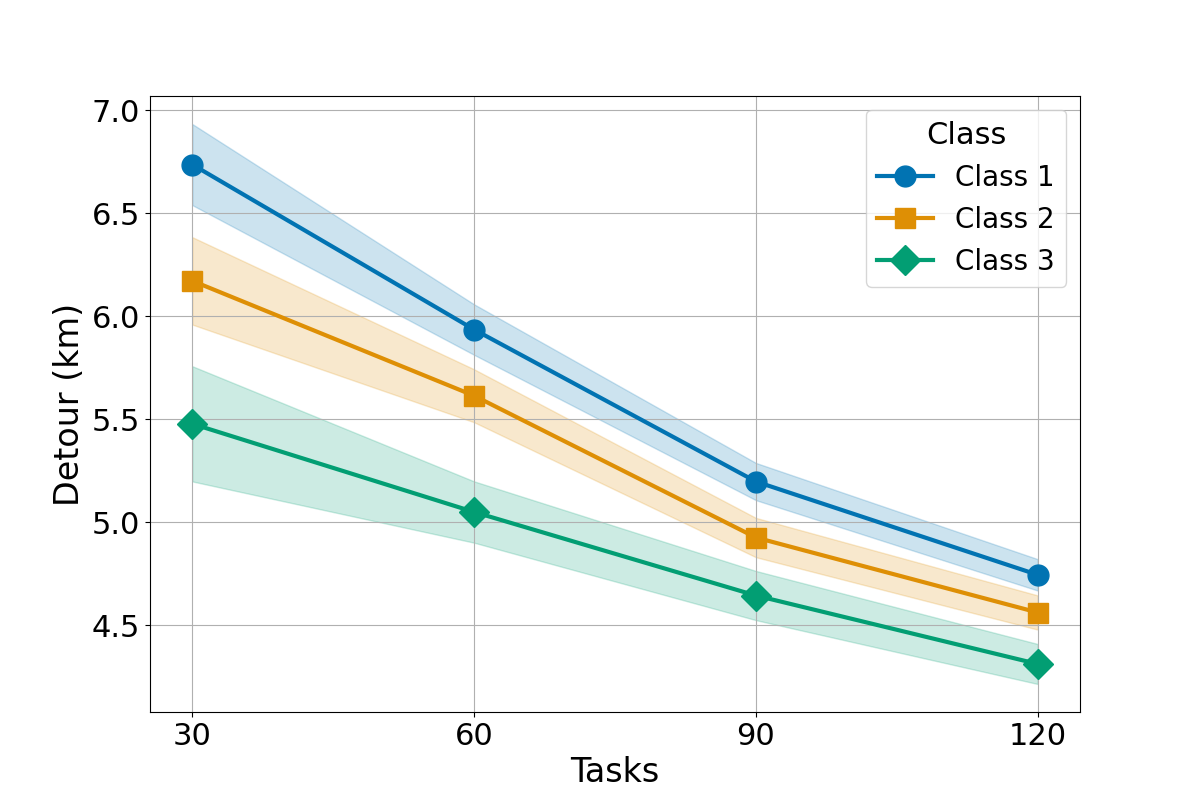}
\caption{Detour.} 
\label{fig:md_het_task_detour}
\end{subfigure}
\caption{\blue{Mean values and their 95\% confidence intervals of different metrics for different numbers of tasks and behavioral classes in a multi-depot setting with heterogenous capacities.}}
\label{fig:md_het_task_all_metrics}
\end{figure}

The general trends in \cref{fig:md_het_driver_all_metrics}, which reports the impact of the driver-to-task ratio under heterogeneous capacities, also carry over from its homogeneous counterpart in \cref{fig:md_hom_driver_all_metrics}, but with several differences. First, under homogeneous capacities acceptance probabilities drop from 0.1 to 0.2 drivers per task and then flatten, whereas under heterogeneous capacities they are essentially flat for Class~1 and mildly increasing for Classes~2 and~3, the only exception being a slight decrease for Class~3 from 0.3 to 0.4 drivers per task. Second, Class~3's compensation declines steadily, Class~1's is flat after an initial drop from 0.1 to 0.2 drivers per task, and Class~2's continues to drift downward mildly in contrast to the consistently decreasing compensations in the homogenous setting. Finally, bundle size decreases consistently for Class~3, while for Classes~1 and~2 it first decreases from 0.1 to 0.2 drivers per task and increases afterwards, in contrast to the monotone decrease observed for all classes under homogeneous capacities. These results support the conclusion that the model prioritizes the drivers with larger capacities, but that it still accounts for the prioritization across behavioral classes while doing so.

Finally, the general trends observed when comparing single-class and mixed-class cases under homogeneous capacities (\cref{tab:md_hom_mean_metrics}) carry over to the heterogeneous setting (\cref{tab:md_het_mean_metrics}). The direction of every mixed-class versus single-class difference is preserved, for all four metrics and all three classes. Class~1 continues to be better off in mixed-class instances in terms of bundle size, compensation, and detour, whereas Classes~2 and~3 receive smaller bundles, lower compensation, and shorter detours than in their single-class counterparts. The nominal values of bundle size, compensation, and detour are uniformly higher under heterogeneous capacities, while acceptance probabilities are slightly lower for Classes~1 and~2. Taken together with the preserved compensation ordering discussed above, these observations support the conclusion that, when the drivers have heterogeneous capacities, the model concentrates tasks on those with larger capacities, which results in larger bundles, longer detours, and larger compensation to keep the additional effort acceptable to the drivers.
}

\begin{figure}
\captionsetup[subfigure]{aboveskip=-0.5pt,belowskip=-0.5pt}
\begin{subfigure}[t]{0.49\textwidth}
\centering
\includegraphics[width=0.9\textwidth]{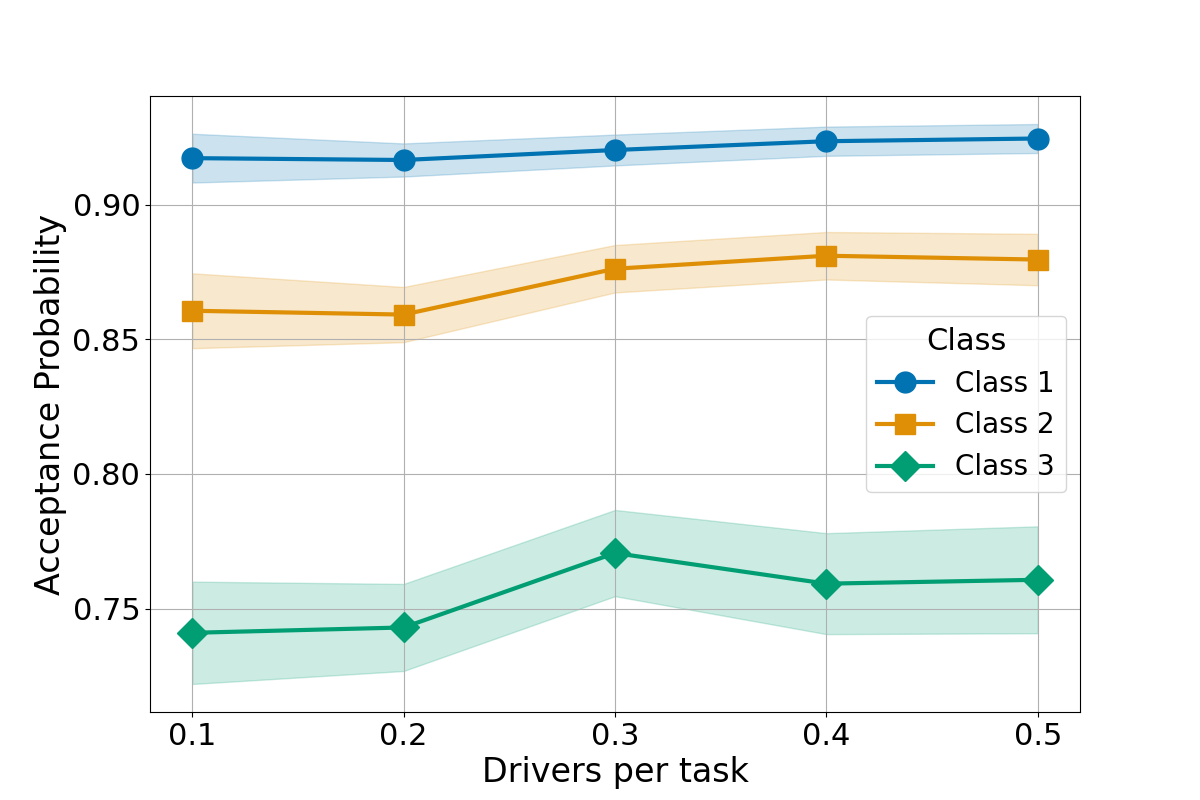}
\caption{Acceptance probability} 
\label{fig:md_het_driver_acceptance}
\end{subfigure}	\hfill	
\begin{subfigure}[t]{0.49\textwidth}
\centering
\includegraphics[width=0.9\textwidth]{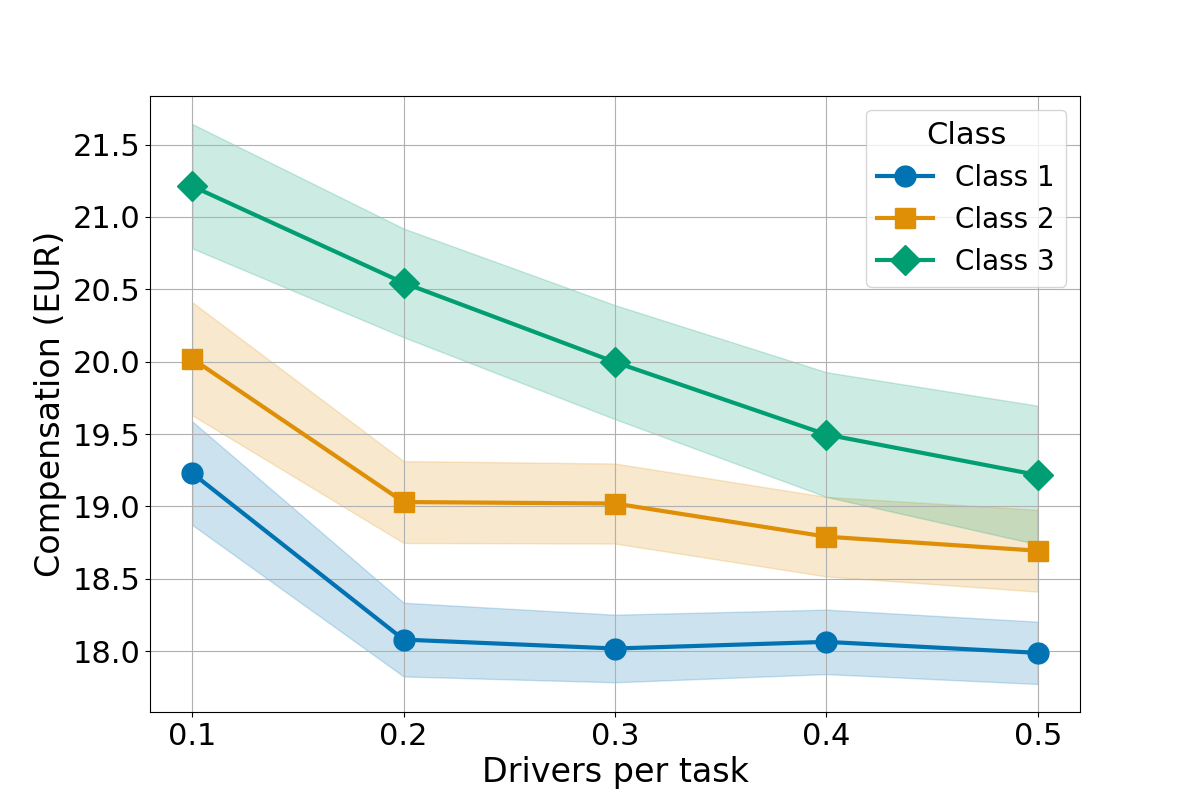}
\caption{Compensation} 
\label{fig:md_het_driver_compensation}
\end{subfigure}
\begin{subfigure}[t]{0.49\textwidth}
\centering
\includegraphics[width=0.9\textwidth]{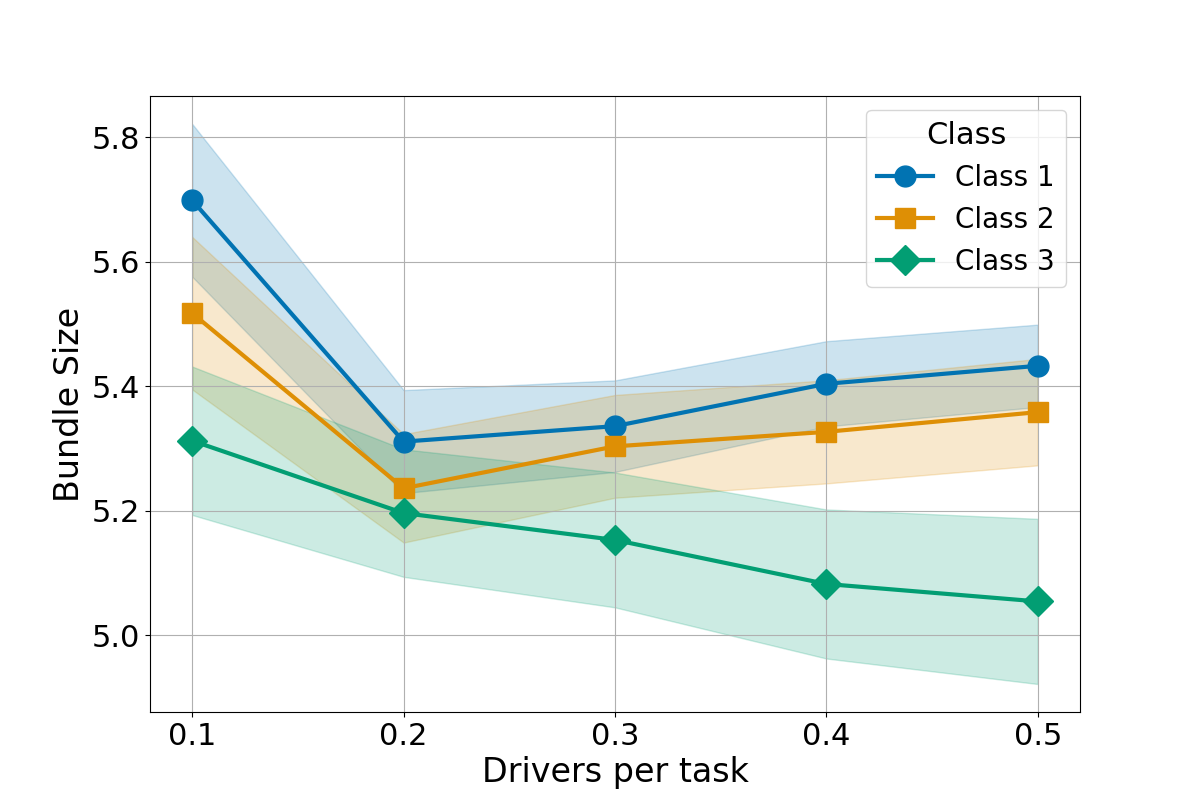}
\caption{Bundle size} 
\label{fig:md_het_driver_bundle}
\end{subfigure}	\hfill	
\begin{subfigure}[t]{0.49\textwidth}
\centering
\includegraphics[width=0.9\textwidth]{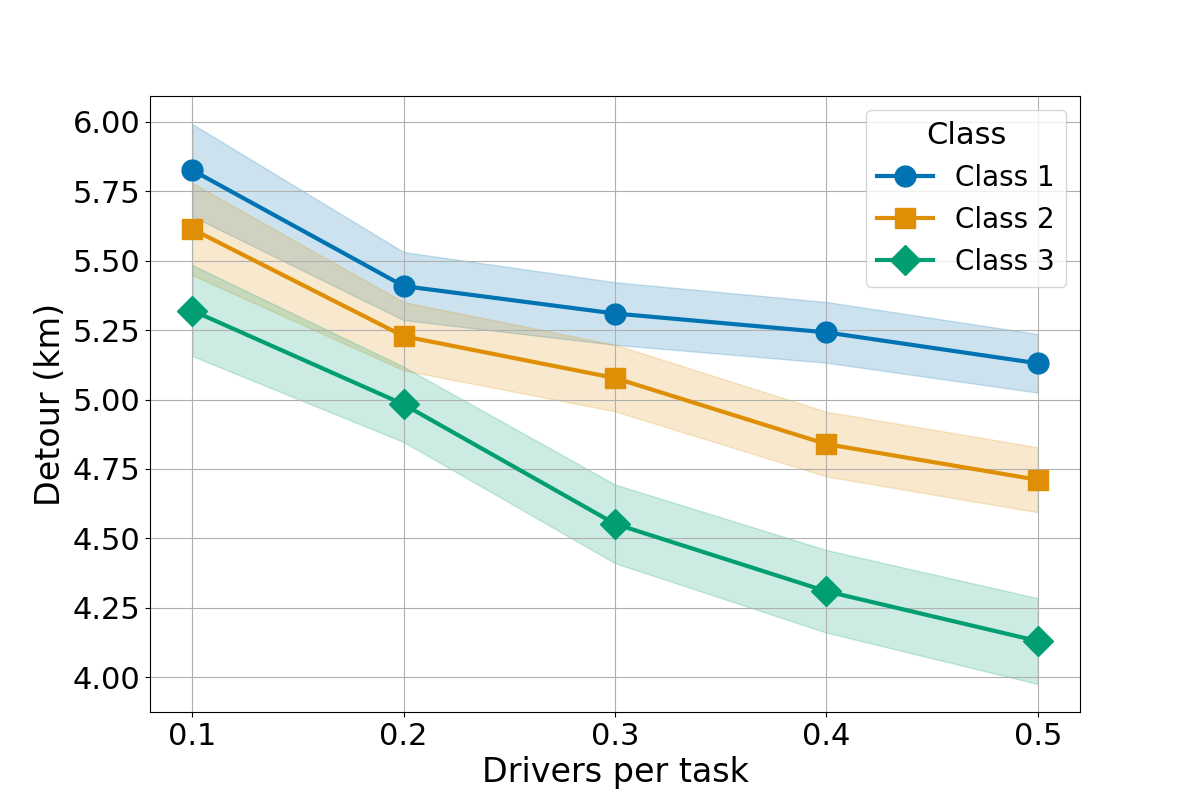}
\caption{Detour} 
\label{fig:md_het_driver_detour}
\end{subfigure}

\caption{\blue{Mean values and their 95\% confidence intervals of different metrics for different driver to task ratios and behavioral classes in a multi-depot setting with heterogenous capacities.}} 
\label{fig:md_het_driver_all_metrics}
\end{figure}

\begin{table}
\caption{\blue{Mean values of different metrics for different behavioral classes in a multi-depot setting with heterogenous capacities.}}
\label{tab:md_het_mean_metrics}
\resizebox{\textwidth}{!}{%
\blue{
\begin{tabular}{lrrrrrr}
\toprule
& \multicolumn{2}{c}{Class 1} & \multicolumn{2}{c}{Class 2} & \multicolumn{2}{c}{Class 3} \\ \cmidrule(lr){2-3} \cmidrule(lr){4-5} \cmidrule(lr){6-7}
& Mixed-Class & Single-Class & Mixed-Class & Single-Class & Mixed-Class & Single-Class \\
\midrule
Acceptance Probability & 0.92 & 0.92 & 0.87 & 0.88 & 0.72 & 0.78 \\
Bundle Size & 5.44 & 5.34 & 5.29 & 5.40 & 4.91 & 5.40 \\
Compensation & 18.30 & 17.81 & 18.91 & 19.21 & 19.23 & 20.95 \\
Detour & 5.40 & 5.12 & 5.01 & 5.08 & 4.48 & 4.85 \\
\bottomrule
\end{tabular}
}
}
\end{table}

\blue{
\section{Absolute driver utility assessment}\label{sec:absolute-driver-utility}

To complement the cross-class utility ranges reported in \cref{tab:welfare_money_metric_sensitivity}, \cref{tab:absolute-driver-utility} reports the absolute levels of the money-metric utility proxy for every combination of $\lambda_{\Delta}$ and $\lambda_b$. For each method and calibration, we summarize the three class-specific utilities by their minimum, maximum, and average. Under the monetary interpretation of $\lambda_{\Delta}$ and $\lambda_b$, these values represent the expected net benefit per offer in monetary units. 
}

\begin{table}[h]
	\centering
	\small
	\setlength{\tabcolsep}{4pt}
	\renewcommand{\arraystretch}{1.15}
	\caption{\blue{
		Absolute money-metric utility by method, calibration class, and utility-proxy parameters. Each cell reports the minimum--maximum across the three driver classes, with the average in parentheses.}
	}
	\label{tab:absolute-driver-utility}
	\blue{
	\begin{tabular}{ccccccccc}
		\toprule
		& & & \multicolumn{3}{c}{\textbf{CC calibration}} & \multicolumn{3}{c}{\textbf{OSO calibration}} \\
		\cmidrule(lr){4-6}\cmidrule(lr){7-9}
		$\boldsymbol{\lambda}_{\boldsymbol{\Delta}}$ & $\boldsymbol{\lambda}_b$ & \textbf{POD} & \textbf{Class 1} & \textbf{Class 2} & \textbf{Class 3} & \textbf{Class 1} & \textbf{Class 2} & \textbf{Class 3} \\
		\midrule
		0.15 & 0.75 & \makecell{7.97--10.52\\(9.60)} & \makecell{9.66--11.42\\(10.56)} & \makecell{9.56--11.30\\(10.44)} & \makecell{9.27--10.90\\(10.08)} & \makecell{0.82--9.66\\(5.39)} & \makecell{3.66--11.20\\(8.44)} & \makecell{10.90--12.69\\(12.03)} \\
		0.23 & 0.75 & \makecell{7.84--10.32\\(9.41)} & \makecell{9.46--11.22\\(10.36)} & \makecell{9.36--11.11\\(10.25)} & \makecell{9.08--10.73\\(9.90)} & \makecell{0.80--9.46\\(5.27)} & \makecell{3.59--10.99\\(8.28)} & \makecell{10.73--12.50\\(11.84)} \\
		0.30 & 0.75 & \makecell{7.72--10.14\\(9.25)} & \makecell{9.27--11.05\\(10.18)} & \makecell{9.19--10.95\\(10.08)} & \makecell{8.92--10.58\\(9.75)} & \makecell{0.78--9.27\\(5.17)} & \makecell{3.53--10.80\\(8.14)} & \makecell{10.58--12.33\\(11.68)} \\
		\midrule
		0.15 & 1.25 & \makecell{6.54--8.30\\(7.58)} & \makecell{7.39--9.37\\(8.38)} & \makecell{7.29--9.25\\(8.26)} & \makecell{7.02--8.89\\(7.94)} & \makecell{0.63--7.39\\(4.12)} & \makecell{2.89--8.82\\(6.65)} & \makecell{8.89--10.33\\(9.80)} \\
		0.23 & 1.25 & \makecell{6.41--8.10\\(7.39)} & \makecell{7.18--9.17\\(8.17)} & \makecell{7.09--9.06\\(8.06)} & \makecell{6.84--8.72\\(7.76)} & \makecell{0.61--7.18\\(4.00)} & \makecell{2.82--8.61\\(6.49)} & \makecell{8.72--10.14\\(9.62)} \\
		0.30 & 1.25 & \makecell{6.30--7.92\\(7.23)} & \makecell{7.00--9.00\\(7.99)} & \makecell{6.91--8.90\\(7.89)} & \makecell{6.68--8.57\\(7.60)} & \makecell{0.59--7.00\\(3.90)} & \makecell{2.76--8.42\\(6.35)} & \makecell{8.57--9.98\\(9.46)} \\
		\midrule
		0.15 & 2.00 & \makecell{4.29--4.97\\(4.55)} & \makecell{3.98--6.29\\(5.10)} & \makecell{3.87--6.18\\(4.99)} & \makecell{3.65--5.88\\(4.72)} & \makecell{0.34--3.98\\(2.22)} & \makecell{1.73--5.25\\(3.96)} & \makecell{5.88--6.80\\(6.46)} \\
		0.23 & 2.00 & \makecell{4.06--4.77\\(4.36)} & \makecell{3.77--6.09\\(4.90)} & \makecell{3.67--5.99\\(4.79)} & \makecell{3.47--5.71\\(4.54)} & \makecell{0.32--3.77\\(2.10)} & \makecell{1.66--5.04\\(3.80)} & \makecell{5.71--6.61\\(6.28)} \\
		0.30 & 2.00 & \makecell{3.85--4.59\\(4.20)} & \makecell{3.59--5.92\\(4.72)} & \makecell{3.50--5.82\\(4.62)} & \makecell{3.31--5.56\\(4.39)} & \makecell{0.31--3.59\\(2.00)} & \makecell{1.60--4.86\\(3.67)} & \makecell{5.56--6.44\\(6.11)} \\
		\bottomrule
	\end{tabular}
	}
\end{table}

\blue{
At the central specification $(\lambda_{\Delta},\lambda_b)=(0.23,1.25)$, POD produces class-specific expected utilities between 6.41 and 8.10, with an average of 7.39. The three CC calibrations yield higher average utilities, ranging from 7.76 to 8.17. This difference is driven primarily by a redistribution toward Class~3. Relative to POD, CC substantially increases the utility of Class~3 while slightly reducing that of Class~1 and, for some calibrations, Class~2. Thus, POD's operator-performance advantage does not arise from maximizing absolute driver utility.

The OSO results illustrate the importance of the assumed representative class. OSO calibrated to Class~3 produces the highest absolute driver utilities at the central specification, ranging from 8.72 to 10.14 with an average of 9.62. Its relatively generous offers benefit drivers but result in lower operator performance than POD. In contrast, OSO calibrated to Class~1 produces a minimum utility of only 0.61 and an average of 4.00. OSO calibrated to Class~2 is intermediate, but its minimum of 2.82 still indicates a substantial miscalibration cost for the least compatible driver class. Consequently, average utility alone can conceal very poor outcomes for an individual class, making the minimum an important complementary measure.

Increasing the monetary value assigned to bundle handling lowers all absolute utility levels mechanically. For example, at $\lambda_{\Delta}=0.23$, POD's average utility decreases from 9.41 for $\lambda_b=0.75$ to 7.39 for $\lambda_b=1.25$ and 4.36 for $\lambda_b=2.00$. Changes in $\lambda_{\Delta}$ have a smaller effect. At the highest bundle-handling valuation, POD has a lower average than CC but a higher minimum and a substantially smaller cross-class range. For example, at $(\lambda_{\Delta},\lambda_b)=(0.23,2.00)$, POD has a minimum of 4.06, compared with 3.47--3.77 for the CC calibrations, while the CC averages remain higher. Hence, POD's distributional advantage at high effort valuations reflects a more compressed utility distribution, whereas CC and OSO calibrated to Class~3 can provide greater average driver surplus. Overall, the results reveal a trade-off among operator performance, absolute driver utility, and the distribution of utility across behavioral classes.
}

{\color{black}
\raggedbottom
\section{Validation with externally estimated and perturbed coefficients}
\label{sec:stress-test}

The behavioral classes in the main experiments are stylized profiles designed to isolate systematic differences in effort and compensation sensitivity. To complement this controlled analysis, we conduct validation tests using independently estimated coefficients from \citet{mohri_modeling_2024} (\cref{sec:mohri_test}), coefficients fitted from an operational Meituan dataset released for the INFORMS TSL Data-Driven Research Challenge \citep{wang2025meituan} (\cref{sec:meituan}), and systematically perturbed class coefficients (\cref{sec:perturbed}). 

\subsection{Coefficients based on \citet{mohri_modeling_2024}} \label{sec:mohri_test}

The latent-class choice model of \citet{mohri_modeling_2024} distinguishes leisurely (LC1), avid (LC2), and skeptical (LC3) public-transport users across detour, weight and incentive coefficients, which are reported in \cref{tab:mohri-coefficients}. We use their weight coefficients as a proxy to the bundle size coefficients. As also recognized in their paper, the incentive coefficients of LC1 and LC2 are not significantly different from zero. For these classes, the optimal compensations carry the boundary compensation $C^*_{wb}=0$. This outcome should not be interpreted as evidence that these users would perform deliveries without compensation. The choice experiment of \citet{mohri_modeling_2024} considers only positive incentives and therefore did not identify behavior at zero compensation. 
Hence, we do not use these solutions and report those results where all drivers receive the LC3 coefficients in the remainder. We use the same single-depot benchmark instances used in the main study.

\begin{table}[H]
	\color{black}
	\centering
	\small
	\setlength{\tabcolsep}{4.5pt}
	\renewcommand{\arraystretch}{1.12}
	\caption{Class-specific acceptance-model estimates from \citet{mohri_modeling_2024}.}
	\label{tab:mohri-coefficients}
	\begin{tabular}{llccccc}
		\toprule
		\textbf{Class} & \textbf{Label} & \textbf{Share} & \textbf{Intercept} & \textbf{Detour} & \textbf{Weight} & \textbf{Incentive} \\
		\midrule
		LC1 & Leisurely & 19\% &
		2.042 &
		$-3.120$ &
		$-0.452$ &
		0.054 \\
		LC2 & Avid & 53\% &
		5.671 &
		$-2.846$ &
		$-0.283$ &
		0.080 \\
		LC3 & Skeptical & 28\% &
		$-3.908$ &
		$-4.893$ &
		$-0.355$ &
		1.616 \\
		\bottomrule
	\end{tabular}
\end{table}


\cref{fig:mohri-task-sensitivity} shows how outcomes vary with task availability. As the number of tasks increases from 30 to 120, mean acceptance rises from 0.86 to 0.95 and mean bundle size from 3.7 to 4.9, while mean detour decreases from 2.00 to 1.50. Mean compensation initially increases from 10.59 to 11.08 before declining to 9.92. Although the values differ, these trends closely mirror those discussed in \cref{sec:sensitivity:tasks}. In particular, \cref{fig:mohri-task-compensation} reproduces the nonmonotonic compensation pattern observed for the stylized Class~3 profile in \cref{fig:task_compensation}.



\begin{figure}[H]
	\color{black}
	\captionsetup[subfigure]{aboveskip=-0.5pt,belowskip=-0.5pt}
	\begin{subfigure}[t]{0.49\textwidth}
		\centering
		\includegraphics[width=0.9\textwidth]{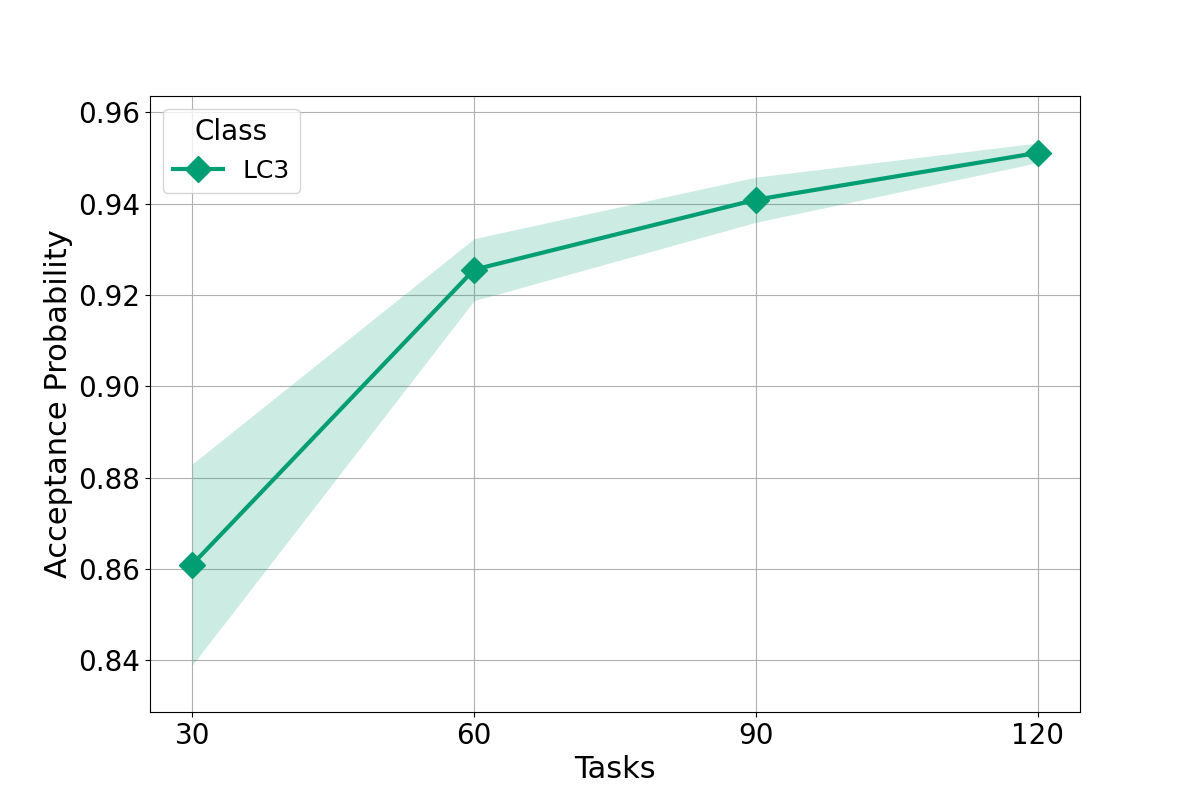}
		\caption{Acceptance probability.}
		\label{fig:mohri-task-acceptance}
	\end{subfigure}\hfill
	\begin{subfigure}[t]{0.49\textwidth}
		\centering
		\includegraphics[width=0.9\textwidth]{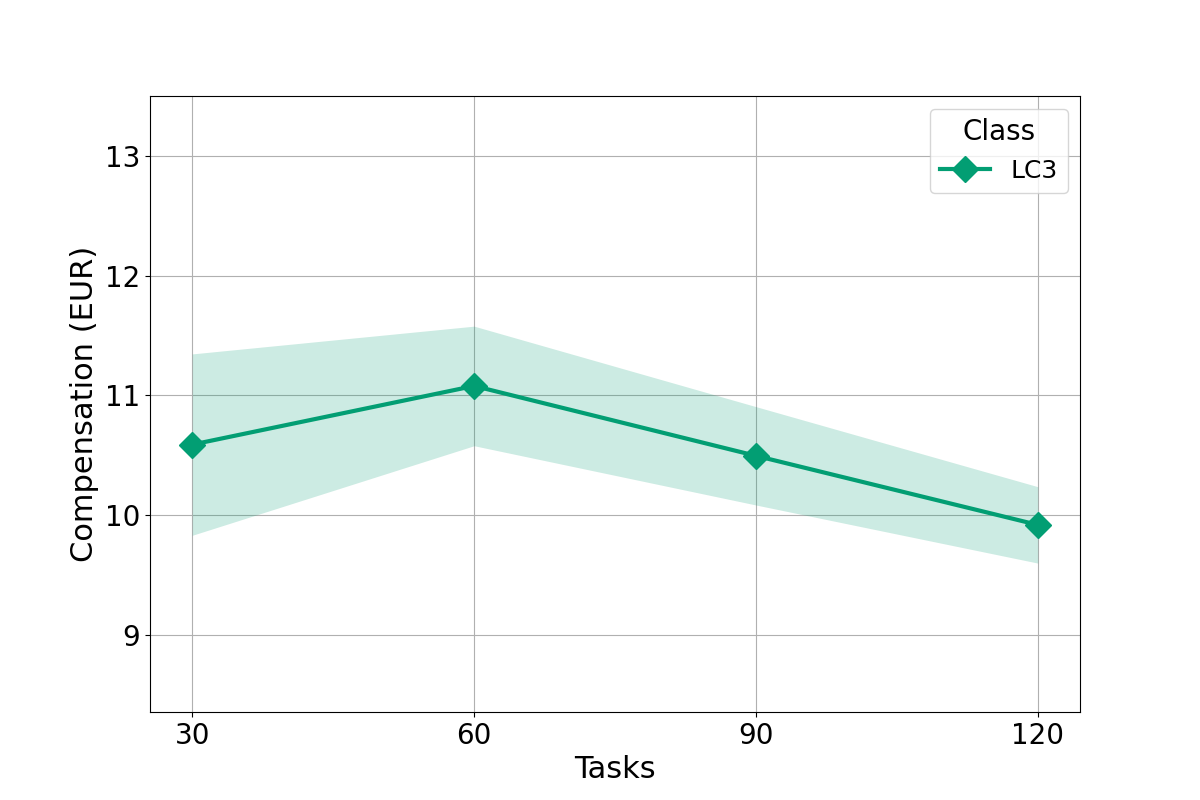}
		\caption{Compensation.}
		\label{fig:mohri-task-compensation}
	\end{subfigure}
	\begin{subfigure}[t]{0.49\textwidth}
		\centering
		\includegraphics[width=0.9\textwidth]{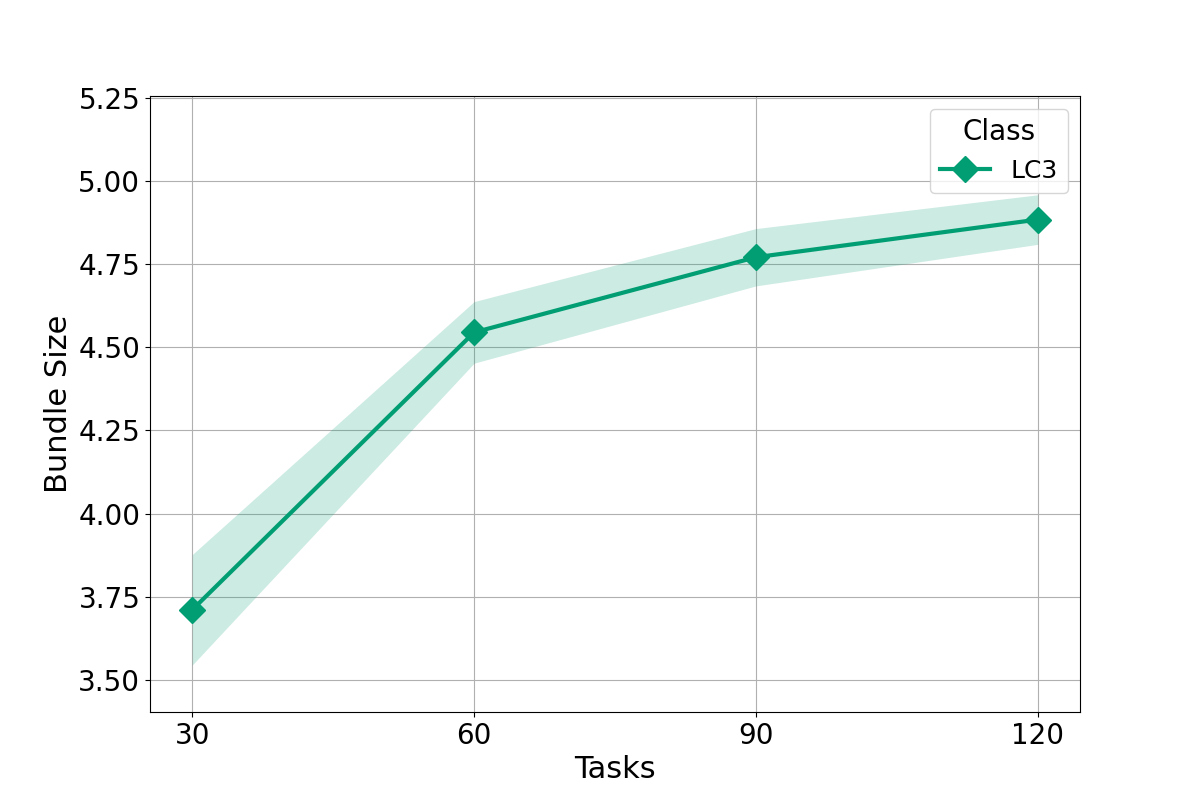}
		\caption{Bundle size.}
		\label{fig:mohri-task-bundle}
	\end{subfigure}\hfill
	\begin{subfigure}[t]{0.49\textwidth}
		\centering
		\includegraphics[width=0.9\textwidth]{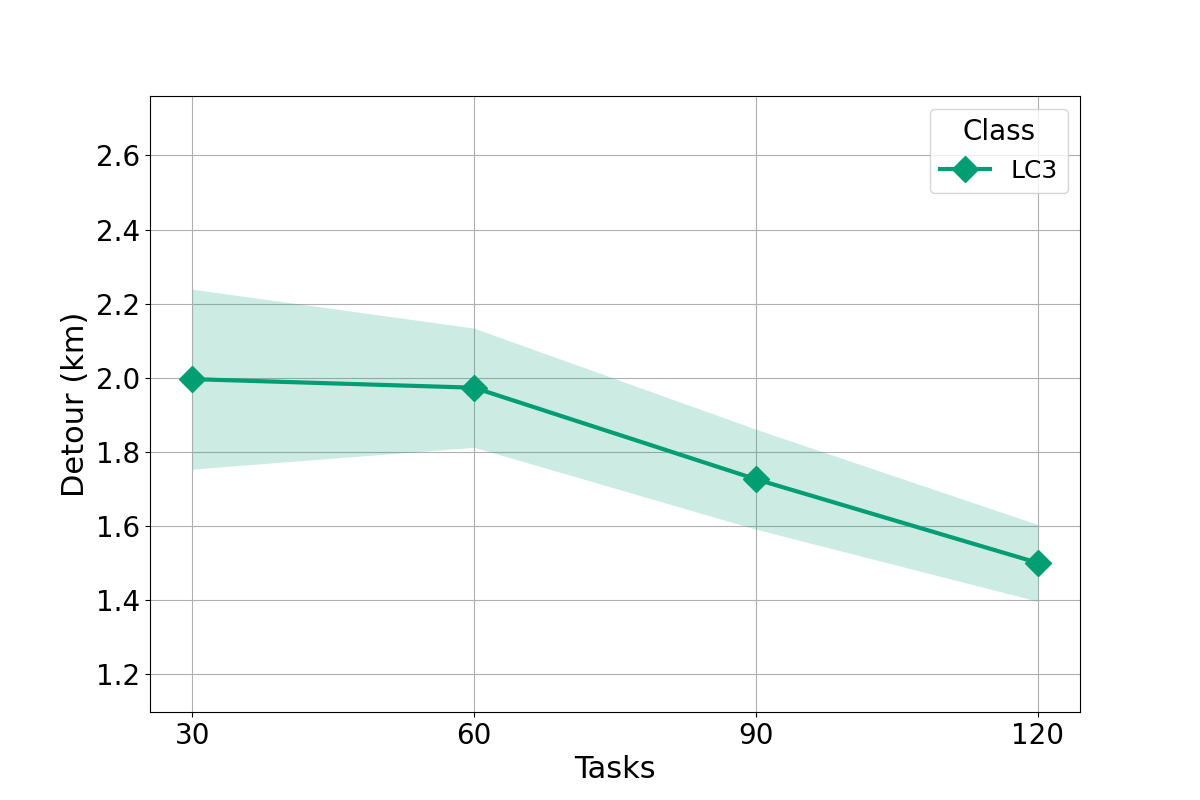}
		\caption{Detour.}
		\label{fig:mohri-task-detour}
	\end{subfigure}
	\caption{Mean LC3 offer characteristics and their 95\% confidence intervals by number of tasks.}
	\label{fig:mohri-task-sensitivity}
\end{figure}

Similarly, the trends with increasing driver availability mirror those in \cref{sec:sensitivity:drivers}, as illustrated in \cref{fig:mohri-driver-sensitivity}. As the driver-to-task ratio increases from 0.1 to 0.5, mean compensation decreases from 12.47 to 9.15, detour decreases from 2.40 to 1.36, and bundle size decreases from 4.95 to 4.25. Meanwhile, acceptance remains comparatively stable between 0.912 and 0.924, mimicing Classes~1 and 2 in \cref{fig:driver_acceptance}. 

\begin{figure}[H]
	\color{black}
	\captionsetup[subfigure]{aboveskip=-0.5pt,belowskip=-0.5pt}
	\begin{subfigure}[t]{0.49\textwidth}
		\centering
		\includegraphics[width=0.9\textwidth]{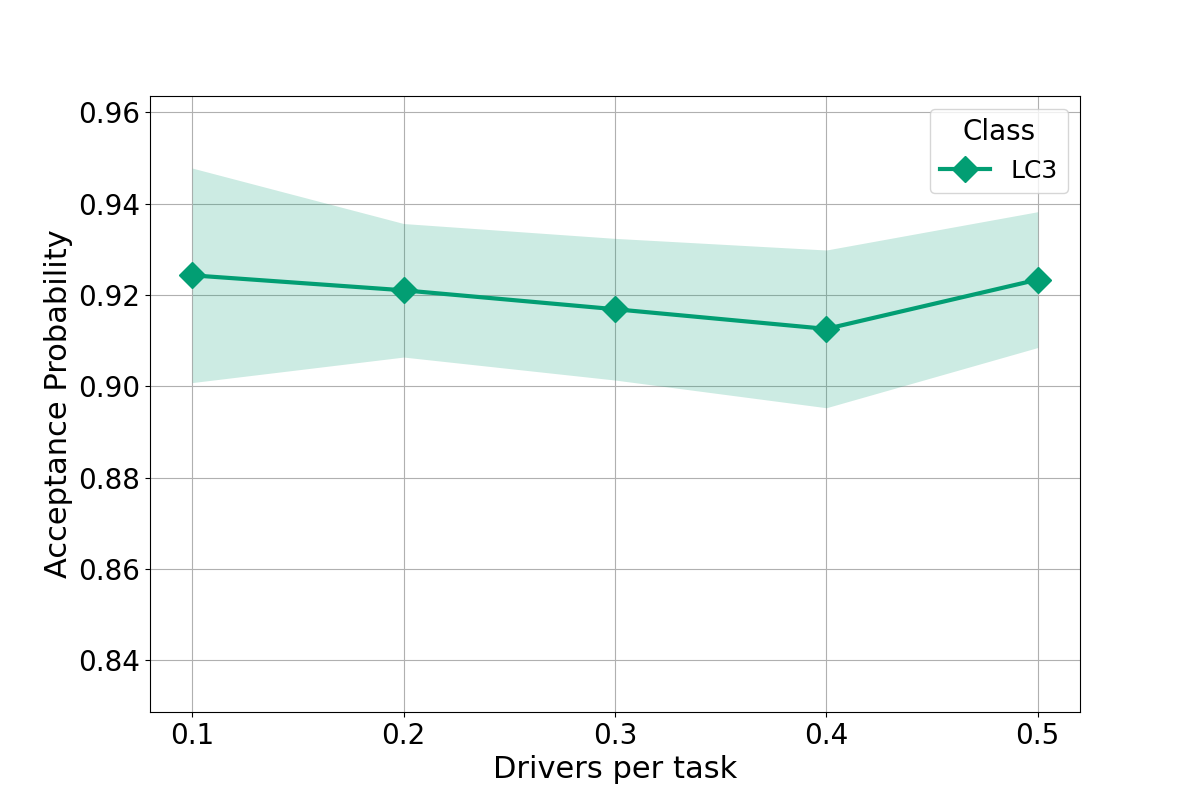}
		\caption{Acceptance probability.}
		\label{fig:mohri-driver-acceptance}
	\end{subfigure}\hfill
	\begin{subfigure}[t]{0.49\textwidth}
		\centering
		\includegraphics[width=0.9\textwidth]{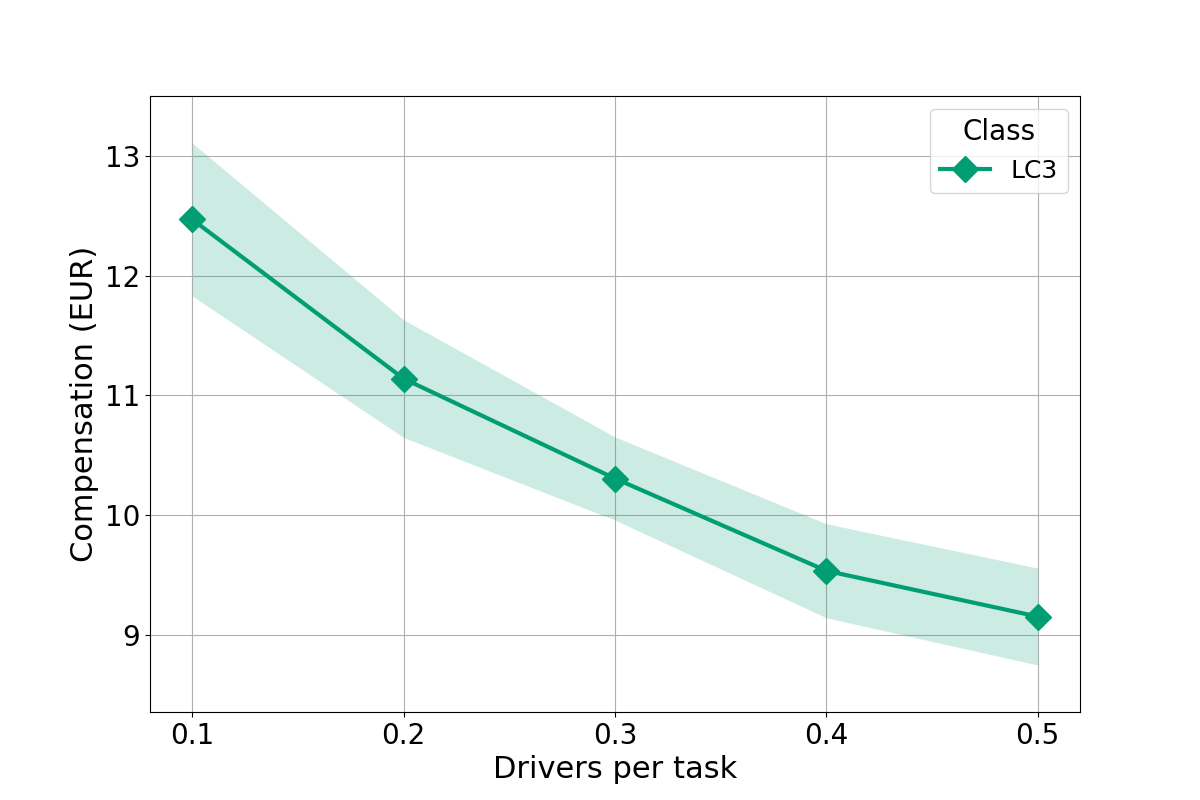}
		\caption{Compensation.}
		\label{fig:mohri-driver-compensation}
	\end{subfigure}
	\begin{subfigure}[t]{0.49\textwidth}
		\centering
		\includegraphics[width=0.9\textwidth]{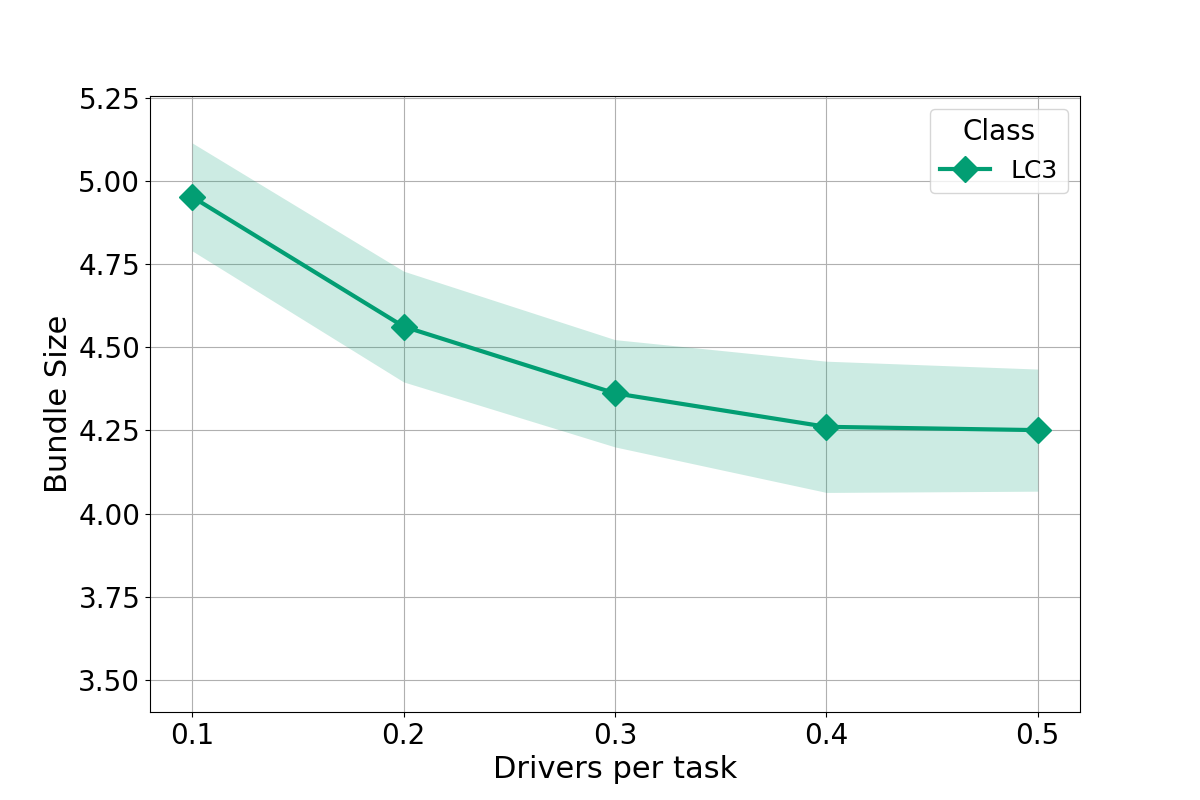}
		\caption{Bundle size.}
		\label{fig:mohri-driver-bundle}
	\end{subfigure}\hfill
	\begin{subfigure}[t]{0.49\textwidth}
		\centering
		\includegraphics[width=0.9\textwidth]{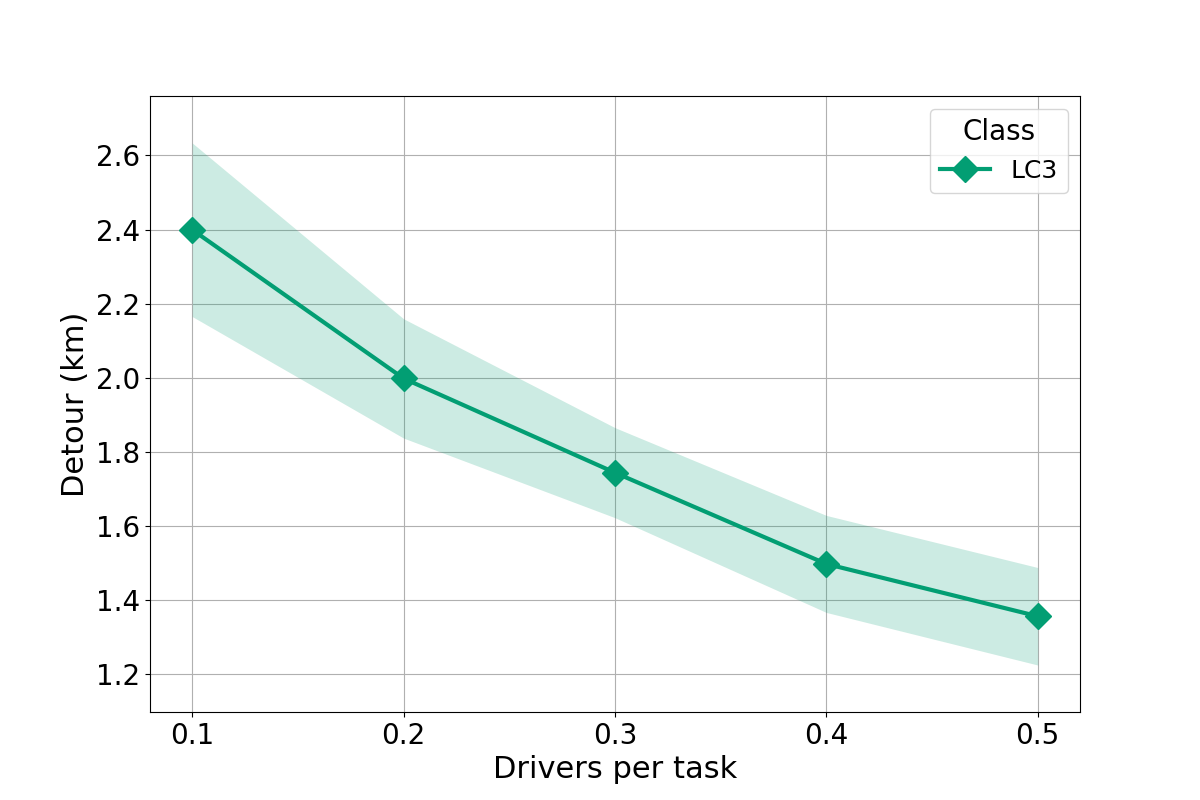}
		\caption{Detour.}
		\label{fig:mohri-driver-detour}
	\end{subfigure}
	\caption{Mean LC3 offer characteristics and their 95\% confidence intervals by driver-to-task ratio.}
	\label{fig:mohri-driver-sensitivity}
\end{figure}

\cref{tab:mohri-computational-results,tab:mohri-heuristic-gaps} show the runtime performances of H-C, H-DDC and E-DDC, and the optimality gaps for H-C and H-DDC, respectively. These results largely support the findings in \cref{sec:results:exact,sec:results:heuristic}. The similarities include that the runtimes increase with task and driver availability, H-C is consistently the fastest method, and H-DDC maintains small optimality gaps while H-C improves as tasks and drivers become more abundant. The main difference is that the externally calibrated LC3 instances are substantially easier to solve. E-DDC solves all instances, including the largest ones, much faster than those in the main experiments. Furthermore, solution quality remains comparable. H-DDC gaps never exceed 0.36\%, while H-C exhibits the same stronger sensitivity to instance size and driver density observed in \cref{tab:corridor_heuristic_gaps}.

\begin{table}[H]
	\color{black}
	\centering
	\small
	\setlength{\tabcolsep}{8pt}
	\renewcommand{\arraystretch}{1.12}
	\caption{Mean runtimes in seconds under the LC3 specification across varying instance sizes.}
	\label{tab:mohri-computational-results}
	\resizebox{\textwidth}{!}{%
	\begin{tabular}{c ccccc ccccc ccccc}
		\toprule
		\multirow{3}{*}{\textbf{Tasks}}
		& \multicolumn{5}{c}{\textbf{H-C}}
		& \multicolumn{5}{c}{\textbf{H-DDC}}
		& \multicolumn{5}{c}{\textbf{E-DDC}} \\
		\cmidrule(lr){2-6}\cmidrule(lr){7-11}\cmidrule(lr){12-16}
		& \multicolumn{5}{c}{\textbf{Driver-to-task ratio}}
		& \multicolumn{5}{c}{\textbf{Driver-to-task ratio}}
		& \multicolumn{5}{c}{\textbf{Driver-to-task ratio}} \\
		\cmidrule(lr){2-6}\cmidrule(lr){7-11}\cmidrule(lr){12-16}
		& \textbf{0.1} & \textbf{0.2} & \textbf{0.3} & \textbf{0.4} & \textbf{0.5}
		& \textbf{0.1} & \textbf{0.2} & \textbf{0.3} & \textbf{0.4} & \textbf{0.5}
		& \textbf{0.1} & \textbf{0.2} & \textbf{0.3} & \textbf{0.4} & \textbf{0.5} \\
		\midrule
		30
		& 0.00 & 0.01 & 0.01 & 0.17 & 0.19
		& 0.01 & 0.02 & 0.05 & 0.21 & 0.23
		& 0.01 & 0.03 & 0.07 & 0.23 & 0.25 \\
		60
		& 0.02 & 0.07 & 0.30 & 0.51 & 0.34
		& 0.10 & 0.35 & 0.87 & 1.09 & 0.93
		& 0.14 & 0.47 & 1.18 & 1.34 & 1.23 \\
		90
		& 0.43 & 1.00 & 1.53 & 1.84 & 2.05
		& 1.72 & 4.33 & 5.01 & 5.71 & 4.98
		& 2.46 & 7.31 & 8.96 & 10.57 & 9.23 \\
		120
		& 1.43 & 4.42 & 6.77 & 8.15 & 10.05
		& 6.02 & 18.18 & 19.08 & 19.39 & 20.98
		& 9.42 & 96.70 & 106.95 & 99.65 & 90.69 \\
		\bottomrule
	\end{tabular}
	}
\end{table}

\begin{table}[H]
	\color{black}
	\centering
	\small
	\setlength{\tabcolsep}{8pt}
	\renewcommand{\arraystretch}{1.12}
	\caption{Mean gaps to the proven-optimal E-DDC objective under the LC3 specification (in \%) for different driver-to-task ratios.}
	\label{tab:mohri-heuristic-gaps}
	\begin{tabular}{c ccccc ccccc}
		\toprule
		\multirow{3}{*}{\textbf{Tasks}}
		& \multicolumn{5}{c}{\textbf{H-C}}
		& \multicolumn{5}{c}{\textbf{H-DDC}} \\
		\cmidrule(lr){2-6}\cmidrule(lr){7-11}
		& \multicolumn{5}{c}{\textbf{Driver-to-task ratio}}
		& \multicolumn{5}{c}{\textbf{Driver-to-task ratio}} \\
		\cmidrule(lr){2-6}\cmidrule(lr){7-11}
		& \textbf{0.1} & \textbf{0.2} & \textbf{0.3} & \textbf{0.4} & \textbf{0.5}
		& \textbf{0.1} & \textbf{0.2} & \textbf{0.3} & \textbf{0.4} & \textbf{0.5} \\
		\midrule
		30
		& 28.45 & 16.69 & 16.05 & 9.72 & 7.23
		& 0.00 & 0.00 & 0.01 & 0.00 & 0.00 \\
		60
		& 13.17 & 6.35 & 3.03 & 2.51 & 1.49
		& 0.06 & 0.10 & 0.07 & 0.15 & 0.03 \\
		90
		& 9.41 & 2.59 & 1.36 & 0.73 & 0.62
		& 0.12 & 0.23 & 0.11 & 0.24 & 0.25 \\
		120
		& 4.14 & 1.40 & 1.11 & 0.29 & 0.35
		& 0.17 & 0.36 & 0.35 & 0.24 & 0.22 \\
		\bottomrule
	\end{tabular}
\end{table}

\subsection{Coefficients based on INFORMS TSL data-driven research challenge dataset} \label{sec:meituan}

The operational dataset released for the INFORMS TSL Data-Driven Research Challenge \citep{wang2025meituan} records the order-dispatching process of Meituan, the largest Chinese on-demand food-delivery platform: 654{,}343 individual delivery offers (waybills) placed to 4{,}955 couriers over an eight-day window (17--24 October 2022), each carrying pickup and delivery coordinates, timestamps, and a binary grab flag recording whether the courier accepted the offer. The dataset does not, however, directly contain the three quantities that our acceptance-model~\eqref{eq:acceptance_model} requires. Compensation is absent entirely, as the data does not contain compensation information and Meituan does not publicly announce their compensation calculation approach. Bundles in our all-or-nothing sense do not exist either, instead the platform assigns, and lets couriers refuse, each waybill \emph{individually}, so a batch does not have to be accepted or rejected as single unit. Detour is likewise undefined, as Meituan couriers are dedicated gig workers (a segment we explicitly exclude from our scope) who follow no preplanned personal trip against which an incremental distance could be measured. We therefore use proxies for all three attributes. Bundle size is proxied by the number of offers made to a courier at the same dispatch instant, detour by the individual delivery distance of the offered orders, and compensation by a piecewise-linear function of bundle size and detour, which is calibrated to the reported average Meituan pay and bonus levels \citep{wu2022framework,sun2025crowdsourced}. Because all three regressors are proxies rather than observed variables, and the resulting data are transformed substantially, partitioning couriers into distinct behavioral classes along these reconstructed axes is not realistic.

This is reinforced by the couriers' acceptance behavior, which shows little of the heterogeneity a clustering would exploit. Acceptance is near-uniform at 89.1\%, 30.7\% of couriers never decline a single offer over the eight days, and the median courier registers only four rejections in total. We therefore regard a single \emph{pooled} estimate as the more realistic specification, and fit one binary logit by maximum likelihood over all couriers, mapping its intercept, detour, bundle-size, and compensation coefficients onto $\alpha$, $\beta_1$, $\beta_2$, and $\gamma$ of the acceptance-model~\eqref{eq:acceptance_model}. \cref{tab:meituan-coefficients} reports the resulting estimate, which carries the expected signs (effort deters acceptance, compensation encourages it) and is applied to every driver in the benchmark. The large positive intercept is consistent with the near-uniform acceptance observed in the data, and the low sensitivity to compensation is partly explained by collinearity. The compensation proxy is itself a function of bundle size and detour, which reduces its explanatory power.

\begin{table}[H]
	\color{black}
	\centering
	\small
	\setlength{\tabcolsep}{6pt}
	\renewcommand{\arraystretch}{1.12}
	\caption{Pooled acceptance-model estimate fitted on the INFORMS TSL Meituan dataset.}
	\label{tab:meituan-coefficients}
	\begin{tabular}{lcccc}
		\toprule
		\textbf{Class} & \(\alpha\) (Intercept) & \(\beta_1\) (Detour) & \(\beta_2\) (Bundle Size) & \(\gamma\) (Compensation) \\
		\midrule
		Meituan & 3.2633 & $-0.6977$ & $-2.2453$ & 0.5153 \\
		\bottomrule
	\end{tabular}
\end{table}

\cref{fig:meituan-task-sensitivity} reports the mean characteristics of selected offers as task availability varies. As the number of tasks increase, the mean acceptance probability rises from 0.66 to 0.69, mean compensation from 2.36 to 6.31, and mean bundle size from 1.51 to 2.43, while reducing the mean detour from 0.585 to 0.329. The trends in acceptance probability, bundle size, and detour resemble those for Class~3 in \cref{fig:task_acceptance,fig:task_bundle,fig:task_detour}. The lower absolute acceptance levels compared to the values reported in \cref{sec:sensitivity} could be explained by the low sensitivity to compensation. Compensation itself, however, increases monotonically with task availability. One possible explanation is that greater task availability enables the model to construct larger bundles with shorter detours. The resulting efficiency gains permit higher compensation at larger instances, as illustrated in \cref{fig:meituan-task-compensation}. 

\cref{fig:meituan-driver-sensitivity} shows the mean characteristics of the selected offers across increasing driver availability. As the driver-to-task ratio rises from 0.1 to 0.5, mean compensation decreases from 7.92 to 2.02 and mean bundle size decreases from 2.78 to 1.45. Mean acceptance probability remains nearly constant, ranging from 0.679 to 0.683 across the five driver-to-task ratios. Mean detour decreases from 0.472 at a ratio of 0.1 to its minimum of 0.397 at a ratio of 0.4, before increasing slightly to 0.412 at a ratio of 0.5. The trends in acceptance probability, compensation, and bundle size mirror those observed for Classes~1 and 2 in \cref{fig:driver_acceptance,fig:driver_compensation,fig:driver_bundle}. The flat detour values, by contrast, likely reflect the fact that absolute detour levels are already low, leaving little room for further reduction.

\begin{figure}[H]
	\color{black}
	\captionsetup[subfigure]{aboveskip=-0.5pt,belowskip=-0.5pt}
	\begin{subfigure}[t]{0.49\textwidth}
		\centering
		\includegraphics[width=0.9\textwidth]{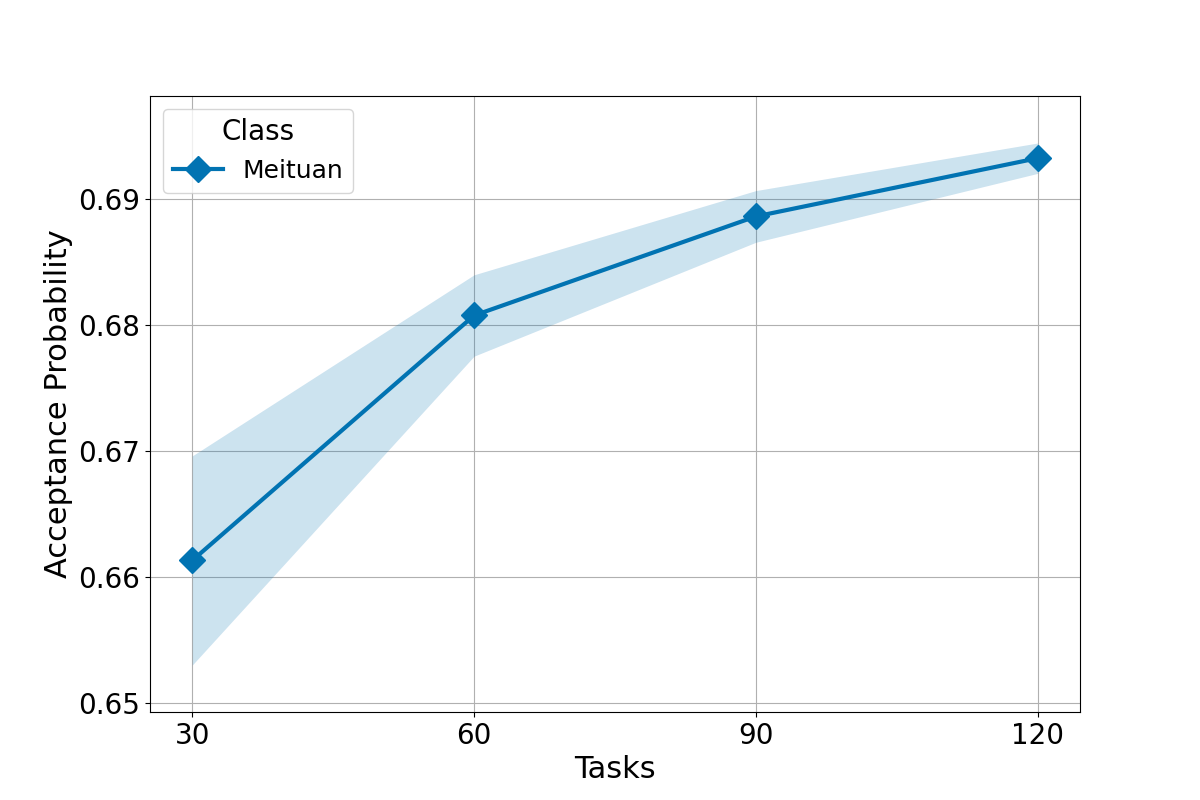}
		\caption{Acceptance probability.}
		\label{fig:meituan-task-acceptance}
	\end{subfigure}\hfill
	\begin{subfigure}[t]{0.49\textwidth}
		\centering
		\includegraphics[width=0.9\textwidth]{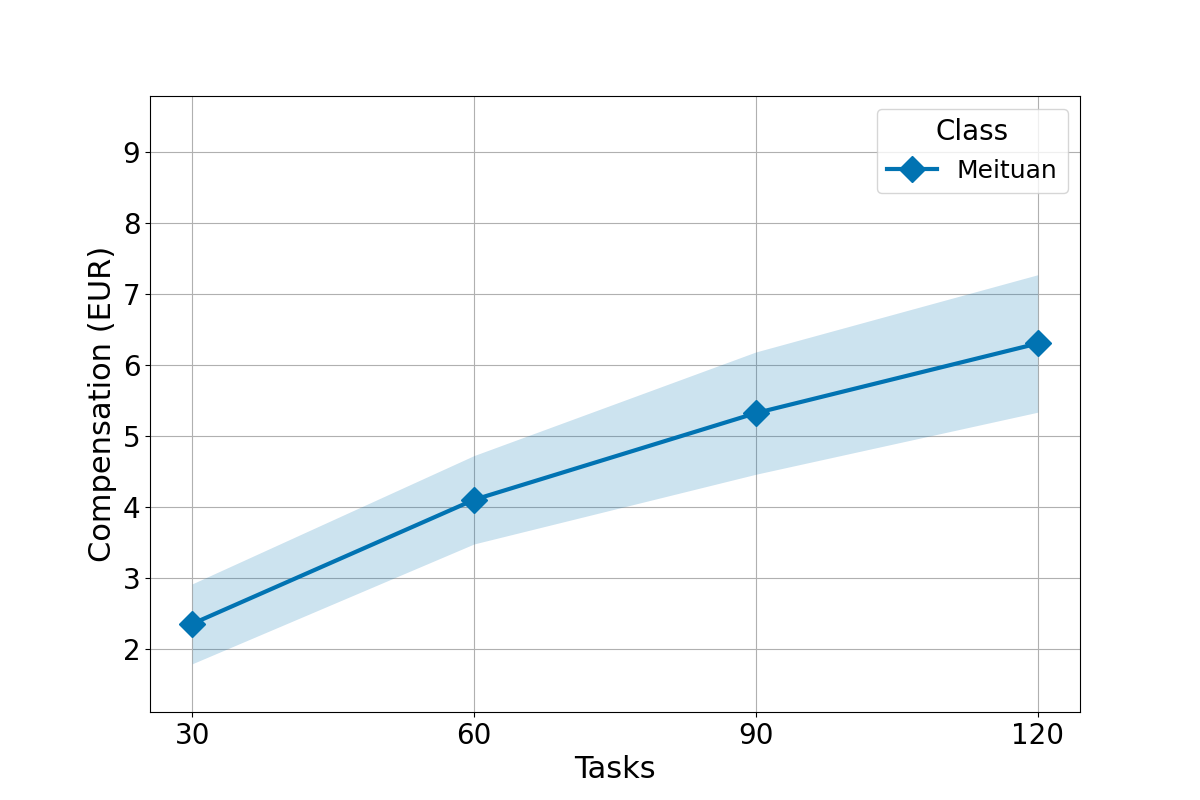}
		\caption{Compensation.}
		\label{fig:meituan-task-compensation}
	\end{subfigure}
	\begin{subfigure}[t]{0.49\textwidth}
		\centering
		\includegraphics[width=0.9\textwidth]{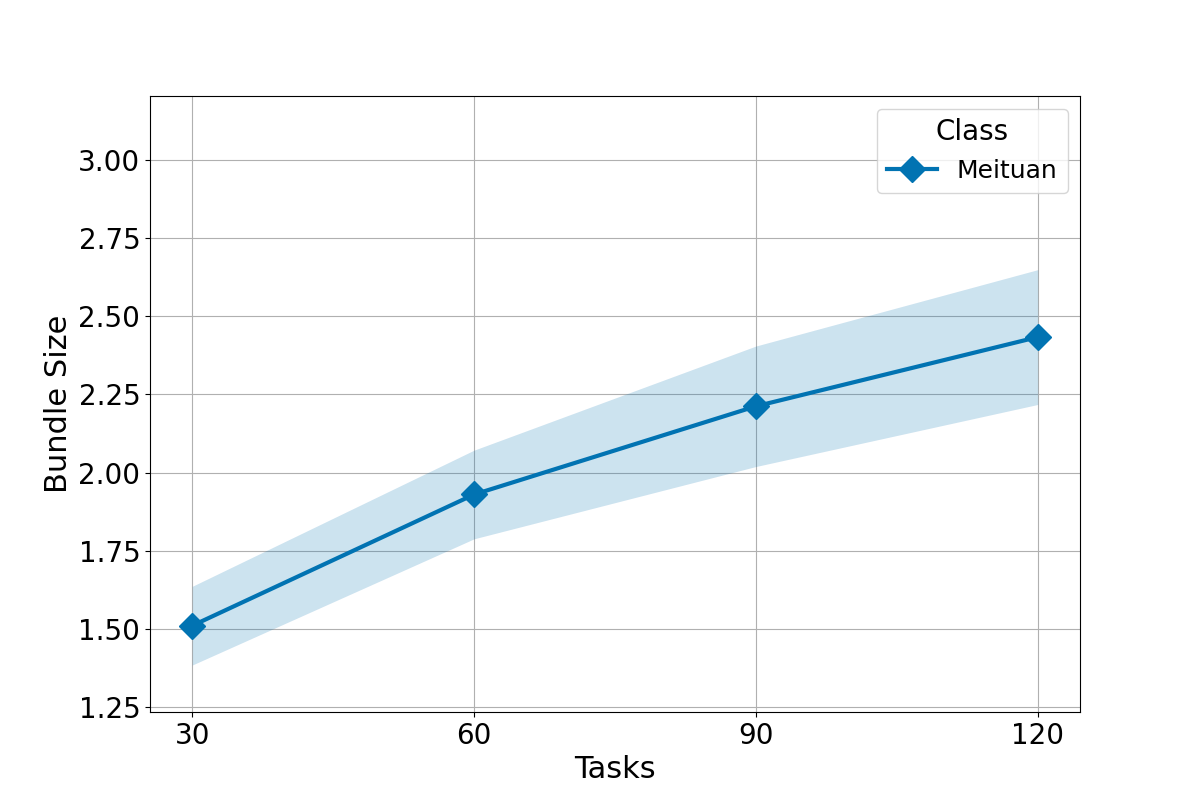}
		\caption{Bundle size.}
		\label{fig:meituan-task-bundle}
	\end{subfigure}\hfill
	\begin{subfigure}[t]{0.49\textwidth}
		\centering
		\includegraphics[width=0.9\textwidth]{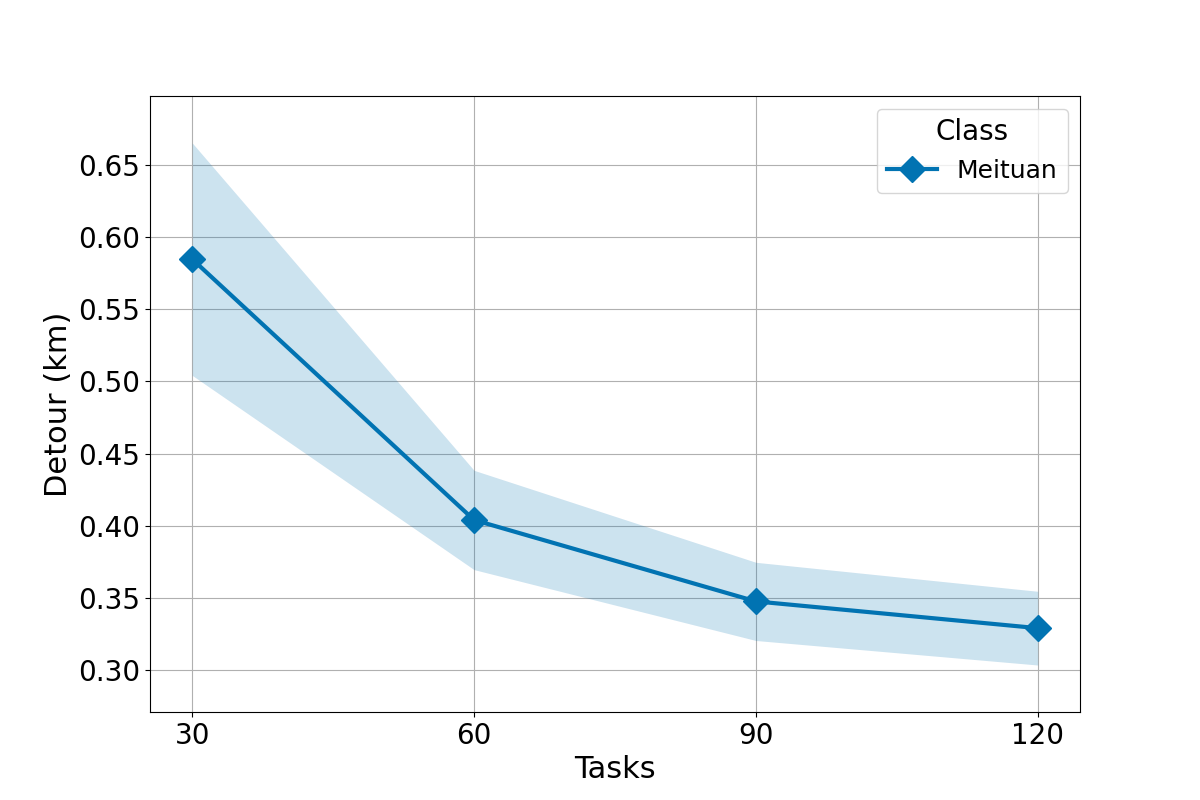}
		\caption{Detour.}
		\label{fig:meituan-task-detour}
	\end{subfigure}
	\caption{Mean offer characteristics and their 95\% confidence intervals under the fitted Meituan specification by number of tasks.}
	\label{fig:meituan-task-sensitivity}
	\vspace{-2em}
\end{figure}

\begin{figure}[H]
	\vspace{-1em}
	\color{black}
	\captionsetup[subfigure]{aboveskip=-0.5pt,belowskip=-0.5pt}
	\begin{subfigure}[t]{0.49\textwidth}
		\centering
		\includegraphics[width=0.9\textwidth]{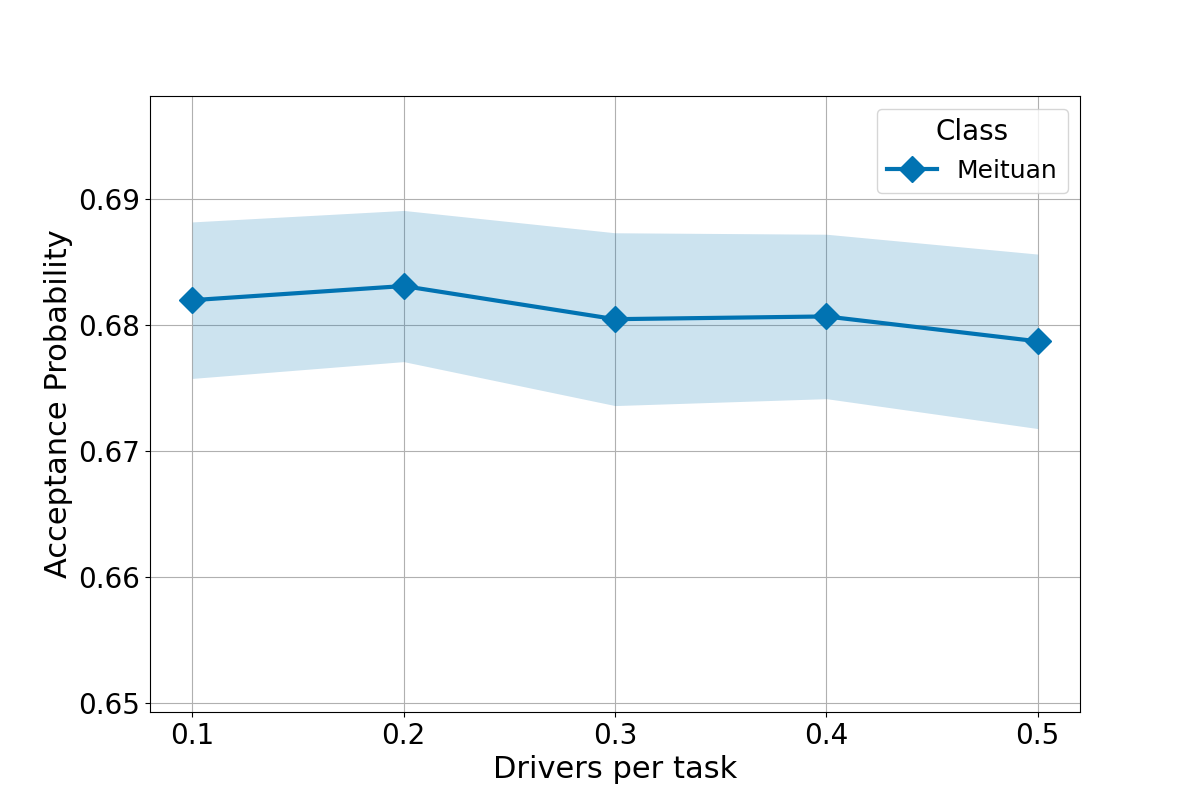}
		\caption{Acceptance probability.}
		\label{fig:meituan-driver-acceptance}
	\end{subfigure}\hfill
	\begin{subfigure}[t]{0.49\textwidth}
		\centering
		\includegraphics[width=0.9\textwidth]{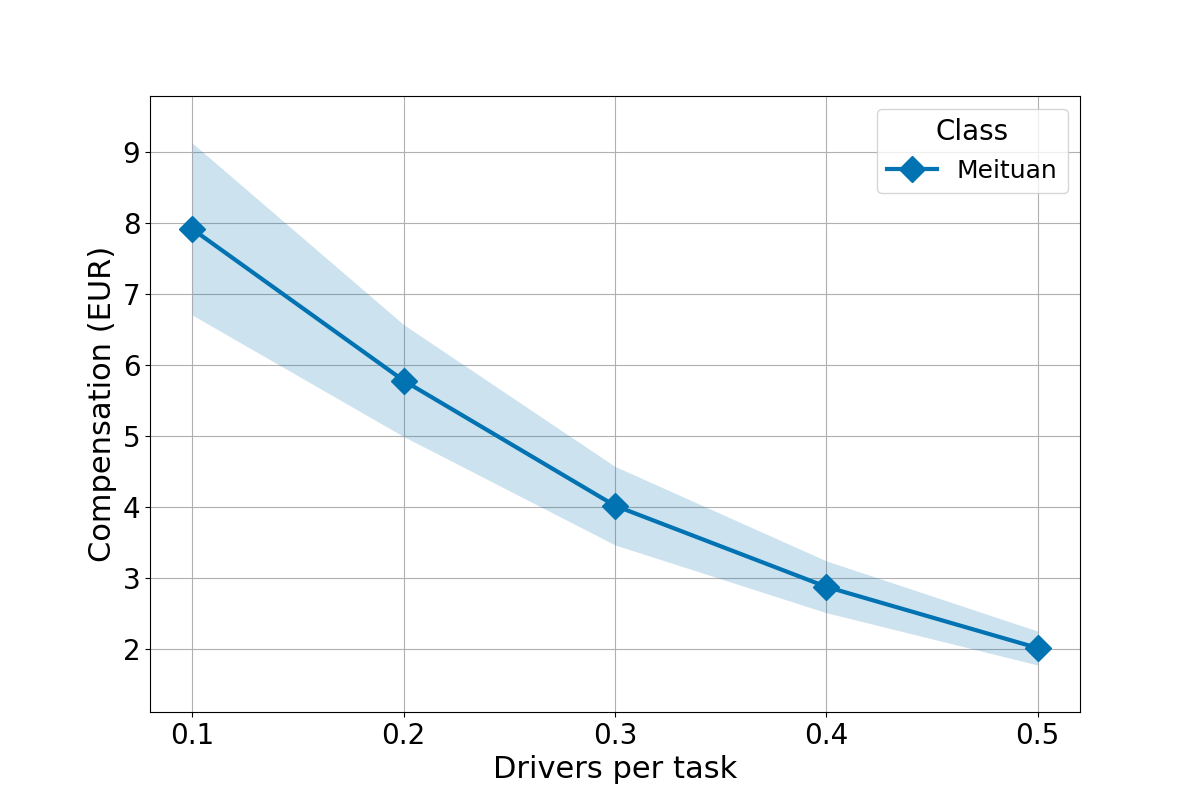}
		\caption{Compensation.}
		\label{fig:meituan-driver-compensation}
	\end{subfigure}
	\begin{subfigure}[t]{0.49\textwidth}
		\centering
		\includegraphics[width=0.9\textwidth]{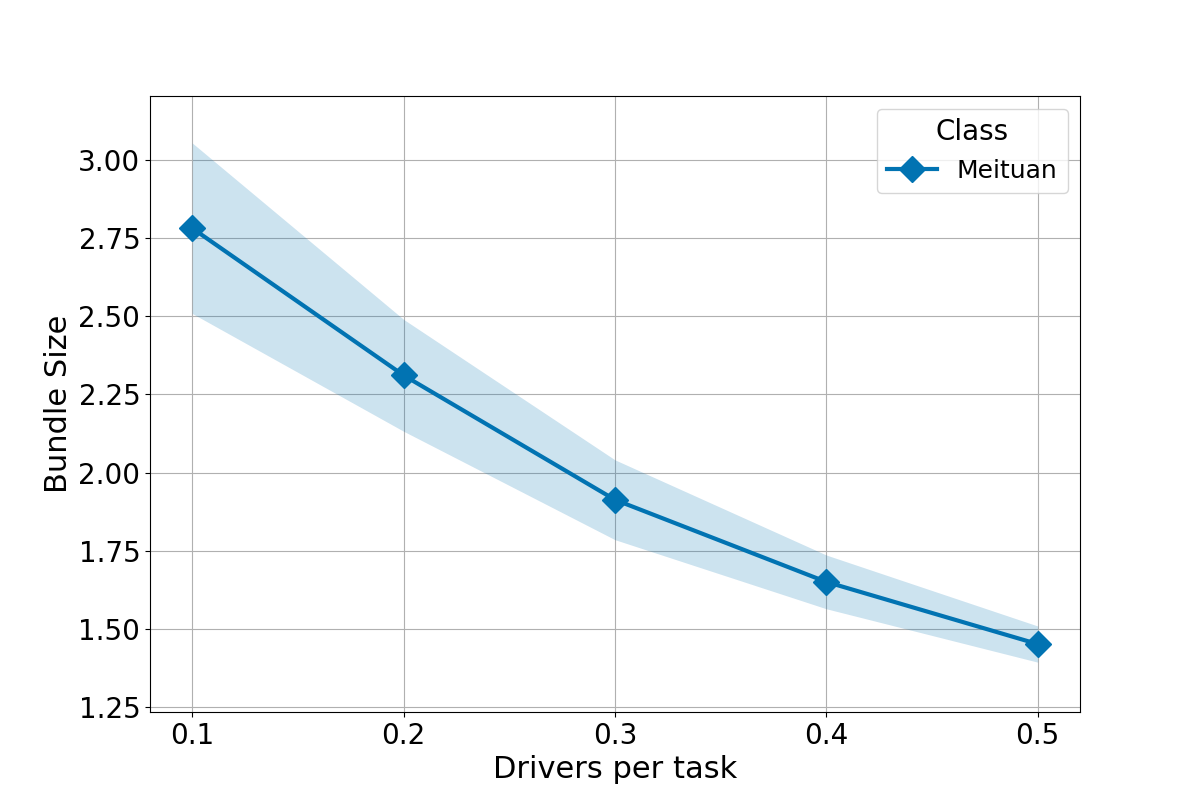}
		\caption{Bundle size.}
		\label{fig:meituan-driver-bundle}
	\end{subfigure}\hfill
	\begin{subfigure}[t]{0.49\textwidth}
		\centering
		\includegraphics[width=0.9\textwidth]{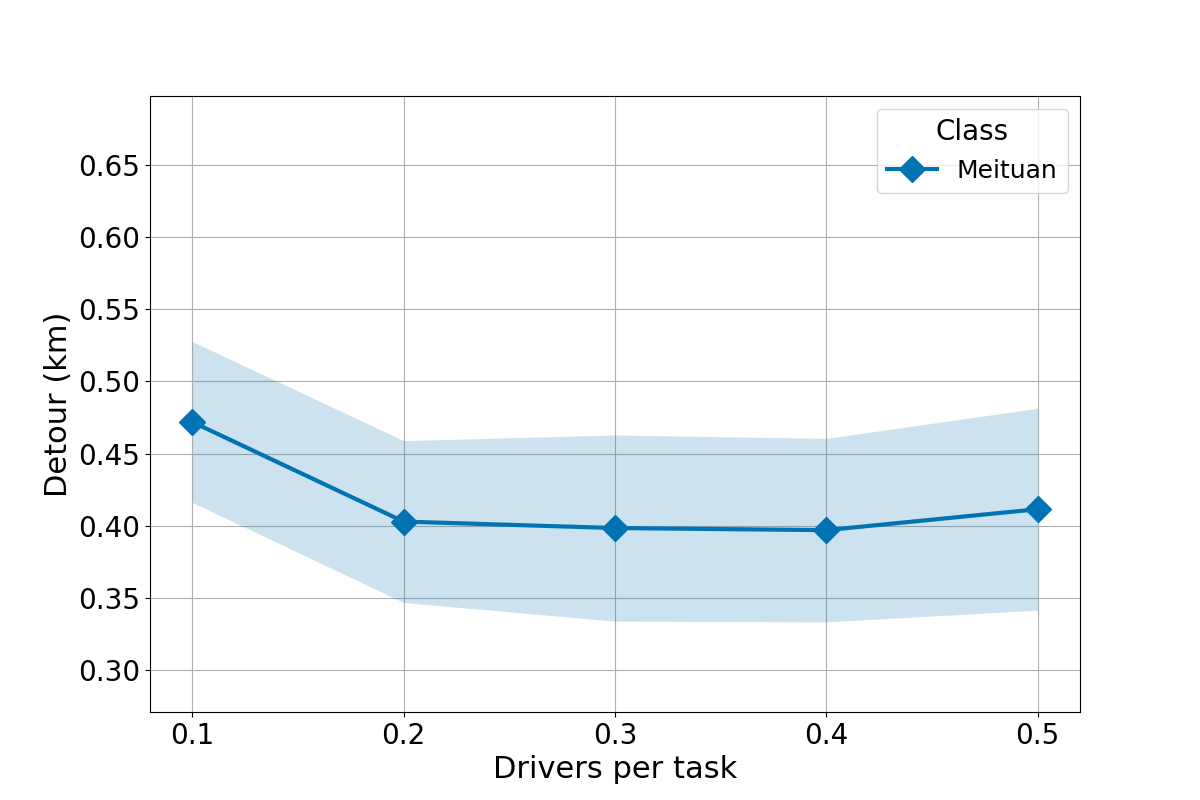}
		\caption{Detour.}
		\label{fig:meituan-driver-detour}
	\end{subfigure}
	\caption{Mean offer characteristics and their 95\% confidence intervals under the fitted Meituan specification by driver-to-task ratio.}
	\label{fig:meituan-driver-sensitivity}
\end{figure}

\cref{tab:meituan-computational-results,tab:meituan-heuristic-gaps} report the mean runtimes and optimality gaps, respectively. All 200 E-DDC runs terminate with a proven-optimal solution, and even the longest finishes in 1.89 seconds. H-C is the fastest method for every instance, while the additional optimization stages in H-DDC and E-DDC incur only small computational costs. H-DDC is also nearly exact, with a mean gap of 0.001\% over all instances and a maximum instance-level gap of 0.098\%. These small runtime and gap differences are consistent with a low compensation sensitivity that narrows the optimal solution space, as reflected in the low absolute levels of acceptance probability, detour, bundle size, and compensation reported above. H-C has a larger overall mean gap of 1.75\% and a maximum of 18.64\%, but this gap generally shrinks with task availability, from 3.54--4.86\% for 30-task instances to 0.43--0.69\% for 120-task instances, which is also observed in the main experiments. 

\begin{table}[H]
	\color{black}
	\centering
	\small
	\setlength{\tabcolsep}{7pt}
	\renewcommand{\arraystretch}{1.12}
	\caption{Mean runtimes in seconds under the fitted Meituan specification across varying instance sizes.}
	\label{tab:meituan-computational-results}
	\resizebox{\textwidth}{!}{%
	\begin{tabular}{c ccccc ccccc ccccc}
		\toprule
		\multirow{3}{*}{\textbf{Tasks}}
		& \multicolumn{5}{c}{\textbf{H-C}}
		& \multicolumn{5}{c}{\textbf{H-DDC}}
		& \multicolumn{5}{c}{\textbf{E-DDC}} \\
		\cmidrule(lr){2-6}\cmidrule(lr){7-11}\cmidrule(lr){12-16}
		& \multicolumn{5}{c}{\textbf{Driver-to-task ratio}}
		& \multicolumn{5}{c}{\textbf{Driver-to-task ratio}}
		& \multicolumn{5}{c}{\textbf{Driver-to-task ratio}} \\
		\cmidrule(lr){2-6}\cmidrule(lr){7-11}\cmidrule(lr){12-16}
		& \textbf{0.1} & \textbf{0.2} & \textbf{0.3} & \textbf{0.4} & \textbf{0.5}
		& \textbf{0.1} & \textbf{0.2} & \textbf{0.3} & \textbf{0.4} & \textbf{0.5}
		& \textbf{0.1} & \textbf{0.2} & \textbf{0.3} & \textbf{0.4} & \textbf{0.5} \\
		\midrule
		30
		& 0.02 & 0.02 & 0.05 & 0.19 & 0.02
		& 0.02 & 0.02 & 0.06 & 0.20 & 0.02
		& 0.02 & 0.02 & 0.06 & 0.20 & 0.03 \\
		60
		& 0.03 & 0.02 & 0.05 & 0.06 & 0.07
		& 0.04 & 0.04 & 0.07 & 0.10 & 0.10
		& 0.04 & 0.04 & 0.07 & 0.11 & 0.12 \\
		90
		& 0.03 & 0.07 & 0.13 & 0.12 & 0.11
		& 0.06 & 0.12 & 0.22 & 0.19 & 0.19
		& 0.07 & 0.14 & 0.24 & 0.21 & 0.21 \\
		120
		& 0.08 & 0.13 & 0.23 & 0.44 & 0.30
		& 0.17 & 0.23 & 0.40 & 0.67 & 0.45
		& 0.20 & 0.29 & 0.47 & 0.75 & 0.52 \\
		\bottomrule
	\end{tabular}
	}
\end{table}

\begin{table}[H]
	\color{black}
	\centering
	\small
	\setlength{\tabcolsep}{7pt}
	\renewcommand{\arraystretch}{1.12}
	\caption{Mean gaps to the proven-optimal E-DDC objective under the fitted Meituan specification (in \%) for different driver-to-task ratios.}
	\label{tab:meituan-heuristic-gaps}
	\begin{tabular}{c ccccc ccccc}
		\toprule
		\multirow{3}{*}{\textbf{Tasks}}
		& \multicolumn{5}{c}{\textbf{H-C}}
		& \multicolumn{5}{c}{\textbf{H-DDC}} \\
		\cmidrule(lr){2-6}\cmidrule(lr){7-11}
		& \multicolumn{5}{c}{\textbf{Driver-to-task ratio}}
		& \multicolumn{5}{c}{\textbf{Driver-to-task ratio}} \\
		\cmidrule(lr){2-6}\cmidrule(lr){7-11}
		& \textbf{0.1} & \textbf{0.2} & \textbf{0.3} & \textbf{0.4} & \textbf{0.5}
		& \textbf{0.1} & \textbf{0.2} & \textbf{0.3} & \textbf{0.4} & \textbf{0.5} \\
		\midrule
		30
		& 4.73 & 3.54 & 4.65 & 4.86 & 3.83
		& 0.00 & 0.00 & 0.00 & 0.00 & 0.00 \\
		60
		& 2.00 & 1.39 & 1.08 & 1.43 & 1.47
		& 0.00 & 0.00 & 0.00 & 0.00 & 0.00 \\
		90
		& 0.78 & 0.40 & 0.33 & 0.79 & 1.13
		& 0.01 & 0.00 & 0.00 & 0.00 & 0.00 \\
		120
		& 0.69 & 0.43 & 0.44 & 0.52 & 0.52
		& 0.01 & 0.01 & 0.00 & 0.00 & 0.00 \\
		\bottomrule
	\end{tabular}
\end{table}

\subsection{Perturbed coefficients} \label{sec:perturbed}

The observations of \cref{sec:results:sensitivity:behavioral} are based on the three stylized coefficient sets of \cref{tab:class_coefficients}. To verify that they do not depend only on these particular values, we perturb all four behavioral coefficients. Each coefficient ($\alpha$, $\beta_1$, $\beta_2$, $\gamma$) is tested at three levels, 80\%, 100\%, and 120\% of its value in \cref{tab:class_coefficients}. Combining all levels of all coefficients yields a full-factorial grid of $3^4=81$ configurations. The same multiplier is applied to a coefficient in all three classes, which preserves the decreasing effort and compensation sensitivity from Class~1 to Class~3, which defines the behavioral classes. Every configuration is solved on the subset of instances with 90 tasks and 45 drivers, including 50 mixed-class and 30 single-class instances (ten per class), and the class-level means underlying \cref{tab:mean_metrics} are recomputed. In \cref{tab:perturbation_means}, which mirrors the structure of \cref{tab:mean_metrics}, each cell reports the mean of the 81 configuration-level means together with their 95\% confidence interval describing their variability across the grid.

\begin{table}[H]
\caption{Mean metric levels and their 95\% confidence intervals (in brackets) across the 81 coefficient configurations.}
\label{tab:perturbation_means}
\centering
\resizebox{\textwidth}{!}{%
\blue{
\begin{tabular}{lcccccc}
	\toprule
	& \multicolumn{2}{c}{Class 1} & \multicolumn{2}{c}{Class 2} & \multicolumn{2}{c}{Class 3} \\ \cmidrule(lr){2-3} \cmidrule(lr){4-5} \cmidrule(lr){6-7}
	& Mixed-Class & Single-Class & Mixed-Class & Single-Class & Mixed-Class & Single-Class \\
	\midrule
	Acceptance Probability & \makecell[r]{0.95 \\ {[0.95, 0.96]}} & \makecell[r]{0.95 \\ {[0.95, 0.96]}} & \makecell[r]{0.92 \\ {[0.90, 0.93]}} & \makecell[r]{0.93 \\ {[0.92, 0.93]}} & \makecell[r]{0.69 \\ {[0.65, 0.74]}} & \makecell[r]{0.86 \\ {[0.84, 0.88]}} \\[1.4ex]
	Bundle Size & \makecell[r]{4.87 \\ {[4.86, 4.88]}} & \makecell[r]{4.70 \\ {[4.68, 4.71]}} & \makecell[r]{4.63 \\ {[4.61, 4.65]}} & \makecell[r]{4.71 \\ {[4.70, 4.73]}} & \makecell[r]{3.67 \\ {[3.58, 3.77]}} & \makecell[r]{4.68 \\ {[4.65, 4.72]}} \\[1.4ex]
	Compensation & \makecell[r]{13.62 \\ {[13.11, 14.13]}} & \makecell[r]{13.02 \\ {[12.51, 13.52]}} & \makecell[r]{14.19 \\ {[13.72, 14.67]}} & \makecell[r]{14.16 \\ {[13.64, 14.67]}} & \makecell[r]{12.96 \\ {[12.74, 13.18]}} & \makecell[r]{15.66 \\ {[15.20, 16.11]}} \\[1.4ex]
	Detour & \makecell[r]{1.91 \\ {[1.85, 1.97]}} & \makecell[r]{1.64 \\ {[1.60, 1.67]}} & \makecell[r]{1.82 \\ {[1.77, 1.87]}} & \makecell[r]{1.60 \\ {[1.56, 1.65]}} & \makecell[r]{1.75 \\ {[1.69, 1.81]}} & \makecell[r]{1.52 \\ {[1.46, 1.58]}} \\
	\bottomrule
\end{tabular}
}
}
\end{table}

\cref{tab:perturbation_means} reproduces every directional pattern of \cref{tab:mean_metrics}. Mean acceptance probability decreases from Class~1 to Class~3 in both pool types, with more pronounced decreases in the mixed pools. Bundle sizes in the single-class pools are nearly identical across classes (4.68--4.71, compared with 4.82--4.84 in \cref{tab:mean_metrics}), while in the mixed pools they decrease from Class~1 to Class~3. Relative to its single-class pool, Class~1 receives larger bundles, higher compensation, longer detours, and slightly higher acceptance probability in the mixed pools, whereas Class~3 receives smaller bundles, lower compensation, and markedly lower acceptance probability, which exactly mirrors the reallocation pattern in \cref{tab:mean_metrics}. Single-class compensation increases from Class~1 to Class~3 (13.02, 14.16, 15.66, compared with 13.84, 14.98, 16.57), and mean detour decreases from Class~1 to Class~3 within each pool type, while being higher in the mixed pools for every class, again as in \cref{tab:mean_metrics}. In conclusion, these agreements indicate that the content of \cref{tab:mean_metrics} is not tied to the particular coefficient values of \cref{tab:class_coefficients}. Thus, the conclusions made in \cref{sec:results:sensitivity:behavioral} are further supported by the perturbed experiment results.

}

\end{document}